\documentclass[11pt]{book}
\advance\footskip by 15pt
\usepackage{amsmath}
\usepackage{amssymb}
\usepackage{amsthm}
\usepackage{fancybox}

\newtheorem{theorem}{Theorem}[section]
\newtheorem{lemma}[theorem]{Lemma}
\newtheorem{corollary}[theorem]{Corollary}
\newtheorem{proposition}[theorem]{Proposition}

\newtheorem{remark}[theorem]{Remark}
\newtheorem{example}[theorem]{Example}

\newtheorem{definition}[theorem]{Definition}

\newenvironment{Proof}{\removelastskip\vskip12pt plus 1pt\noindent\em Proof.
\rm}{\hspace*{\fill}$\Box$\vskip12pt plus 1pt}
\newcommand{\m}{\mathfrak{m}}

\newcommand{\q}{\mathfrak{q}}
\newcommand{\codim}{\text{codim}}
\newcommand{\la}{\lambda}

\newcommand{\w}{\widetilde}
\newcommand{\M}{\mathbb{M}}

\def\depth{\operatorname{depth}}

\def\cocoa
{\mbox{\rm C\kern-.13em o\kern-.07
em C\kern-.13em o\kern-.15em A}}

\def\R{{\cal R}}

\def\deg{\operatorname{deg}}

\def\dim{\operatorname{dim}}
\def\depth{\operatorname{depth}}

\makeindex 

\makeatletter
\def\cleardoublepage{\clearpage\if@twoside%
\ifodd\c@page\else\thispagestyle{empty}\hbox{}\newpage\if@twocolumn\hbox{}\newpage\fi\fi\fi}
\makeatother

\newfam\msbfam

\newfont{\tuelvemsb}{msbm10 scaled\magstep1}

\font\tenmsb=msbm10
\font\sevenmsb=msbm7

\textfont\msbfam=\tuelvemsb
\scriptfont\msbfam=\tenmsb
\scriptscriptfont\msbfam=\sevenmsb

\def\dim{{\textsf{dim}\,}}

\def\codim{{\textsf{codim}\,}}

\def\depth{{\textsf{depth}\,}}

\def\deg{{\textsf{deg}\,}}

\def\min{{\textsf{min}\,}}

\def\deg{{\textsf{deg}\,}}

\title{  \bf \huge  Hilbert Functions  of Filtered  Modules   }

\author{\large   Maria Evelina Rossi \and 
\large Giuseppe Valla \
 }

\date{ }
 
\begin{document}

\maketitle

\par \vskip 4mm 

\newpage

\newpage

\setcounter{page}{3}
\tableofcontents

\addcontentsline{toc}{chapter}{\bf Introduction \hfill}

\chapter*{Introduction}

The notion of Hilbert function  is  central   in commutative algebra and is becoming increasingly important  in  algebraic geometry and in computational algebra.  In this presentation we shall deal  with  some aspects of the theory of    Hilbert functions of modules over local rings,  and we intend to guide  the reader along one of the  possible routes  through   the last three decades of progress in this  area of dynamic  mathematical activity.   

  Motivated by the ever increasing interest  in this field, our goal is to gather together many new developments of this theory  into one place, and to present them  using a unifying approach which gives  self-contained and easier  proofs.     In this text  we shall discuss many  results   by different  authors, following essentially the direction typified    by the pioneering  work of J. Sally (see \cite{  S1, S2, S, S6,  S7, S3, S4, S5}).    Our personal view of the subject is  most visibly expressed by the presentation of Chapters 1 and  2  in which we discuss the use of the superficial elements and related devices. 
  
  Basic techniques will be stressed with the aim of reproving   recent results by using a more elementary approach. This choice was made at the expense of certain  results and  various interesting aspects of the topic that, in this presentation, must remain peripherical.  We apologize to those whose work we may have failed to cite properly.  
\vskip 4mm

The material is intended for graduate students and   researchers who are interested in  results on  the Hilbert function and the Hilbert polynomial   of a local ring,   and applications of these.   The aim   was not to write a book on the subject, but rather to collect results and problems  inspired by    specialized lecture courses and schools recently delivered by the authors.  We hope the reader will   appreciate  the   large number of examples  and the   rich bibliography. 
\vskip 4mm

Starting from classical results of D. Northcott, S. Abhyankar, E. Matlis and  J. Sally, many papers have been written on this topic which   is  considered    an important  part of  the theory of blowing-up rings. 
This is because the Hilbert function of the local ring $(A,\m)$ is by definition the numerical function
$H_A(t):=\dim_k(\m^t/\m^{t+1}), $ hence it coincides with the classical Hilbert function of the standard graded algebra $gr_{\m}(A):=\oplus_{t\ge 0}\m^t/\m^{t+1}.$ The problems arise because, in  passing from $A$ to $gr_{\m}(A),$ we may  lose many good properties, such as being a complete intersection,  being Cohen-Macaulay or Gorenstein. 
 
  Despite the fact that the Hilbert function of a standard graded  algebra $A$ is well understood when  $A$ is Cohen-Macaulay, very little is known when it is a  local Cohen-Macaulay ring. 

 The Hilbert function of a local ring $(A,\m)$ is a classical invariant which gives information on the corresponding singularity. The reason is that 
the graded algebra $gr_{\m}(A)$ corresponds to an important geometric construction: 
namely, if $A$ is the localization at the origin of the coordinate ring of an affine variety $V$ passing through $0$, then $gr_{\m}(A)$ is the coordinate ring of the {tangent cone} of $V$, that  is the cone composed of all lines that are limiting positions of secant lines to $V$ in $0$. The {\it Proj} of this algebra can also be seen as the  {exceptional set} of the  {blowing-up} of $V$ in 0.

Other graded algebras come into the picture for  different reasons,  for example the Rees algebra, the Symmetric algebra, the Sally module  and the Fiber Cone.  All these algebras are doubly interesting  because on one side  they have a deep  geometrical  meaning, on the other side they are employed for detecting basic numerical characters  of the ideals in the local ring $(A,\m)$. Therefore, much attention has been paid  in the past to determining  under which circumstances these objects  have a good structure.

In some cases  the natural extension   of these results to $\m$-primary ideals has been achieved, starting from the fundamental  work of P. Samuel on multiplicities. More recently the generalization to the case of a descending multiplicative filtration of ideals of the local ring $A$ has now become of crucial importance. For example, the Ratliff-Rush filtration (cfr. papers   by   S. Huchaba, S. Itoh, T. Marley, T. Puthenpurakal, M.E. Rossi, J. Sally, G.~Valla) and the filtration given by the integral closure of the powers of an ideal (cfr. papers    by  A. Corso, S. Itoh, C. Huneke, C. Polini, B. Ulrich, W.~Vasconcelos) are fundamental tools in much of the recent work on blowing-up rings.

Even though of intrinsic interest, the extension to  modules  has   been  largely overlooked,  probably because, even in the classical case,   many problems were already so difficult. Nevertheless,  a number of results have been obtained in this direction: some of the work done by  D. Northcott, J.Fillmore, C. Rhodes, D. Kirby, H.  Meheran and,  more recently, T. Cortadellas and S.~Zarzuela,   A.V. Jayanthan and J.Verma,  T. Puthenpurakal has been  carried over to  the  general setting. 
\vskip 3mm
We remark that the graded algebra $gr_{\m}(A)$ can also be seen as the graded algebra associated to an ideal filtration of the ring itself, namely the $\m$-adic filtration $\{\m^j\}_{j\ge 0}$. This gives an indication of a possible natural extension of the theory to general filtrations of a finite module over the local ring  $(A,\m)$.

 Let $A$ be a commutative noetherian local ring with maximal ideal $\frak{m}$ and let $M$ be a finitely generated $A$-module.  Let $\q$ be an ideal of $A;$  a  $\frak{q}$-filtration ${\mathbb{M}} $ of $M  $   is a collection of submodules $M_j$ such that
$$ M=M_0\supseteq M_1 \supseteq \cdots \supseteq M_j \supseteq \cdots. $$
with the property that $\frak{q}M_j\subseteq M_{j+1}$ for each $j\ge 0. $
In the present work we consider only  \textbf{good $\frak{q}$-filtrations} of $M:   $  this means that  $M_{j+1}=\frak{q}M_j$ for all sufficiently large $j. $ A    good   $\frak{q}$-filtration  is also called a  \textbf{stable   $\frak{q}$-filtration.} 
For example, the $\frak{q}$-adic filtration on $M$ defined by $M_j:=  \q^jM $
is clearly a good  $\frak{q}$-filtration.  

We define the \textbf{associated graded ring} of $A$ with respect to $\q$ to be the graded ring $$gr_{\frak{q}}(A)=\bigoplus_{j\ge 0}(\frak{q}^j/\frak{q}^{j+1}).$$ Given a $\q$-filtration   $\mathbb{M}=\{M_j\}$ on the module $M,$ we consider the  \textbf{associated graded module} of $ M $ with respect to ${\mathbb{M}} $ $$gr_{\mathbb{M}}(M):=\bigoplus_{j\ge 0}(M_j/M_{j+1})$$ and for any $\overline{a}\in \q^n/\q^{n+1},$ $ \overline{m}\in M_j/M_{j+1}$ we define $\overline{a}\ \overline{m}:=\overline{am}\in M_{n+j}/M_{n+j+1}.$
The assumption that $\mathbb{M}$ is a $\q$-filtration ensures that this is well defined so that $gr_{\mathbb{M}}(M)$  has a natural structure as a graded module over the graded ring $gr_{\q}(A).$

Denote by $\la( *) $ the length of an $A$-module. If  $\la(M/ \q M$) is finite, then we can define the  {\bf Hilbert function} of the filtration $\mathbb{M}$, or 
of the filtered module $M$ with respect to the filtration $\mathbb{M}$. It is the numerical  function $$H_{\mathbb{M}}(j):=\la (M_j/M_{j+1}).$$ Its generating function is the power series $$P_{\mathbb{M}}(z):=\sum_{j\ge 0}H_{\mathbb{M}}(j)z^j.$$ which is called
the {\bf Hilbert series} of the filtration $\mathbb{M}.$  
 By the Hilbert-Serre theorem  we know that the series is of the form  $$ P_{\mathbb{M}}(z)=\frac{h_{\mathbb{M}}(z)}{(1-z)^r}$$
where  $h_{\mathbb{M}}(z)\in \mathbb{Z}[z],$ $h_{\mathbb{M}}(1)\not= 0$ and $r$ is the Krull dimension of $M.$ The polynomial $h_{\mathbb{M}}(z)$ is called the {\bf h-polynomial} of $\mathbb{M}.$

This implies that, for $n\gg 0 $
$$H_{\mathbb{M}}(n)=p_{\mathbb{M}}(n) $$  where the polynomial $p_{\mathbb{M}}(z)$ has rational coefficients, degree $r-1$ and is called 
the \textbf{Hilbert polynomial} of $\mathbb{M}.$

We can write $$p_{\mathbb{M}}(X):=\sum_{i=0}^{r-1}(-1)^ie_i(\mathbb{M})\binom{X+r-i-1}{r-i-1}$$ where we denote  for every integer $q \ge 0 $ $$ \binom{X+q}{q} := \frac{(X+q)(X+q-1) \dots (X+1)}{q!}.$$ The  coefficients  $e_i(\mathbb{M})$ are integers which will be called the \textbf{Hilbert coefficients}  of $\mathbb{M}$.

\vskip 3mm

It is clear that different Hilbert functions can have the same Hilbert polynomial;  but in many cases it happens that ``extremal" behavior  of some of the $e_i $   forces the filtration to have a specified Hilbert function. The trivial case is when the multiplicity is one: if this happens, then  $A$ is a regular local ring and $P_A(z)=\frac{1}{(1-z)^r}.$ 
Also the case of multiplicity 2 is easy, while the first non trivial result  along  this line was proved by J. Sally in \cite{S2}.
If $(A,\m)$ is a Cohen-Macaulay local ring of dimension $r$ and embedding dimension $v:=\dim_k(\m/\m^2),$ we let $h:=v-r$, the embedding codimension of $A.$     It is a result of Abhyankar that a  lower bound for $e_0$ is given by $$e_0\ge h+1.$$  This result   extends to the local Cohen-Macaulay rings the well known lower bound for the degree of a variety $ X $    in  ${\mathbb{P}^n}: $  
$$\deg X \ge \codim X +1. $$   The varieties for which the bound is attained are called  varieties of minimal degree and they are  completely classified.  In particular, they are always arithmetically Cohen-Macaulay.
In the local case, Sally proved that if the equality $e_0=h+1$ holds, then $gr_{\m}(A)$ is Cohen-Macaulay and $P_A(z)=\frac{1+hz}{(1-z)^r}.$  
\vskip 2mm

The next case, varieties satisfying $ \deg\  X= \codim\  X +2, $  is considerably more difficult. In particular such varieties are not necessarily arithmetically Cohen-Macaulay.  Analogously, in  the case $e_0=h+2,$ in \cite{S3} it was shown that $gr_{\m}(A)$ is not necessarily Cohen-Macaulay, the  exceptions lie among  the local rings of maximal Cohen-Macaulay type $\tau(A)=e_0-2.$ In the same paper Sally made the conjecture that, in the critical case, the depth of $gr_{\m}(A)$ is at least $r-1.$ This conjecture was  proved in \cite{W1} and \cite{RV1}  by using deep properties of the Ratliff-Rush filtration on the maximal ideal of $A$. Further in \cite{RV1}, all the possible Hilbert functions have been described:  they are of the form
$$ P_A(z)=\frac{1+hz+z^s}{(1-z)^r}$$ where $2\le s\le h+1.$ 

The next case, when $e_0=h+3, $ is  more complicated  and indeed still largely  open.  J. Sally,  in another  paper, see \cite{S6},   proved that if $A$ is Gorenstein and $e_0=h+3,$  then $gr_{\m}(A)$ is Cohen-Macaulay and $$P_A(z)=\frac{1+hz+z^2+z^3}{(1-z)^r}.$$
 If the Cohen-Macaulay type  $ \tau(A) $ is bigger than 1, then $ gr_{\m}(A)$ is no longer Cohen-Macaulay. Nevertheless, if $ \tau(A)<h, $    in \cite{RV3} the authors proved that $\depth(gr_{\m}(A))$ $\ge r-1$ and the Hilbert series is given by $$P_A(z)=\frac{1+hz+z^2+z^s}{(1-z)^r}$$ where $2\le s \le \tau(A)+2.$  This gives a new  and shorter  proof of the result of Sally, and it  points to the remaining  open question: what are the possible Hilbert functions  for a Cohen-Macaulay $r$-dimensional local ring with $e_0=h+3$ and $\tau(A)\ge h?$

It is clear that, moving  away from the minimal value of the multiplicity, things soon become   very difficult, and we do not  have any idea   what are the possible Hilbert functions of a one-dimensional Cohen-Macaulay local ring. 
The conjecture made by M.E. Rossi,  that this  function is non decreasing when $A$ is Gorenstein, is very much open, even for   coordinate rings  of monomial curves. 
 
Intuitively,    involving the higher Hilbert coefficients should  give stronger results. Indeed  if $(A,\m)$ is Cohen-Macaulay and we consider the $\m$-adic filtration $\mathbb{M}$ on $A$, then D. Northcott proved in \cite{N1} that  $e_1\ge e_0-1$ and,  if $e_1=e_0-1,$ then $P_{\mathbb{M}}(z)=\frac{1+hz}{(1-z)^r}$ while if $e_1=e_0$ then $P_{\mathbb{M}}(z)=\frac{1+hz+z^2}{(1-z)^r}$ (where $h$ is the embedding codimension of $A$). 

Results of this kind are quite remarkable because, in principle,  $e_0$ and $e_1$ give only partial information  on the Hilbert polynomial which depends  on the  asymptotic behavior of  the Hilbert function. 

In this presentation  we discuss further results   by several authors along  this line.

\vskip 4mm 
  
Over the past few years   several   papers have appeared which extend   classical results on the theory of Hilbert functions of Cohen-Macaulay local rings to the case of a filtration of a module. Very often, because of this  increased generality,   deep obstructions arise which can be overcome only by bringing new ideas to bear.   Instead, in this paper we illustrate how    a suitable and natural recasting     of the main basic tools of the classical theory is often enough to obtain  the required extensions.

More precisely what one needs is to make available  in the generalized  setting a few basic tools of the classical  theory, such as 
{\bf Superficial sequences, the Valabrega-Valla criterion, Sally's machine, Singh's formula}.  

Once these fundamental results have been established, the approach followed in this note gives a simple and  clean method  which applies  uniformly to many cases. 

In this way we make use  of the usual  machinery  to  get    easier proofs,   extensions of known results as well as numerous entirely new results. 
We mention   two nice examples of this philosophy:

$1. $ The problem of the existence of elements which are superficial simultaneously  with respect to a finite number of $\q$-filtrations on the same module has a natural solution in the module-theoretical approach (see Remark \ref{conti}), while it is rather complicated  in the ring-theoretical setting (see \cite{GN}, Lemma 2.3).

$2. $ In the literature two different definitions  of minimal multiplicity are given (see Section 2.1.). Here they are unified,  being just instances of the more general   concept of minimal multiplicity with respect to different filtrations of the same ring.

\vskip3mm The notion of {\bf superficial  element}  is a fundamental tool in our work. The original  definition was given by Zariski and Samuel, \cite{ZS}, pg.285. There it is shown how to use this concept for devising proofs by induction and reducing problems to lower dimensional ones.  We are concerned only with    the purely algebraic meaning of this notion, even if  superficial elements play an important  role also in  Singularity Theory, as shown by R. Bondil and Le Dung Trang in \cite {Bo} and \cite{BoTr}. 

We know that superficial sequences of order one always exist if the residue field is infinite, a condition which is not so restrictive. We make a lot of use of   this, often reducing a  problem to the one-dimensional case where things are  much easier.

The main consequence of this  strategy  is that  our arguments are quite elementary and, for example, we are able  to avoid  the more sophisticated  homological methods used in other papers.

The extension of the theory to the case of general filtrations on a module  has one more important  motivation. Namely, we have   interesting applications to the study of graded algebras which are not associated to a filtration. Here we have in mind  the  Symmetric algebra $S_A(\q)$, the Fiber cone  $F_\m (\q) $  and the Sally-module $ S_J(\q) $ of an ideal $\q$ of $A$ with respect to a minimal reduction $J.$ 
These  graded algebras  have been studied for their intrinsic interest; however,  since the rich theory of filtrations apparently does not apply, new and complicated  methods have been developed. 

We show here that  each of these algebras fits into certain  short  exact sequences  together with  algebras associated to filtrations. Hence we  can study the Hilbert function and the depth of these algebras with the aid of  the know-how we got  in the  case of a filtration.

This strategy   has been  already used  in \cite{HRV} to study the depth of the Symmetric Algebra $S_A(\m) $ of the maximal ideal $\m$ of a local ring $A. $ Also, in \cite{CZ, Co},   T. Cortadellas and S. Zarzuela used similar ideas to study  the Cohen-Macaulayness of the Fiber cone.  \vskip 3mm

In the last two chapters, we present   selected results  from the recent literature on  the Fiber cone  $F_\m (\q) $  and the Sally-module $ S_J(\q).$   We  have chosen  not  to pursue  here the study of the Symmetric algebra of the maximal ideal of a local ring, even if  we think it could be interesting and fruitful   to apply the ideas  of this paper also to that  problem.

In  developing this work,  one needs to consider filtrations on modules which are not necessarily Cohen-Macaulay. 
This opens up  a new and interesting terrain because most of the research  done on Hilbert functions  has been carried out  in the framework of Cohen-Macaulay local rings. Recently,  S. Goto, K. Nishida, A. Corso, W. Vasconcelos  and others have discovered  interesting results on the Hilbert function of general local rings. Following their methods, we tried to develop our theory as for  as possible but without the  strong assumption of Cohen-Macaulayness . But soon things became so difficult that we had to return  quickly to the classical assumption.

\vskip 3mm

We finish this introduction by giving a brief summary of each chapter. A longer description can be found  at the beginning of each chapter.

In the first chapter we introduce and discuss the notion of a good  $\q$-filtration of a module  over a local ring. The corresponding associated graded module is defined,  and a criterion for detecting regular sequences on this module is 
presented. Next we define  the notions  of Hilbert function and polynomial of a filtration  and  describe the  relationship with superficial elements.  Finally, we give a natural upper bound for the Hilbert function of a filtration   in terms of a maximal superficial sequence.

In the second chapter we give several upper bounds for the first two Hilbert coefficients of the Hilbert polynomial of a filtration;  it turns out that modules which are extremal with respect to these bounds have good associated graded modules and fully determined Hilbert functions. In particular,  we present a notable generalization of   Northcott's classical bound. We present several results for modules which do not require  Cohen-Macaulayness. Here   the theory of the $1$-dimensional case plays a crucial role. 

The third chapter deals with the third Hilbert coefficient;    some upper and lower bounds  are discussed. We introduce the Ratliff-Rush filtration associated to a $\q$-adic filtration  and show some  applications to the study of   border cases. The proofs become more sophisticated  because the complex structure  of   local rings of dimension at least two comes into play.  

In the fourth chapter  we give a proof of  Sally's conjecture in a very general context, thus greatly  extending the classical case. Several applications  to the first Hilbert coefficients are discussed. The main result  of this chapter is a bound on the reduction number of a filtration which has  unexpected  applications.

In chapters  five and six we explore the depth and the Hilbert coefficients of  the Fiber Cone and the Sally module respectively.   In spite of the fact  that   the Fiber Cone and the Sally module  are not    graded modules  associated  to a filtration,   the aim of this section is to show how one can  deduce their properties   as a consequence of the   theory on    filtrations.   In particular we will get short proofs of  several recent results   as an easy consequence  of certain classical results on the  associated graded rings to special filtrations.  

\vskip 4mm

Finally, we would like to take this  opportunity to thank sincerely Judith Sally because  her work  has had such a strong influence on our research into these subjects.  In particular several problems, techniques and  ideas presented in this text, came from a careful reading of her papers,  which are always rich in examples and motivating applications. 

Let us also not forget the many other colleagues who over the years have shared their ideas on these topics with us. Some of them were directly involved as co-authors in joint research  reported here, while others gave a substantial contribution via their publications and discussions.   

Since neither of us is a native English speaker,  we apologize for our numerous  linguistic infelicities.  We hope that this lack of expertise does not spoil the enjoyment   of the mathematical part.

\chapter{Preliminaries}
In this chapter  we present   the basic tools of the classical theory of filtered modules, in particular we introduce the  machinery we shall use throughout this  work: ${\mathbb{M}}$-superficial elements and their interplay with Hilbert Functions, the Valabrega-Valla  criterion, which is a basic tool for studying  the depth of  gr$_{\mathbb{M}}(M), $ ${\mathbb{M}}$-superficial sequences for an ideal $\q$ and their  relevance  to  Sally's machine, which is a very important device  to reduce dimension in questions relating to depth properties of blowing-up rings and  local rings.  

\section{Notation} Let $A$ be a commutative noetherian local ring with maximal ideal $\frak{m}$ and let $M$ be a finitely generated $A$-module. We will denote by $\la(\cdot )$ the {  length} of an $A$-module.    An (infinite) chain 
$$ M=M_0\supseteq M_1 \supseteq \cdots \supseteq M_j \supseteq \cdots $$ 
where the $M_n$ are submodules of $M$ is called a {\bf filtration} of $M,$ and denoted by ${\mathbb{M}} =\{M_n\}.$ Given an ideal  $\frak{q}$ in $A,$ ${\mathbb{M}} $  is a $\frak{q}$-{\bf filtration}  if  $\frak{q}M_j\subseteq M_{j+1}$ for all $j,$ and a
 {\bf good}   $\frak{q}$-filtration if $M_{j+1}=\frak{q}M_j$ for all sufficiently large  $j.$ Thus for example $\{\q^nM\}$ is a good $\q$-filtration.  In the literature a    good   $\mathfrak{q}$-filtration  is sometimes  called a  stable   $\frak{q}$-filtration.   We say that $ {\mathbb{M}}$ is nilpotent if $M_n =0 $ for $n \gg 0.$ Thus a  good ${\q}$-filtration ${\mathbb{M}} $ is nilpotent if and only if $\q \subseteq \sqrt{Ann\, M}.$ 
\vskip 2mm
\noindent From now on $ {\mathbb{M}}$ will denote always a good $\q$-filtration on the finitely generated $A$-module $M.$
\vskip 2mm
 \noindent We will assume  that the ideal $\frak{q}$ is proper. As a consequence 
$\cap _{i=0}^{\infty}M_i=\{0_M\}.$ 

\noindent Define $$gr_{\q}(A)=\bigoplus_{j\ge 0}(\frak{q}^j/\frak{q}^{j+1}).$$ This is a graded ring, in which the multiplication is defined as follows: if $a\in \q^i,$ $b\in \q^j$ define $\overline{a}\overline{b}$ to be $\overline{ab}$, i.e. the image of $ab$ in $\q^{i+j}/\q^{i+j+1}.$ This ring is the {\bf associated graded ring} of the ideal $\q.$

Similarly, if $M$ is an $A$-module and $\mathbb{M}=\{M_j\}$ is a $\q$-filtration of $M$, define

$$gr_{\mathbb{M}}(M)=\bigoplus_{j\ge 0}(M_j/M_{j+1})$$ which is a graded 
$gr_{\frak{q}}(A)$-module in a natural way. It is called the {\bf associated graded module} to the $\q$-filtration ${\mathbb{M}} =\{M_n\}.$ 

To avoid triviality we shall assume that $gr_{\mathbb{M}}(M)$ is not zero or equivalently $M \not = 0.$   
Each element $a \in A$ has a natural image, denoted by $a^* \in 
 gr_{\frak{q}}(A).$  If $a=0, $ then $a^*=0,$ otherwise $a^*=\overline{a}\in \q^t/\q^{t+1}$ where $t$ is the unique integer such that $a\in \q^t,$ $a\notin \q^{t+1}.$

If $N$ is a submodule of $M,$  by the Artin-Rees Lemma, the collection  $\{N\cap M_j\ | \ j\ge 0\}$ is a good $\q$-filtration of $N$.  Since  $$(N\cap M_j)/(N\cap M_{j+1})\simeq (N\cap M_j+M_{j+1})/M_{j+1}$$   $gr_{\mathbb{M}}(N)$ is a graded submodule of $gr_{\mathbb{M}}(M).$

On the other hand  it is clear that $\{(N+M_j)/N \ | \ j\ge 0\}$ is a good $\q$-filtration of $M/N$ which we denote  by $\mathbb{M}/N.$ These graded modules are related by the graded isomorphism $$gr_{\mathbb{M}/N}(M/N)\simeq gr_{\mathbb{M}}(M)/gr_{\mathbb{M}}(N).$$

If $a_1,\cdots,a_r$ are elements in $\q,  $   $ \not \in \q^2$ and $I=(a_1,\cdots,a_r)$ it is clear that $$\left [ (a_1^*,\cdots,a_r^*)gr_{\mathbb{M}}(M)\right ]_j=(IM_{j-1}+M_{j+1})/M_{j+1} $$ for each $j \ge 1.$
By the Artin-Rees Lemma one immediately gets 
that the following conditions are equivalent:
\vskip 2mm
1.  $gr_{\mathbb{M}/IM}(M/IM)\simeq gr_{\mathbb{M}}(M)/(a_1^*,\cdots,a_r^*)gr_{\mathbb{M}}(M).$

\vskip 2mm
2. $gr_{\mathbb{M}}(IM)=(a_1^*,\cdots,a_r^*)gr_{\mathbb{M}}(M).$

\vskip 2mm
3. $IM\cap M_j=IM_{j-1}\ \  \forall j\ge 1.$

\vskip 2mm
An interesting case in which  the above equalities  hold is  when the elements $a_1^*,\cdots,a_r^*$ form a regular sequence on $gr_{\mathbb{M}}(M).$ For example, if $r=1$ and $I=(a),$  then $\overline{a}\in \q/\q^2$ is regular on $gr_{\mathbb{M}}(M)$ 
if and only if the map  $$M_{j-1}/M_j {\overset {\overline a}{\to} } M_j/M_{j+1}$$ is injective for every $j\ge 1.$ This is equivalent to the equalities $M_{j-1}\cap(M_{j+1}:~a)=M_j$ for every $j\ge 1.$  An easy  computation shows the following result.

\begin{lemma} {\label{reg}} Let $a\in \q.$ The following conditions are equivalent:
\vskip 2mm
1. $\overline a\in \q/\q^2$ is a regular element on $gr_{\mathbb{M}}(M).$
\vskip 2mm
2. $M_{j-1}\cap(M_{j+1}:a)=M_j$ for every $j\ge 1.$
\vskip 2mm
3. $a$ is a regular element on $M$ and $aM\cap M_j=aM_{j-1}$ for every $j\ge 1.$
\vskip 2mm
4. $M_{j+1}:a=M_j$ for every $j\ge 0.$
\end{lemma}

This result leads us  to the  Valabrega-Valla criterion, a tool which has been very useful in the study of the depth of blowing-up rings (see \cite{VV}).

Many authors have discussed this topic recently.
For example  in \cite{HM2} Huckaba and Marley gave an extension of the classical result to the case of filtrations of ideals, by giving a completely new proof based on some deep investigation of a modified  Koszul complex.

Instead, in \cite{P1}, Puthenpurakal  extended the result to the case of $\q$-adic filtrations of a module by using the device of idealization of a module and then applying the classical result. 

This would suggest  that the original proof does not work in the more general setting. But, after looking at it carefully, we can say that, in order to prove the following very general statement, one does not need any new idea:  a straightforward adaptation of the  dear old  proof does the job.

\begin{theorem} \label{vv} {\rm{(Valabrega-Valla)}}  Let $a_1,\cdots,a_r$ be elements in $\q, \not \in \q^2, $ and $I$ the ideal they generate. Then $a_1^*,\cdots,a_r^*$ form a regular sequence on $gr_{\mathbb{M}}(M)$ if and only if $a_1,\cdots,a_r$ form a regular sequence on $M$ and $IM\cap M_j=IM_{j-1}\ \  \forall j\ge 1.$
\end{theorem}

\section{Superficial elements}  A  fundamental tool in local algebra is the notion of superficial element. This notion goes back to P. Samuel and our methods are also related to the construction given by Zariski and Samuel (\cite{ZS} p.296).

\begin{definition} An element $a\in \q,$ is called \textbf{$\mathbb{M}$-superficial}
for $\q$ if  there exists a non-negative integer $c$ such that $$(M_{n+1}:_Ma)\cap M_c=M_n$$ for $n\ge c.$
\end{definition}

\noindent 

For every $a \in \q$ and $n \ge c, $  $M_n$ is contained in $(M_{n+1}:_Ma)\cap M_c.$ Then it is the other inclusion that makes superficial elements special. 
It is clear that if $\mathbb{M}$ is nilpotent, then every element is superficial.  If the length $\lambda(M/\q M) $ is finite, then $\mathbb{M}$ is nilpotent if and only if $\dim M =0.$ Hence in the following,    when we deal with superficial elements, {\bf{we shall assume that }} $\dim M \ge 1.$

If this is the case, as a consequence of the definition, we deduce that  {$\mathbb{M}$-superficial} elements $ a $ 
for $\q$  have order one, that is $a \in \q \setminus \q^2. $  With a slight modification of the given definition, superficial elements can be introduced of every order, but  in the following we need superficial element of order one because they have a better  behaviour in  studying Hilbert functions.

Hence  { {superficial element }} always  means a superficial element of order one. 
 It is well known that superficial elements do not always exist, but their  existence is guaranteed if the residue field is infinite (see Proposition 8.5.7. \cite{HSw}).  By passing, if needed, to the  faithfully  flat extension   $A[x]_{\m A[x]} $ ($x$ is a variable over $A$) we may assume that the residue field is infinite. 
  
If $\q $ contains a regular  element on $M, $ it  is easy to see that every $\mathbb{M}$-superficial element of $\q$  is regular  on $M.$   
 
 \vskip 3mm
 Given  $A$-modules  $M $ and $N$, let ${\mathbb{M}} $ and ${\mathbb{N}}$ be  good $\frak{q}$-filtrations  of $M$ and $N $ respectively.  We define the new filtration  as follows  $$  \mathbb{M}\oplus \mathbb{N}: \ \ \  M \oplus N \supseteq M_1 \oplus N_1 \supseteq \dots \supseteq M_n \oplus N_n \supseteq \dots $$ 
It is easy to see that  $\mathbb{M}\oplus \mathbb{N} $ is a good $\q$-filtration on the $A$-module  $ M \oplus N.$ Of course this construction can be extended to any finite number of modules.

The following  remark is due to David Conti in his thesis (see \cite{Conti}).

\begin{remark}\label{conti}  {\rm{ Let $ \mathbb{M}_1, \dots, \mathbb{M}_n $ be $\q$-filtrations of $M$ and let $a \in \q.$ Then $a$ is $ \mathbb{M}_1\oplus  \dots \oplus \mathbb{M}_n$-superficial  for $\q$ if and only if $a$ is $\mathbb{M}_i$-superficial  for $\q$  for every $i=1, \dots, n.$

  }}
\end{remark}

\vskip 2mm \par \noindent This result   is an easy consequence of the good behaviour of    intersection  and    colon of ideals  with respect to  direct sum  of modules. As a consequence we deduce that, if the residue field is infinite,  we can always find an element $a \in \q$ which is superficial for a finite number of $\q$-filtrations  on $M.$    As  mentioned in the introduction,  we want to apply the general theory of  the filtrations on a module  to the study of certain  blowing-up rings which are not necessarily associated graded rings  to a single filtration. Since they are related   to different filtrations,  the above remark will be relevant  in our approach.   

David Conti also remarked   that    $\mathbb{M}$-superficial elements of order  $s \ge 1$ for $\q$  can be seen as superficial 
elements of order  one for a suitable filtration which is strictly related to $\mathbb{M}$.  Let   $\mathbb{N}$ be the $\q^s$-filtration :
{\small { $$ M \oplus M_1 \oplus \dots \oplus M_{s-1} \supseteq  M_s\oplus M_{s+1} \oplus \dots \oplus M_{2s-1}  \supseteq  M_{2s }\oplus M_{2s +1} \oplus \dots \oplus M_{3s-1} \supseteq \dots $$}}
Then it is easy to see that an element $a$ is $\mathbb{M}$-superficial of degree $ s $ for $\q$ if and only if $a$ is $\mathbb{N}$-superficial of degree one  for $\q^s.$ This remark could be useful in studying properties which well behave  with  the direct sum.

\vskip 3mm

  We  give now  equivalent conditions for an element to be $\mathbb{M}$-superficial for $\q.$ Our development of the theory of superficial elements is basically the same as that given by Kirby in \cite{K} for the case $M=A$ and $\{M_j\}=\{\q^j\} $ the $\q$-adic filtration on $A.$ 

\noindent If there is no confusion, we let $$G:=gr_{\mathbb{M}}(M), \ \ \ Q:=\bigoplus_{j\ge 1}(\frak{q}^j/\frak{q}^{j+1}).$$

Let $$\{0_M\}=P_1\cap\cdots\cap P_r\cap P_{r+1}\cap\cdots\cap P_s$$ be an irredundant primary decomposition of $\{0_M\}$ and let $\{\wp_i\}=Ass(M/P_i).$ Then $Ass(M)=\{\wp_1,\cdots, \wp_s\}$ and $\wp_i=\sqrt{0:(M/P_i)}.$

\noindent We may  assume   $\q\nsubseteq \wp_i$ for $i=1,\cdots,r$ and 
$\q\subseteq \wp_i$ for $i=r+1,\cdots,s.$ Then we let $$N:=P_1\cap\cdots\cap P_r.$$

\noindent It is clear that  $$N=\{x\in M\ | \ \exists \ n,  \q^n x =0_M\}$$ and that $N\cap M_j=\{0\}$ for all large $j$ and $Ass(M/N)=\bigcup _{i=1}^r\wp_i.$
\vskip4mm 
Similarly  we denote by  $H$ the   homogeneous submodule of $G$ consisting of the elements $\alpha\in G$ such that $ Q^n\alpha =0_G$, hence $$H=\{\alpha \in G \ | \ \exists \ n,  Q^n\alpha =0_G\}.$$ 

\noindent If $$\{0_G\}=T_1\cap\cdots\cap T_m\cap T_{m+1}\cap\cdots\cap T_l$$ is an irredundant primary decomposition of $\{0_G\}$ and we let $\{\frak{P}_i\}=Ass(G/T_i),$ then $Ass(G)=\{\frak{P}_1,\cdots, \frak{P}_l\}$ and $\frak{P}_i=\sqrt{0:(G/T_i)}.$ 

Further, if we assume that $Q\nsubseteq \frak{P}_i$ for $i=1,\cdots,m$ and 
$Q\subseteq \frak{P}_i$ for $i=m+1,\cdots,l,$ then $$H=T_1\cap\cdots\cap T_m$$ and $Ass(G/H)=\bigcup _{i=1}^m\frak{P}_i.$

\begin{theorem} \label{super} Let $a\in \q \setminus \q^2, $ the following conditions are equivalent:
\vskip 2mm
1. $a$ is $\mathbb{M}$-superficial
for $\q.$
\vskip 2mm
2.  $a^*\notin \bigcup_{i=1}^m \frak{P_i}$
\vskip 2mm
3. $H:_Ga^*=H.$  
\vskip 2mm
4. $N : a=N$ and $M_{j+1}\cap  aM=aM_j$ for all large $j.$
\vskip 2mm
5. $(0:_Ga^*)_j=0$ for all large  $j.$
\vskip 2mm
6. $M_{j+1}:a=M_j+(0:_Ma)$ and $M_j\cap (0:_Ma)=0$ for all large $j.$

\end{theorem}

We note with  \cite{ZS}  that if the residue field $A/\m$ is infinite, condition 2. of the above theorem insures  the existence of $\mathbb{M}$-superficial elements for $\q $ as we have already mentioned. Moreover condition 5. says that $a$ is $\mathbb{M}$-superficial
for $\q $ if and only if $a^* $ is an homogeneous  {\it filter-regular  } in $G. $  We may refer to \cite{T} concerning the definition  and the properties of the homogeneous filter-regular elements. Filter-regular elements were also  introduced   in the local contest by N.T. Cuong, P. Schenzel and N.V. Trung in \cite{CST}.  One of the main result in \cite{CST}  says  that a local ring is
generalized Cohen-Macaulay  if and only if 
every system of parameters is filter-regular. The notion of filter-regular element in a local ring is weaker than superficial element and in general it does not behave well with Hilbert  functions.  A superficial element is filter-regular, but the converse does not hold. It is enough to recall that if $a$ is $M$-regular,  then accordingly with \cite{CST},   $a $ is a filter-regular element, but it is not necessarily a superficial element. 
\vskip 2mm

A sequence of elements $a_1,\cdots,a_r$  will be called a 
\textbf{
$\mathbb{M}$-superficial sequence
for $\q$
} if, for $i=1,\cdots,r,$  $a_i$ is an $\left (\mathbb{M}/(a_1,\cdots,a_{i-1})M\right )$-superficial element  for $\q.$

In order to prove properties of superficial sequences often we can argue by induction on the number of elements, the above theorem giving the first step of the induction.  For example, since $\depth_\q(M)\ge 1$ implies $N=0$, condition 4. gives the following result. Here  $\depth_\q(M) $ denotes   the common cardinality of all the maximal $M$-regular sequences of elements in $\q.$

\begin{lemma}\label{m} Let $a_1,\cdots,a_r$  be    an
$\mathbb{M}$-superficial sequence for $\q.$ Then   $a_1,\cdots,a_r$  is   a regular 
 sequence on $M$  if and only if  $\depth_\q(M)\ge r$.
\end{lemma}

In the same way, since $\depth_Q (gr_{\mathbb{M}}(M))\ge 1$ implies $H=0$, condition 3. and Theorem
 \ref{vv}  give the following result which shows the relevance of superficial elements in the study of the depth of blowing-up rings. 

\begin{lemma}\label{sup} Let $a_1,\cdots,a_r$  be    an
$\mathbb{M}$-superficial sequence for $\q.$ Then  $a_1^*,\cdots,a_r^*$ is a regular sequence on $ gr_{\mathbb{M}}(M)$
if and only if $\depth_Q (gr_{\mathbb{M}}(M))\ge r.$\end{lemma}

Now we come to  \textbf{Sally's machine} or {\it{Sally's descent}}, a very important device  to reduce dimension in questions relating to depth properties of blowing-up rings.
\begin{lemma}$\rm{(Sally's  \ machine)}$  Let $a_1,\cdots,a_r$  be    an
$\mathbb{M}$-superficial sequence for $\q$ and $I$ the ideal they generate. Then  $\depth_Q (gr_{\mathbb{M}/IM}(M/IM))\ge 1,$ if and only if $\depth_Q (gr_{\mathbb{M}}(M))\ge r+1.$
\end{lemma}

A proof of the  {\it{if}} part   can be obtained by a straightforward adaptation of the original proof given by Huckaba and Marley  in \cite{HM2}, Lemma 2.2. The converse is an easy consequence of Lemma \ref{m} and Theorem \ref{vv}.

\vskip 2mm
 
\vskip 3mm We will also need a property of superficial elements which seems to be neglected in the literature. It is well known  that if $a$ is $M$-regular, then $M/aM$ is Cohen-Macaulay if and only if $M$ is Cohen-Macaulay. We prove,  as a consequence of  the following Lemma, that the same holds for a superficial element,  if the dimension of the module $M$ is at least two.  

In the following we denote by $H_{\q}^i(M)$ the $i$-th local cohomology module  of $M$ with respect to $\q.$ We know that $H_{\q}^0(M):=\cup_{j\ge 0}(0:_M\q^j)=0:_M\q^t$ for every $t\gg 0; $ further  $\min\{i \ | \ H^i_{\q}(M)\neq 0\}=   \depth_\q(M). $  

\begin{lemma} \label{nuovo} Let $ a$ 
 be an $\mathbb{M}$-superficial element for $\q$ and let $j\ge 1$.Then    $\depth_\q(M)  \ge j+1$ if and only if  $\depth_\q(M/aM)\ge j.$
\end{lemma}
\begin{proof} Let     $ \depth_\q(M) \ge j+1$; then  $ \depth  M>0$ so that  $a$ is $M$-regular. This implies   $ \depth_\q(M/aM)=  \depth_\q(M) -1\ge j+1-1=j.$

Let us assume now that    $\depth_\q(M/aM) \ge j.$ Since $j\ge 1,$ this implies  $H^0_{\q}(M/aM)=0.$ Hence $H^0_{\q}(aM)=H^0_{\q}(M),$ so that $H^0_{\q}(M)\subseteq aM.$ We claim that $H^0_{\q}(M)=aH^0_{\q}(M).$ If this is the case, then,  by Nakayama, we get $H^0_{\q}(M)=0$ which  implies    $\depth (M)  >0,$ so that   $a$ is $M$-regular.  Hence 

$$   \depth_\q(M)  =   \depth_\q(M/aM)+1\ge j+1,$$ 
as wanted. 

Let us prove the claim. Suppose by contradiction that $$aH^0_{\q}(M)\subsetneq H^0_{\q}(M)\subseteq aM,$$ and let $ax\in H^0_{\q}(M), x\in M\setminus H^0_{\q}(M).$ This means that for every  $t\gg 0$ we have $$\begin{cases} a\q^tx=0\\ \q^tx\neq 0 \end{cases} $$ We prove  that this implies that $a$ is not
 $\mathbb{M}$-superficial  for $\q.$  Namely, given a positive integer $c,$ we can find an integer   $t\ge c$  and an element $d\in \q^t$ such that $adx=0$ and $dx\neq 0.$ Since $ \cap M_i=\{0\},$ we have $dx\in M_{j-1}\setminus M_j$ for some integer $j.$ Now, $d\in \q^t$ hence $dx\in M_t\subseteq M_c,$ which implies  $j\ge c.$  Finally we have $dx\in M_c,$ $dx\notin M_j$ and $adx=0 \in M_{j+1},$ hence $$(M_{j+1}:_Ma)\cap M_c\supsetneq M_j.$$ The claim and the Lemma are proved.\end{proof}

It is interesting to recall that the integral closure of the ideals generally behaves well  going modulo a superficial element.
For example S. Itoh \cite{I}, pag.648, proved that :

\begin{proposition} \label{Itoh}
If $\q$ is an $\m-$primary ideal which is  integrally closed 
in a Cohen-Macaulay local ring $(A,\m)$ of dimension $r \ge 2$,
then (at least after passing to a faithfully flat extension) 
there exists a superficial element $a \in \q $
such that $\q A/(a)$ is integrally closed in $A/(a)$.
\end{proposition}

This result will be useful in proofs working by induction on $r.$  The
compatibility of integrally closed  ideals with specialization by generic elements can be extended to that
of modules (see \cite{HU}).  
\vskip 2mm

However we will see later, superficial elements do not behave well for Ratliff-Rush 
closed ideals:
there exist many ideals of $A $ all of whose powers are Ratliff-Rush closed,
yet for every superficial element $a \in \q$,
$\q /(a)$ is not Ratliff-Rush closed.  Examples are given by   Rossi and Swanson in  \cite{RS}. Recently  
Puthenpurakal in \cite{Puth}  characterized  local rings  and ideals for which the Ratliff-Rush filtration {\it{behaves well}}\  modulo a superficial element. 
\par 
It is interesting to recall that    Trung  and  Verma in \cite{TV} introduced superficial sequences  with respect to a set of ideals   by  working  in the multigraded contest.

\vskip 6mm 
\section{The Hilbert Function and Hilbert coefficients} 
Let  $\mathbb{M}$ be a good $\q$-filtration on the $A$-module $M.$ 
From now on we shall require the assumption that the length of $M/ \q M$, which we denote by $\la(M/ \q M$), is finite.   In this case there exists an integer $s$ such that $\m^sM \subseteq (\q+(0:_AM))M, $ hence (see [N2]), the ideal $\q+(0:_AM)$ is primary for the maximal ideal $\m$. Also the length of $M/M_j$ is finite for all $j\ge 0.$
\vskip 2mm
\noindent From now on  $\q$ will denote an $\m$-primary ideal of the local ring $(A, \m).$
\vskip 2mm

In this setting we can define the  \textbf{\textbf{Hilbert function}} of the filtration $\mathbb{M}$ or simply 
of the filtered module $M$ if there is no confusion. By definition it is the function $$H_{\mathbb{M}}(j):=\la (M_j/M_{j+1}).$$

\noindent It is also useful to consider the numerical function $$H^1_{\mathbb{M}}(j):=\la (M/M_{j+1})=\sum_{i=0}^jH_{\mathbb{M}}(i)$$ which is called the  \textbf{Hilbert-Samuel function} of the filtration $\mathbb{M}$ or of the filtered module $M.$

The \textbf{Hilbert series}  of the filtration $\mathbb{M}$  is the power series $$P_{\mathbb{M}}(z):=\sum_{j\ge 0}H_{\mathbb{M}}(j)z^j.$$

The power series $$P^1_{\mathbb{M}}(z):=\sum_{j\ge 0}H^1_{\mathbb{M}}(j)z^j$$ is called the \textbf{Hilbert-Samuel series}
of  $\mathbb{M}$. It is clear that $$P_{\mathbb{M}}(z)=(1-z)P^1_{\mathbb{M}}(z).$$

By the Hilbert-Serre theorem, see for example \cite{BH},  we can write $$ P_{\mathbb{M}}(z)=\frac{h_{\mathbb{M}}(z)}{(1-z)^r}$$
where  $h_{\mathbb{M}}(z)\in \mathbb{Z}[z],$ $h_{\mathbb{M}}(1)\not= 0$ and $r$ is the dimension of $M.$

The polynomial $h_{\mathbb{M}}(z)$  is called 
the \textbf{h-polynomial} of $\mathbb{M}$ and we clearly have $$P^1_{\mathbb{M}}(z)=\frac{h_{\mathbb{M}}(z)}{(1-z)^{r+1}}.$$

An easy computation shows that if we let for every $i\ge 0$ $$e_i(\mathbb{M}):=\frac
{h_{\mathbb{M}}^{(i)}(1)}{i!} $$ then for $n\gg 0$ we have $$H_{\mathbb{M}}(n)=\sum_{i=0}^{r-1}(-1)^ie_i(\mathbb{M})\binom{n+r-i-1}{r-i-1}.$$ The polynomial
$$p_{\mathbb{M}}(X):=\sum_{i=0}^{r-1}(-1)^ie_i(\mathbb{M})\binom{X+r-i-1}{r-i-1}$$ has rational coefficients and is  called the \textbf{Hilbert polynomial} of $\mathbb{M};$ it verifies the equality
$$H_{\mathbb{M}}(n)=p_{\mathbb{M}}(n)$$ for $n\gg 0.$
  The integers $e_i(\mathbb{M})$ are called the \textbf{Hilbert coefficients}  of $\mathbb{M}$ or of the filtered module $M$ with respect to the filtration $\mathbb{M}.$ 

In particular $e_0(\mathbb{M}) $ is the  \textbf{multiplicity} of $\mathbb{M}$ and, by Proposition 11.4. in \cite{AM} (also Proposition  4.6.5. in \cite{BH}) 
$e_0(\mathbb{M}) =e_0(\mathbb{N}),  $ for every couple of good $\q$-filtrations $\mathbb{M}$ and $\mathbb{N} $ of $M.$

\vskip 3mm
\noindent Since $P^1_{\mathbb{M}}(z)=\frac{h_{\mathbb{M}}(z)}{(1-z)^{r+1}} $
the polynomial 
$$p^1_{\mathbb{M}}(X):=\sum_{i=0}^{r}(-1)^ie_i(\mathbb{M})\binom{X+r-i}{r-i}$$  verifies the equality
$$H^1_{\mathbb{M}}(n)=p^1_{\mathbb{M}}(n)$$ for $n\gg 0$
 and is  called the \textbf{Hilbert-Samuel polynomial} of $M.$ 
 
 \vskip 2mm In the following we will write the $h$-polynomial of $M$ in the form $$h_{\mathbb{M}}(z)=h_0(\mathbb{M})+h_1(\mathbb{M})z+\cdots +h_s(\mathbb{M})z^s,$$  so that the integers 
 $h_i(\mathbb{M})$ are well defined for every $i\ge 0$ and  we have 
\begin{equation} \label{coeff} e_i(\mathbb{M})=\sum_{k\ge i} {{k}\choose{i}}h_k(\mathbb{M}).
\end{equation}

\noindent   Finally we remark that if $M$ is Artinian, then $e_0(\mathbb{M})=\la (M).$
  
  \noindent In the case of the $\q$-adic filtration on the ring $A, $ we will denote  by $H_{\q}(j)$ the Hilbert function,  by $ P_{\q}(z) $ the Hilbert series, and by $e_i(\q) $ the Hilbert coefficients of the $\q$-adic filtration.  In the case of  modules, that is the $\q$-adic filtration on a module  $M, $ we will replace $\q $ with $\q M.$

\vskip 3mm
The following result  was  called \textbf{Singh's formula} because the corresponding equality in the classical case was obtained by Singh in \cite{Si}.  See also  \cite{ZS}, Lemma 3, Chapter 8, for   the corresponding equality with  Hilbert polynomials.

\begin{lemma} \label{Singh}  Let $a\in \q;$ then for every $j\ge 0$ we have $$H_{\mathbb{M}}(j)=H^1_{\mathbb{M}/aM}(j)-\la(M_{j+1}:a/M_j).$$
\end{lemma}
\begin{proof} The proof is an easy consequence of the following exact sequence:
$$0\to (M_{j+1}:a)/M_j \to M/M_j\overset{a}  {\to} M/M_{j+1}\to M/(M_{j+1}+aM) \to 0.$$
\end{proof}

We remark that, by Singh's formula,   for every $a\in \q $ we have 
\begin{equation} \label{s}
   P_{\mathbb{M}}(z) \le  P^1_{\mathbb{M}/aM}(z) \end{equation}

It is thus interesting to consider elements $a\in \q$ such that $ P_{\mathbb{M}}(z)=  P^1_{\mathbb{M}/aM}(z).$ By the above formula and Lemma \ref{reg}, these are exactly the elements $a\in \q$ such that $\overline a \in \q/\q^2$ is  regular over $gr_{\mathbb{M}}(M).$
\vskip 2mm
As a corollary of Singh's formula we get   useful properties of superficial elements.

\begin{proposition}\label{ei}  Let $a$ be an $\mathbb{M}$-superficial element  for $\q$ and let $ r=\dim M \ge 1.$ Then we have:
\vskip 2mm 
1. $\dim(M/aM)=r-1.$
\vskip 2mm
2. $e_j(\mathbb{M})=e_j(\mathbb{M}/aM)$ for every $j=0,\cdots,r-2.$
\vskip 2mm 
3. $e_{r-1}(\mathbb{M}/aM)=e_{r-1}(\mathbb{M})+(-1)^{r-1}\la(0:a).$
\vskip 2mm 
4. There exists an integer $j$ such that for every $n\ge j-1$ we have
\vskip 2mm 
$e_r(\mathbb{M}/aM)=e_r(\mathbb{M})+(-1)^r\left(\sum_{i=0}^n\la(M_{i+1}:a/M_i)-(n+1)\la(0:a)\right )$  
\vskip 2mm 
5.  $  a^*$ is a regular element on $gr_{\mathbb{M}}(M)$ if and only if $ P_{\mathbb{M}}(z)=  P^1_{\mathbb{M}/aM}(z)=\frac{P_{\mathbb{M}/aM}(z)}{1-z} $  if and only if $a$ is $M$-regular and  $e_r(\mathbb{M})=e_r(\mathbb{M}/aM)$  
\end{proposition}
\begin{proof}  By Lemma \ref{Singh}  we have $$P_{\mathbb{M}}(z)=P^1_{\mathbb{M}/aM}(z)-\sum_{i\ge 0}\la(M_{i+1}:a/M_i)z^i.$$
Since $a$ is superficial, by Theorem \ref{super}, 6, there exists an integer $j$ such that for every $n\ge j$ we have  $$M_{n+1}:a/M_n=M_n+(0:a)/M_n\simeq (0:a)/(0:a)\cap M_n=(0:a).$$

Hence we can write $$P_{\mathbb{M}}(z)=P^1_{\mathbb{M}/aM}(z)-\sum_{i= 0}^{j-1}\la(M_{i+1}:a/M_i)z^i-\frac{\la(0:a)z^j}{(1-z)}.$$
This proves $1.$  (actually $1.$ follows also by  Theorem \ref{super} $(4.)$) and,  as a consequence,  we get
$$h_{\mathbb{M}}(z)=h_{\mathbb{M}/aM}(z)-(1-z)^r\left (\sum_{i= 0}^{j-1}\la(M_{i+1}:a/M_i)z^i\right )-(1-z)^{r-1}\la(0:a)z^j.$$ This gives easily  2., 3. and 4, while 5. follows from  4. and  Lemma \ref{reg}.
\end{proof}

As we can see in the proof of the above result, if $a$ is a  $\mathbb{M}$-superficial element,  there exists an integer $j$ such that $M_{n+1}:a/M_n = (0:a) $ for every $n \ge j.$ Hence,    Proposition \ref{ei}  (4.)can be rewritten as follows  $$e_r(\mathbb{M}/aM)=e_r(\mathbb{M})+(-1)^r\left(\sum_{i=0}^{j-1}\la(M_{i+1}:a/M_i)-j \la(0:a)\right ).$$
 
For example let $A=k[[X,Y,Z]]/(X^3,X^2Y^3,X^2Z^4)=k[[x,y,z]].$ We have  $r= \dim(A)=2$ and if we consider on $M =A$ the $\m$-adic filtration $\mathbb{M}=\{\m^j\}, $    then   $y $ is an $\mathbb{M}$-superficial element for $\q=\m.$ We have 
$$P_{\mathbb{M}}(z)= \frac{1+z+z^2-z^5-z^6+z^9}{(1-z)^2}$$ so that 
$$e_0(\mathbb{M})=2, \ \ e_1(\mathbb{M})=1,  \ \ e_2(\mathbb{M})=12.$$
Also it is clear that 

$$\la(\m^{n+1}:y/\m^n)=\begin{cases}
  0    & \text{for} \ \ n=0,\cdots,4, \\
  1    & \text{for} \ \ n=5,\\
  2 & \text{for} \ \ n=6,\\
  3 & \text{for} \ \ n=7 ,\\
  4=\la (0:y) & \text{for} \ \ n\ge 8
\end{cases}$$ so that $j=8$ in the above proposition. Hence, by using Proposition \ref{ei},  
 $$e_0(\mathbb{M}/(y))=2, \ \ e_1(\mathbb{M}/(y))=-3,  \ \ e_2(\mathbb{M}/(y))=-14$$ 

Notice that   $$P_{\mathbb{M}/(y)}(z)= \frac{1+z+z^2-z^6}{1-z}. $$
 
  \vskip 3mm
 \begin{remark} \label{criterion} {\rm {By  Proposition \ref{ei}, Theorem \ref{super} and Singh's formula, if      $a \in \q $ is  an element which is  $M$-regular,    then  $$a { \mbox { is }}   \mathbb{M}{\mbox {-superficial  for }} \q   \iff  e_i(\mathbb{M})= e_i(\mathbb{M}/aM) {\mbox { for every }} i=0, \dots, r-1.$$
\noindent   It is an useful criterion in order to check if an element is superficial or not. From  the computational point of view  it reduces  infinitely many conditions,  coming a priori  from the definition of superficial element,   to the computation of the Hilbert polynomials 
of $ \mathbb{M} $ and $\mathbb{M}/a M$ which is available for example in CoCoA system. }} 
 \end{remark}
 
 \vskip 2mm

 \vskip 4mm
 \noindent An     \textbf{$ \mathbb{M}$}-superficial sequence $a_1, \dots, a_r $  for $\q $ is said a { \textbf{maximal  {$\mathbb{M}$}-superficial sequence}} if  $\mathbb{M}/(a_1, \dots, a_r)M $ is nilpotent, but $\mathbb{M}/(a_1, \dots, a_{r-1} )M $ is not nilpotent.  By Proposition \ref{ei},    a superficial sequence is in particular a system of parameters, hence if $ \dim M=r, $ every \textbf{$\mathbb{M}$}-superficial sequence $a_1, \dots, a_r $  for $\q $ is  a  maximal superficial sequence.
 \vskip 2mm 
   In our approach maximal  \textbf{$\mathbb{M}$}-superficial sequences  will play a fundamental role.  We also remark that they are strictly related to   minimal reductions.    
 
  Let  $\mathbb{M}$ be a good $\q$-filtration of $M$ and $J \subseteq \q$ an ideal. Then $J$ is said an  $\mathbb{M}$-reduction of $\q $ if $ M_{n+1} =J M_n $ for $n \gg 0.$ We say that $J$ is a {\it{minimal $  \mathbb{M}$-reduction}} of $\q$  if it is minimal with respect to the inclusion. 

 If $\q$ is $\m$-primary and the residue field is infinite, there is a complete correspondence between maximal  \textbf{$\mathbb{M}$}-superficial sequences for $\q$  and minimal \textbf{$\mathbb{M}$}-reductions of $\q.$  Every minimal $\mathbb{M}$-reduction $J$ of $\q $ can be  generated by a maximal  \textbf{$\mathbb{M}$}-superficial sequence, conversely the ideal generated by a maximal  \textbf{$\mathbb{M}$}-superficial sequence  is a minimal $\mathbb{M}$-reduction of $\q $ (see \cite{HSw} and \cite{Conti}).

\noindent  Notice that if   we consider  for example  $\q=(x^7,x^6y,x^3y^4,x^2y^5,y^7)$ in the power series ring $k[[x,y]] $, then $J=(x^7,y^7) $ is a
minimal reduction of $\q, $ but $\{ x^7,y^7\} $ is not a superficial sequence, in particular
$x^7 $ is not superficial for $\q.$  Nevertheless $\{x^7+y^7, x^7-y^7\} $ is a minimal system of generators of $J$ and it   is  a superficial sequence for $\q. $  
 \vskip 2mm
 
In  this presentation     we prefer to handle  $\mathbb{M}$-superficial sequences with respect to minimal reductions because they have a  better  behaviour in studying Hilbert functions.  
 \vskip 3mm
 
 \section{Maximal Hilbert Functions}
 \vskip 2mm
     Superficial sequences play  an important role in the following result where { {maximal Hilbert functions}}  are described.  The result  was proved in the classical case  in \cite{RVV} Theorem 2.2.  and here is extended to the filtrations of a module which   is not necessarily Cohen-Macaulay.  
 
\vskip 2mm
\begin{theorem} \label{nCM} Let $M$ be a module of dimension $ r\ge 1$ and let $J $ be  the ideal generated by a maximal  \textbf{
$\mathbb{M}$}-superficial sequence
for $\q.$ Then $$ P_{\mathbb{M}} (z) \le \frac{\lambda (M/M_1) +
\lambda (M_1/JM) z }{(1-z)^r}. $$
If the equality holds, then   $ gr_{\mathbb M}(M)$ is Cohen-Macaulay and hence $M$ is Cohen-Macaulay. \end{theorem}

\begin{proof}  We induct on $r. $ Let $r=1 $ and $J=(a). $  We have
 $${{\lambda (M/M_1) +
\lambda (M_1/JM) z }\over{(1-z)}} = \lambda (M/M_1) + \lambda(M/aM)\sum_{j\ge 1}z^j
$$
and $H_{\mathbb {M}}(0)=\lambda(M/M_1).$ 

From the diagram
 $$\begin{matrix}
     M & \supseteq & M_n &\supseteq & M_{n+1} \\
      \| &                  &         &                 & \cup   \\
       M& \supseteq & aM  & \supseteq & aM_n
\end{matrix}$$
we get $$\lambda(M/aM) + \lambda (aM/aM_n)=\lambda(M/M_n) + H_{\mathbb{M}}(n)+\lambda (M_{n+1}/aM_n).$$
On the other hand, from the exact sequence 
$$0 \to  (0:a + M_n)/M_n \to M /M_n \overset {a}\to aM/aM_{n}\to 0$$
we get $$ \lambda(M/M_n) = \lambda (aM/aM_n)+\lambda((0:a + M_n)/M_n). $$
It follows that for every $n\ge 1$ 
\begin{eqnarray} \label{basic}
\lambda(M/aM)= 
H_{\mathbb{M}}(n)+\lambda (M_{n+1}/aM_n)+\lambda((0:a + M_n)/M_n).\end{eqnarray}
This proves that  $H_{\mathbb{M}}(n)\le \lambda(M/aM)$ for every $n\ge 1$ 
and  the first assertion of the theorem follows. 

If we have  equality above then,
by
(\ref{basic}),   for every $n\ge 1$ we get $$\lambda (M_{n+1}/aM_n)=\lambda((0:a + M_n)/M_n)=0.$$ 
This clearly implies that $M_{j+1}:a=M_j$ for every $j\ge 1$ so that, by Lemma \ref{reg}, $a^* \in \q/\q^2 $ is regular on $ gr_{\mathbb{M}}(M)$ and $
gr_{\mathbb{M}}(M)$  is Cohen-Macaulay.

Suppose $ r\ge 2, $  $J=(a_1, \dots, a_r) $ and let us consider the good $\q$-filtration $\mathbb M/a_1M$ on $M/a_1M.$ We have $\dim M/a_1M=r-1$ and we know that $a_2,\dots, a_r$ is a maximal $\mathbb M/a_1M$-superficial sequence for $\q.$  By the inductive assumption and since  $\ \ a_1M\subseteq M_1,\ \ $ we get 

{\small{ \begin{eqnarray*}   P_{\mathbb M/a_1M}  (z) \le   \frac{\lambda [(M/a_1M)/(a_1M+M_1/a_1M)]+\lambda[(a_1M+M_1/a_1M)/K(M/a_1M)]z} {(1-z)^{r-1}}  
 \end{eqnarray*}}}
{\small{ $$=\frac{\lambda (M/M_1) + \lambda(M_1/JM) z }{(1-z)^{r-1}}.$$}}
 where we let $K:=(a_2,\dots,a_r).$   
 
 By using (\ref{s}) and since the power series $\frac{1}{1-z}$ is positive, we get
 $$P_{\mathbb {M}}(z) \le P^1_{\mathbb{M}/a_1M}(z)=\frac{P_{\mathbb{M}/a_1M}(z)}{1-z}\le \frac{\lambda (M/M_1) + \lambda(M_1/JM) z }{(1-z)^{r}},$$ as wanted.
 
 If we have equality, then $$P_{\mathbb{M}/a_1M}(z)=\frac{\lambda (M/M_1) + \lambda(M_1/JM) z }{(1-z)^{r-1}}$$ so that  $gr_{\mathbb{M}/a_1M}(M/a_1M)$ is Cohen-Macaulay and hence    $gr_{\mathbb{M}}(M)$ is Cohen-Macaulay as well by Sally's machine. In particular $M$ is Cohen-Macaulay.
 \end{proof}

The above result  says that  if   the $h$-polynomial  is $h_{\mathbb M}(z)= \lambda(M/M_1)+ \lambda(M_1/JM) z,   $  we may conclude that $ gr_{\mathbb M}(M)$  is Cohen-Macaulay even if we do not assume the Cohen-Macaulyness of $M.$  The result cannot be extended to any "short"  $h$-polynomial  $h_{\mathbb M}(z)= h_0+h_1z.  $  For example if $A= k[[x,y]]/(x^2,xy,xz, y^3) $ and we consider the $\m$-adic filtration, then $P_A(z) = \frac{1+2z}{1-z}, $ but $gr_{\m}(A) \simeq A $ is not Cohen-Macaulay.

 In the classical case of the $\m$-adic filtration on a local Cohen-Macaulay ring $A,$  Elias and Valla in  \cite{EV} proved   that    the $h$-polynomial  of the form $h_{\m}(z)= h_0 + h_1 z + h_2 z^2 $   forces $ gr_{\m}(A)$  to be Cohen-Macaulay.  We cannot extend   this result to general filtrations because this is not longer true even if we consider     the $\q$-adic filtration with $\q$ an $\m$-primary ideal of a Cohen-Macalay ring. The following example is due to Sally (\cite{S4} Example 3.3).

\begin{example} \label{huno}   {\rm{Let $A=k[[t^4,t^5,t^6,t^7]] $  and consider $\q=(t^4,t^5,t^6).$ We have   
$$P_{\q}(z)= {{ 2 +z+z^2}\over (1-z)}$$
and    $ gr_{\q}(A) $  is not  Cohen-Macaulay because $ a=t^4 $ is a superficial (regular) element for $\q, $ but  $\q^2 : a \not = \q $ (cfr. Lemma \ref{sup} and Lemma \ref{reg}).}}
\end{example}

 \chapter{Bounds for  $e_0(\mathbb{M})$ and  $e_1(\mathbb{M})$ } 
In this chapter  we prove lower and upper  bounds for the first two coefficients of the Hilbert polynomial which we defined in Chapter 1.   We recall that $ e_0(\mathbb{M}) $  depends only on $\q$ and $M, $  but it does not depend on the good $\q$-filtration $ \mathbb{M}.  $   In contrast  $e_1(\mathbb{M})   $ does depend  on  the filtration $\mathbb{M}. $  It is called by Vasconcelos {\it tracking number} for its {\it tag}   position  among the different filtrations having the same multiplicity. The coefficient  $e_1(\mathbb{M})   $    is
also called the Chern number (see \cite{V1}). \par  In establishing the properties of $e_1(\mathbb{M}), $    we will need    an {\it{ ad hoc}}  treatment of the 
one-dimensional case.  The results will be extended to higher dimension  via  
  inductive arguments  and via Proposition \ref{ei},  which describes the behaviour of the Hilbert coefficients  modulo superficial elements. In this way we  prove and extend   several classical bounds  on $e_0(\mathbb{M}) $ and  $e_1(\mathbb{M})  $ which we are going   to describe. We start with the so called {Abhyankar-Valla} formula, which gives a natural lower bound for the multiplicity of a Cohen-Macaulay filtered module $M$. The study of Cohen-Macaulay local rings of minimal multiplicity with respect to this bound, was carried out  by J. Sally in \cite {S2}. This paper can be considered as the starting point of much of the recent research in this field.  We extend here our interest to the non Cohen-Macaulay case taking advantage of  the fact that the correction term we are going  to introduce  behaves  well  modulo superficial elements. 
  
Concerning $e_1(\mathbb{M}),$ we extend considerably  the inequality
$$ e_1(\m) \ge e_0(\m) -1 $$
proved by D.G. Northcott in \cite{N}. Besides Northcott's inequality, 
 Theorem \ref{MeN}   extends the corresponding inequality  proved   by Fillmore   in \cite {F} in the case of Cohen-Macaulay modules,    by Guerrieri and Rossi in \cite{GR1} for filtration of ideals and later by Puthenpurakal in  \cite{P1} for   $\q$-adic filtrations of Cohen-Macaulay modules. When  the module $M$ is  not necessarily Cohen-Macaulay, we present a new proof  a recent result by Goto  and Nishida in \cite{GN}.  In our general setting we will focus on an upper bound of $e_1(\mathbb{M}) $ which was introduced and studied in the classical case by Huckaba and Marley in \cite{HM2}.

In the last section, we show   that modules which are extremal with respect to the inequalities proved above  have good associated graded modules and  Hilbert functions of very specific shape.   In some cases we shall see that   extremal values  of the integer $e_1(\m) $  necessarily imply that the ring $A$ is Cohen-Macaulay. These  results  can be considered as a confirmation of the general philosophy of  the paper of W. Vasconcelos \cite{V1},  where the Chern number is conjectured to be a measure of  how far $A$ is from being Cohen-Macaulay.  
\vskip 3mm
 
\section{The multiplicity and the first Hilbert coefficient: basic facts}

In \cite{A}, S. Abhyankar proved a nice lower bound for the multiplicity of a Cohen-Macaulay local ring $(A, \m). $ He found that 
$$ e_0(\m) \ge \mu(\m) -r +1 $$
where $r$ is the dimension of $A$ and $\mu(\m)= H_A(1)  $  is the embedding dimension of $A. $

G. Valla extended the formula to $\m$-primary ideals  in \cite{V}.   Guerrieri and  Rossi in \cite{GR1} showed that the result holds for ideal filtrations.   In \cite{P1}  Puthenpurakal proved the formula for Cohen-Macaulay modules and ideal filtrations by using the idealization of the module.   Here we show that the original proof by Valla extends naturally  to our  general setting.   

\vskip 4mm
As before,  we write the $h$-polynomial of $M$ as $$h_{\mathbb{M}}(z)=h_0(\mathbb{M})+h_1(\mathbb{M})z+\cdots +h_s(\mathbb{M})z^s.$$ Hence, if $\dim(M)=r$,  we have $$h_0(\mathbb{M})=\la(M/M_1) \ \ \ \ h_1(\mathbb{M})=\la(M_1/M_2)-r\la(M/M_1).$$

\noindent If $a$ is an $\mathbb{M}$-superficial element  for $\q, $ then 
$$ h_0(\mathbb{M})= h_0(\mathbb{M}/aM) {\mbox{ but in general }}  h_1(\mathbb{M})\not= h_1(\mathbb{M}/aM).$$
If $a$ is a regular element on $M$,  then
$$  h_1(\mathbb{M}/aM)=  h_1(\mathbb{M}) + \la(M/M_1) - \la(M/M_2:a) \ge h_1(\mathbb{M}) $$
and $  h_1(\mathbb{M}/aM)=  h_1(\mathbb{M}) $ if and only if $M_2:a =M_1.$  

\vskip 3mm
 In the classical case of    the $\m-$adic filtration on a local Cohen-Macaulay ring $A$,     $ h_1(\mathbb{M}) $ is the embedding codimension of $A; $   it  is   positive unless $A$ is a regular ring.   In particular $ h_1(\mathbb{M})= h_1(\mathbb{M}/aM). $  The inequality $  h_1(\mathbb{M})  \le h_1(\mathbb{M}/aM)    $ can be strict. In Example \ref{huno}  we have    $ h_1(\mathbb{M})=1 <   h_1(\mathbb{M}/aM)=2.$ 
 
\noindent  This point makes a   crucial difference between the $\m$-adic filtration and more general  filtrations and it   justifies the new invariant  $h(\mathbb{M}):=\la(M_1/JM+~M_2)$  which will be introduced later (see (\ref{h})).  Notice that if $M$ is Cohen-Macaulay and   $J$ is the ideal  generated by a maximal   \textbf{
$\mathbb{M}$}-superficial sequence such that   $M_2 \cap JM =J M_1, $  then
$$ h_1(\mathbb{M})= h_1(\mathbb{M}/JM). $$ 
Which  assumption is valid if we consider $\mathbb{M}$ the $\q$-adic filtration on a local Cohen-Macaulay ring $A$ 
and $\q $ is  integrally closed (see \cite{H} and \cite{I}).

\begin{proposition} \label{ab}Let $M$ be a module of dimension $ r\ge 1$ and let $J=(a_1,\cdots,a_r) $ be  the ideal generated by an \textbf{
$\mathbb{M}$}-superficial sequence
for the $\m$-primary ideal $\q.$ If  $L:=(a_1,\cdots,a_{r-1})$ and $I\supseteq \q$ is an ideal of $A,$ then 

\begin{equation*}
\begin{split}e_0(\mathbb{M})&=h_0(\mathbb{M})+h_1(\mathbb{M})+r\la(M/M_1)-
\la(JM/JM_1) \\ &\quad \quad \quad \quad \quad \quad \quad +\la(M_2/JM_1)-\la(LM:a_r/LM)\\
&=h_0(\mathbb{M})+\la(M_1/I M_1)-\la(JM/I JM)\\ &\quad \quad \quad \quad \quad \quad \quad+\la(IM_1/ I JM)-\la(LM:a_r/LM).
\end{split}
\end{equation*}
\end{proposition}

\begin{proof} By using Proposition \ref{ei}, we get
\begin{equation*}
\begin{split} e_0(\mathbb{M})
&=e_0(\mathbb{M}/LM)=e_0(\mathbb{M}/JM)-\la(LM:a_r/LM)\\
&=\la(M/JM)-\la(LM:a_r/LM)\\
\end{split}
\end{equation*}

The conclusion easily  
  follows by using the following diagrams

$$\begin{matrix}
     M & \supset & M_1 &\supset & JM  \\
       &         & \cup &       & \cup \\
       &         & M_2 & \supset & JM_1
\end{matrix}\ \ \ \ \ \ \ \ \ \ \ \ \ \ \ \ \ \begin{matrix}
     M & \supset & M_1 &\supset & I M_1 \\
       &         & \cup &       & \cup \\
       &         & JM & \supset & I JM
\end{matrix}$$

\end{proof}

\noindent  If  $M$ is Cohen-Macaulay, since $\q$ is $\m$-primary, then   $ \depth_\q(M) =   \depth_\m (M) $ $ =  r   $  and    the elements $a_1,\dots,a_r$ form a regular sequence on $M.$  Since $J=(a_1,\cdots,a_r) $ is generated by a regular sequence and $J M \subseteq M_1, $ we get   $r\la(M/M_1)=
\la(JM/JM_1)$  and $ \la(JM/I JM)=r\la(M/I M). $ \  Moreover  $ \la(LM:a_r/LM)=0.$ 

As a consequence of the above proposition, we get the following  result. 
\begin{corollary} \label{av} Let $M$ be a Cohen-Macaulay module of dimension $ r\ge 1,$   $J $   the ideal generated by a maximal  \textbf{
$\mathbb{M}$}-superficial sequence
for $\q$  and $I$ an ideal containing $\q.$ Then 
\begin{equation*}
\begin{split}e_0(\mathbb{M})&=h_0(\mathbb{M})+h_1(\mathbb{M})+  \la(M_2/JM_1).
\\ 
e_0(\mathbb{M})&=     \la(M/M_1) -r \la(M/I M) +  \la(M_1/ I M_1) +  \la(I M_1/J I M).
\end{split}
\end{equation*}

\end{corollary}

\noindent  The above result shows in particular that,  if $M$ is  Cohen-Macaulay, then  $\la(M_2/JM_1)$  does  not depend on $J.$ The first formula was  proved by Abhyankar  for the $\m$-adic filtration  and by Valla for an  $\m$-primary ideal $\q.$  The second formula  is due to Goto who proved it  in the case of the $\q$-adic filtration and  $I=\m.$

Since J. Sally made the pioneering  work on Cohen-Macaulay local ring of multiplicity as small as possible, inspired by  the above  formula, it is natural to give the following definitions.

Let  $M$ be Cohen-Macaulay and $ \mathbb{M}$ a good $\q$-filtration. We say that the filtration $ \mathbb{M}$ has {\bf minimal multiplicity} if $$e_0(\mathbb{M})=h_0(\mathbb{M})+h_1(\mathbb{M})$$ or equivalently if $M_2=JM_1.$ 

\vskip 2mm Following  \cite{G}, we say that the filtration $ \mathbb{M}$ has {\bf  Goto minimal multiplicity} with respect to the ideal $I$ if $$e_0(\mathbb{M})=     \la(M/M_1) -r \la(M/I M) +  \la(M_1/ I  M_1) $$ or equivalently $IM_1=JIM.$

Let us compare these two definitions. Given  a  good $\q$-filtration $\mathbb{M}$ on the module $M$ and an ideal $I \supseteq \q,$ we define a new filtration $\mathbb{M}^I $ on $M$ as follows:
\begin{equation} \label{mk} \mathbb{M}^I : \ \ M\supseteq IM\supseteq IM_1\supseteq\dots \supseteq I M_n\dots\supseteq \dots \end{equation} It is clear that $\mathbb{M}^I $ is a good $\q$-filtration on $M$ so that $e_0(\mathbb{M})=e_0(\mathbb{M}^I ).$

\begin {proposition}    $\mathbb{M}$ has  $ Goto$ minimal multiplicity with respect to $I $ if and only if $\mathbb{M}^I $ has minimal multiplicity.

\end{proposition}
\begin{Proof} $\mathbb{M}^I $ has minimal multiplicity if and only if $$e_0(\mathbb{M}^I  )=h_0(\mathbb{M}^I )+h_1(\mathbb{M}^I )$$ This means 
$$e_0(\mathbb{M})=\la(M/I M)+\la(I M/I M_1)-r\la(M/I M).$$ Since we have a diagram $$\begin{matrix}
     M & \supset &  I M   \\
           \cup &      & \cup \\
      M_1 &     \supset    & I M_1
\end{matrix}$$ we get that $\mathbb{M}^I$ has minimal multiplicity if and only if $$e_0(\mathbb{M})=\la(M/M_1)+\la(M_1/I M_1)-r\la(M/I M).$$ This is exactly the condition for $\mathbb{M}$ to have Goto minimal multiplicity.
\end{Proof}

If $\q$ is an $\m$-primary ideal of a local Cohen-Macaulay ring $(A,\m) $ of dimension $r, $ and $\mathbb{M}$ is the $\q$-adic filtration, we get that $e_0(\q)$ is minimal if and only if $\q^2=J\q.$ It is clear that this implies that $gr_{\q}(A)$ is Cohen-Macaulay by Valabrega-Valla. 

On the other hand, by its definition,  the $\q$-adic filtration has Goto minimal multiplicity $e_0(\q)$     with respect to $\m$ if and only if $\q\m=J\m.$

The two notions coincide if $\q=\m.$    We remark  that if $\q$ is integrally closed and if the $\q$-adic filtration has Goto minimal multiplicity, then it has minimal multiplicity. In fact  the condition $\q \m =J \m$ implies $\q^2\subseteq J, $ hence $\q^2=J\q $ since  the integrality guarantees    $\q^2 \cap J= J \q $ (see \cite{H} and \cite{I}). The converse is not longer true. 

\noindent In  general   the condition 
$\q \m =J \m$ seems  far from giving any restriction on the Hilbert function of $\q. $ Notice that it does not imply that $gr_{\q}(A)$ is  Cohen-Macaulay, nor even that  $gr_{\q}(A)$ is  Buchsbaum (see \cite{GN}, Theorem 3.1).

We will study the Hilbert function in the case of minimal multiplicity in Theorem \ref{nor} and Corollary \ref{ind}.

\vskip 2mm

We introduce now  the notion of almost minimal multiplicity. Given    a good $\q$-filtration $\mathbb{M}$ of a Cohen-Macaulay module $M, $ we  say that   $\mathbb{M}$ has {\it{ almost minimal multiplicity}}  if 
$$e_0(\mathbb{M})\ =h_0(\mathbb{M})+h_1(\mathbb{M})+1,  $$
or equivalently $\lambda(M_2/JM_1)=1. $ 
\vskip 2mm
Analogously  we will say that $\mathbb{M} $ has  {\bf{ Goto  almost minimal multiplicity}}  if and only if $\mathbb{M}^I $ has almost minimal multiplicity, equivalently $\lambda(IM_1/JM_1)=1 $ for every $J$ generated by a maximal $\mathbb{M}$-superficial sequence  for $\q.$

We will investigate the Hilbert function of $   \mathbb{M}$ when it has  almost minimal multiplicity. The problem is far more difficult; it amounts to the so called Sally conjecture which was open for several years and finally solved by Wang and independently by Rossi and Valla.
 \vskip 2mm  
We notice  that the concept of almost minimal multiplicity introduced by Jayanthan and Verma (see \cite{JV1}), is equivalent to saying  that the filtration $\mathbb{M}^I  $  has almost minimal multiplicity.  
\vskip 4mm
That's all for the moment, as far as we are concerned with bounds for $e_0.$ Instead, let us come to the first  Hilbert coefficient $e_1.$

 When  $\q$ is an $\m$-primary ideal of the  Cohen-Macaulay local ring $(A,\m)$, $M=A$ and $\mathbb{M}$ the $\q$-adic filtration, in order to prove that $e_1(\q)\ge 0$, Northcott proved a basic  lower bound for $e_1,$ namely $e_1(\q)\ge e_0(\q)-\la (A/\q)\ge 0$ (see \cite{N1}). Fillmore extended it to Cohen-Macaulay modules in \cite {F} (see also \cite{P1}) and later Huneke (see \cite{H}) and Ooishi (see \cite{O}) proved that $e_1(\q)= e_0(\q)-\la(A/\q)$ if and only if $\q^2 = J\q,$ where $J$ is the ideal generated by a maximal superficial sequence for $\q.$   When this is the case, by the Valabrega-Valla criterion,  the associated graded ring is Cohen-Macaulay and the Hilbert function is easily described. This result has been extended to   ideal filtrations of Cohen-Macaulay rings   by Guerrieri and Rossi in \cite{GR1}.  Recently Goto and Nishida in \cite{GN} generalized the inequality, with suitable correction terms, to any local ring not necessarily Cohen-Macaulay  and they  studied the equality  in  the Buchsbaum case. 

The result of Northcott  was improved by Elias and Valla in \cite{EV} where, for the  maximal ideal  of a Cohen-Macaulay local ring $(A,\m),$ one can find a proof of the inequality   $$e_1(\m)\ge 2e_0(\m)-h-2,$$ where $h=\mu(\m) -\dim(A)$ is the embedding codimension of $A.$ Since by Abhyankar 
$$ 2 e_0(\m) -h -2 \ge e_0(\m) -1, $$ this is also an extension of Northcott's inequality. Further, it has been proved that, if equality holds, then the associated graded ring is Cohen-Macaulay and the Hilbert function is determined. This result   was   extended to Hilbert  filtrations of ideals  by Guerrieri and Rossi in \cite{GR1} and rediscovered, only for the  $\q$-adic filtration associated to  an $\m$-primary ideal $\q$, by Corso, Polini and Vasconcelos in \cite{CPV}, 2.9. When equality holds, they need an extra assumption on the Sally module, in order to get a Cohen-Macaulay associated graded ring. We note  that this extra assumption is not essential as already proved in \cite{GR1}, Theorem 2.2. Recently  Corso  in \cite{C} was able to remove  the Cohen-Macaulayness assumption.  Here we extend and improve at the same time all these results in our general setting, by using a very simple inductive argument.

Other notable bounds will be presented.

\section{The $1$-dimensional case}
In establishing the properties of the Hilbert coefficients of a filtered module $M$, it will be convenient to use  induction on the dimension of the module. To start the induction, we need an {\it ad hoc}\  treatment of the one-dimensional case. One dimensional Cohen-Macaulay rings had been extensively studied in the classical case by Matlis in \cite{M}. We present here a general approach for filtered modules which are not necessarily Cohen-Macaulay.

\vskip 3mm
Given a  module {\bf{$M$ of dimension one}}  and a good $\q$-filtration $\mathbb{M}=\{M_j\}_{j\ge 0}$ of $M$, we know that for large $n$ we have $e_0(\mathbb{M})=H_{\mathbb{M}}(n),$ so that we define for every $j\ge 0$ the integers
\begin{equation}\label{uj} u_j(\mathbb{M}):=e_0(\mathbb{M})-H_{\mathbb{M}}(j).\end{equation}

\begin{lemma}\label{uij} Let   $M$ be a module of dimension one. If $a$ is an $\mathbb{M}$-superficial element for $\q,$ then for every $j\ge 0$ we have $$u_j(\mathbb{M})=\la(M_{j+1}/aM_j)-\la(0:_{M_j}a).$$
\end{lemma}
\begin{proof} By Proposition \ref{ei} (3.),  we have \begin{equation*}\begin{split} e_0(\mathbb{M})& =e_0(\mathbb{M}/aM)-\la(0:_Ma)=\la(M/aM)-\la(0:_Ma)\\ &= \la(M/aM_j)-\la(aM/aM_j)-\la(0:_Ma).\\
\end{split}\end{equation*}

By using the following exact sequence $$0\to (0:_Ma)/(0:_{M_j}a)\to M/M_j\to  aM/aM_j\to 0$$ we get 
\begin{equation*}\begin{split}
e_0(\mathbb{M})&=\la(M/aM_j)-\la(M/M_j)+\la((0:_Ma)/(0:_{M_j}a))-\la(0:_Ma)\\
\end{split}\end{equation*} and finally 
\begin{equation*}\begin{split} u_j(\mathbb{M})&=e_0(\mathbb{M})-\la(M_j/M_{j+1})\\ &
= \la(M/aM_j)-\la(M/M_j)-\la(0:_{M_j}a)-\la(M_j/M_{j+1})\\
&=\la(M_{j+1}/aM_j)-\la(0:_{M_j}a).\\ \end{split} \end{equation*} \end{proof}

 It follows that, when  $M$  is one dimensional and Cohen-Macaulay, then $u_j(\mathbb{M})=\la(M_{j+1}/aM_j)$ is  non negative and we have, for every $j\ge 0$, $H_{\mathbb{M}}(j)=e_0(\mathbb{M})-\la(M_{j+1}/aM_j)\le e_0(\mathbb{M}).$

It will be useful to write down the Hilbert coefficients in  terms of  the integers $ u_j(\mathbb{M}).$ 

\begin{lemma}\label{eu} Let  $M$ be a module of dimension one. Then for every $j\ge 1$ we have $$e_j(\mathbb{M})=\sum_{k\ge j-1}\binom{k}{j-1}u_k(\mathbb{M}).$$
\end{lemma}
\begin{proof} We have  $$P_{\mathbb{M}}(z)=\frac{h_{\mathbb{M}}(z)}{1-z}=\sum_{j\ge 0}H_{\mathbb{M}}(j)z^j.$$ Hence, if we write $h_{\mathbb{M}}(z)=h_0(\mathbb{M})+h_1(\mathbb{M})z+\dots +h_s(\mathbb{M})z^s,$ then we get for every $k\ge 1$ $$h_k(\mathbb{M})=H_{\mathbb{M}}(k)-H_{\mathbb{M}}(k-1)=u_{k-1}(\mathbb{M})-u_k(\mathbb{M}).$$ 
Hence we can compute the Hilbert series of $M$
\begin{equation} \label{PZ} P_{\mathbb{M}}(z) = \frac{e_0(\mathbb{M})-u_0(\mathbb{M}) + \sum_{k \ge 1} (u_{k-1}(\mathbb{M})- u_k(\mathbb{M}))z^k}{(1-z)} \end{equation}
Finally \begin{equation*}\begin{split}e_j(\mathbb{M})&=\frac{h_{\mathbb{M}}^{(j)}(1)}{j!}=\sum_{k\ge j}\binom{k}{j}h_k(\mathbb{M}) =\sum_{k\ge j}\binom{k}{j}(u_{k-1}(\mathbb{M})-u_k(\mathbb{M}))\\ & =\sum_{k\ge j-1}\binom{k}{j-1}u_k(\mathbb{M}).\\
\end{split} \end{equation*}
\end{proof}

If we apply the above Lemma when $M$ is Cohen-Macaulay, and using the fact that the integers $u_k(\mathbb{M})$ are non negative,  we get $$e_1(\mathbb{M})=\sum_{k\ge 0}u_k(\mathbb{M})\ge u_0(\mathbb{M})+u_1(\mathbb{M})\ge u_0(\mathbb{M}).$$ Since we have $u_0(\mathbb{M})=e_0(\mathbb{M})-\la(M/M_1),$ and $u_0(\mathbb{M})+u_1(\mathbb{M})=2e_0(\mathbb{M})-\la(M/M_2),$ we trivially  get 
$$e_1(\mathbb{M})\ge e_0(\mathbb{M})-\la(M/M_1),$$ and $$e_1(\mathbb{M})\ge 2e_0(\mathbb{M})-\la(M/M_2),$$ which are the bounds   proved by Northcott and by Elias-Valla, in the one dimensional Cohen-Macaulay case.

\vskip 3mm
If we do not assume that $M$ is Cohen-Macaulay, then the integers $u_i(\mathbb{M})$ can be negative, so that the equality $e_1(\mathbb{M})=\sum_{k\ge 0 }u_k(\mathbb{M})$ no longer imply   $e_1(\mathbb{M})\ge u_0(\mathbb{M}).$  For instance, if we consider for example the local ring $A=k[[X,Y]]/(X^2,XY)$  endowed with  the $\m$-adic filtration, we have $e_0=1,$ $e_1=-1$ and  $u_0=0.$ 
\vskip 3mm
Hence we need to change this formula by introducing   some  correction terms which vanish   in the Cohen-Macaulay case.
\vskip 3mm

Given  a good $\q$-filtration $\mathbb{M}=\{M_j\}_{j\ge 0}$  of the $r$-dimensional module $M,$ let  $a_1,\cdots,a_r $ be  an  $\mathbb{M}$-superficial sequence  for $\q;$ further  let $J:=(a_1,\cdots,a_r)$ and $$\mathbb{N}:=\{J^jM\} $$  be the $J$-adic filtration on $M$ which is clearly $J$-good.
 It is not difficult to prove that also the original filtration $\mathbb{M}$ is $J$-good and  this implies  that $e_0(\mathbb{M})=e_0(\mathbb{N}).$

When  $M$ is Cohen-Macaulay, the elements $a_1,\dots,a_r $ form a regular sequence on $M$ and, as a consequence, one can prove  $$J^iM/J^{i+1}M\simeq (M/JM)^{\binom{r+i-1}{i}}.$$ This implies that  the Hilbert Series of $\mathbb{N}$ is $P_{\mathbb{N}}(z)=\frac{\la(M/JM)}{(1-z)^r}$ and thus   $e_i(\mathbb{N})=0$ for every $i\ge 1.$ This proves that   these integers  give a good measure of how $M$ differs from being  Cohen-Macaulay.
 \vskip 2mm

We prove that $ e_1(\mathbb{N}) \le  0$  in the one dimensional case by relating  the integer $e_1(\mathbb{N})$ to  the 0-th local cohomology module of $M.$

 The following Lemma is a partial confirmation  of a conjecture stated by Vasconcelos in \cite{VW4} and \cite{V1}  concerning the negativity of $e_1(J) $ in the higher dimensional case.  Interesting results concerning  this problem have been recently proved in \cite{GHV} and in  \cite{GGHOV}.
\vskip 2mm
  
  In the following we denote by $W$ the $0$-th local cohomology module $H^0_{\m}(M)$ of $M$  with respect to $\m.$ We know that $H_{\m}^0(M):=\cup_{j\ge 0}(0:_M\m^j)=0:_M\m^t$ for every $t\gg 0.  $ 
  
  In the $1$-dimensional case we have the following nice formula (see \cite{GN}, Lemma 2.4).
  \begin{lemma} \label{W} Let $M$ be a finitely generated $A$-module of dimension one  and let $a$ be a parameter for $M.$  Then for every $t \gg 0$ we have $W=0:_Ma^t$  and, if we denote by $\mathbb{N}$  the $(a)$-adic filtration on $M,$  then  $$\la(W)=-e_1(\mathbb{N}).$$
\end{lemma}
\begin{proof} Since there is an integer $j$ such that  $\m^jM\subseteq aM=((a)+0:M)M$, the ideal $(a)+0:M$ is $\m$-primary and therefore $\m^s\subseteq  (a)+0:M$ for some $s;$ this implies     $$\m^{ts}\subseteq (a)^t+0:M$$ for every $t.$ On the other hand, $W=0:_M\m^t$ for every  integer $t \gg 0$, so that 
$$W=0:_M\m^t\subseteq 0:_Ma^t\subseteq 0:_M\m^{ts}=W.$$
We denote by $\mathbb{N}^n$ the $(a^n)$-adic filtration on $M.$  Now   for $n \gg 0, $ it is easy to see that $ n e_0(\mathbb{N})= e_0( \mathbb{N}^n)=\la(M/a^nM)  - \la(0:_M a^n)$ and  the result follows because $ \la(M/a^nM) = n e_0(\mathbb{N}) - e_1(\mathbb{N}). $ 

\end{proof}

Given a good $\q$-filtration of the module $M$ (any dimension), we consider now the corresponding filtration of the saturated module $M^{sat}:=M/W.$ This is the filtration $$\mathbb{M}^{sat}:=\mathbb{M}/W=\{M_n+W/W\}_{n\ge 0.}$$ Since $W$ has finite length and $\cap_{i\ge 0}  M_i=\{0\},$ we have $M_i\cap W=\{0\}$ for every $i\gg 0.$ This implies  $p_{\mathbb{M}}(X)=
p_{\mathbb{M}^{sat}}(X).$

\vskip 2mm Further, it is clear that for every $j\ge 0$ we have an exact sequence $$ 0\to W/(M_{j+1}\cap W)\to M/M_{j+1}\to M/(M_{j+1}+W)\to 0$$ so that for every $j\gg 0$ we have $$\la ( M/M_{j+1})=\la [M/(M_{j+1}+W)]+\la(W)$$ which implies \begin{equation} \label {sat} p^1_{\mathbb{M}}(X)=p^1_{\mathbb{M}^{sat}}(X)+\la(W).\end{equation} This proves the following result:
\begin{proposition}\label{e1(sat)} Let $M $ be a module of dimension $r. $ Denote   $W:=H^0_{\m}(M) $   and $\mathbb{M}^{sat}:=\mathbb{M}/W. $ Then  $$ e_i(\mathbb{M})=e_i({\mathbb{M}}^{sat}) \ \    0\le i \le r-1,\ \ \ \ \ \ \ \ \
e_r(\mathbb{M})=e_r({\mathbb{M}}^{sat})+(-1)^d\la(W).$$
\end{proposition}

We remark that, if $\dim(M)\ge 1$, the module $M/W$   always has  positive depth. This is the reason why, sometimes,  we move our attention from the module $M$ to the module $M/W.$ This will be  the strategy of the proof of  the next proposition which gives,  in the one dimensional case, the promised upper bound for $e_1.$

\begin{proposition}\label{e1<=} Let $\mathbb{M}=\{M_j\}_{j\ge 0} $ be a good $\q$-filtration of a module $M$ of dimension one. If $a$ is an  $\mathbb{M}$-superficial element for $\q$ and $\mathbb{N}$ the $(a)$-adic filtration on $M,$ then
$$e_1(\mathbb{M})-e_1(\mathbb{N})\le \sum_{j\ge 0}\la(M_{j+1}/aM_j).$$  If  $W\subseteq M_1$ and  equality holds above, then  $M$ is Cohen-Macaulay.\end{proposition}
\begin{proof} By Proposition \ref{e1(sat)}  and Proposition \ref{W} we have $$e_1(\mathbb{M})=e_1(\mathbb{M}^{sat})-\la(W)=e_1(\mathbb{M}^{sat})+e_1(\mathbb{N})$$ so that we need to prove that $e_1(\mathbb{M}^{sat})\le \sum_{j\ge 0}\la(M_{j+1}/aM_j).$

Now $M/W$ is Cohen-Macaulay and $a$ is regular on $M/W,$ hence by Lemma \ref{eu} and Lemma \ref{uij}, we get 
\begin{equation*}
\begin{split} 
e_1(\mathbb{M}^{sat})&= 
\sum_{j\ge 0}u_j(\mathbb{M}^{sat})=\sum_{j\ge 0}\la(M_{j+1}^{sat}/aM_j^{sat})\\
&=\sum_{j\ge 0}\la\left[\frac{M_{j+1}+W}{aM_j+W}\right]
=\sum_{j\ge 0}\la\left[\frac{M_{j+1}}{aM_j+M_{j+1}\cap W}\right]  \\
& \le     \sum_{j\ge 0}\la(M_{j+1}/aM_j).
\end{split}
\end{equation*} The first assertion  follows.  In particular equality holds if and only if $   \ \ M_{j+1} \ \cap W~\subseteq ~aM_j$ for every $j\ge 0.$ Let as assume $W\subseteq M_1$ and equality above; then we have $W=W\cap M_1\subseteq aM.$ Now recall that $W=0:_Ma^t$ for $t\gg 0,$ hence if $c\in W$ then $c=am$ with $a^tc=a^{t+1}m=0.$ This implies  $m\in 0:_Ma^{t+1}=W$ so that $W\subseteq aW$ and, by Nakayama, $W=0.$
\end{proof}

\vskip 2mm
  
   \vskip 5mm
We turn out  to describing  lower bounds on the first Hilbert coefficient thus extending the classical result proved by Northcott.

\begin{proposition} \label{d=1} Let $\mathbb{M}=\{M_j\}_{j\ge 0}$ be a good $\q$-filtration of a module $M$ of dimension one.  If $a $ is an $\mathbb{M}$-superficial  element for $\q$ and $s\ge 1$ a given  integer,  then for every $n\gg 0$ we have 

\begin{equation*}
\begin{split}e_1(\mathbb{M})-e_1(\mathbb{N})&= s \, e_0(\mathbb{M})-\la (M/M_s)+\la(M_s+W/M_s)+\la(M_n/a^{n-s}M_s)
\\ &=\sum_{j=0}^{s-1}u_j(\mathbb{M})+\la(M_s+W/M_s)+\la(M_n/a^{n-s}M_s).\end{split}
\end{equation*}

\end{proposition}
\begin{proof} We have for  every $n\gg 0$ the following equalities:
$$\la(M/M_n)=p^1_{\mathbb{M}}(n-1)=e_0(\mathbb{M})n-e_1(\mathbb{M})$$ 
$$\la(M/a^{n-s}M)=p^1_{\mathbb{N}}(n-s-1)=e_0(\mathbb{N})(n-s)-e_1(\mathbb{N}).$$ Since $e_0(\mathbb{M})=e_0(\mathbb{N})$ , we get   
$$e_1(\mathbb{M})-e_1(\mathbb{N})=s\, e_0(\mathbb{M})-\la(M/M_n)+\la(M/a^{n-s}M).$$

From the diagram $$\begin{matrix}
     M & \supset & M_n  \\
      \cup &  & \cup \\
     a^{n-s}M & \supset & a^{n-s}M_s
\end{matrix}$$ we get 
$$e_1(\mathbb{M})-e_1(\mathbb{N})=s\, e_0(\mathbb{M})+\la(M_n/a^{n-s}M_s)-\la(a^{n-s}M/a^{n-s}M_s).$$ By using  the exact sequence $$0\to (M_s+0:_Ma^{n-s}/M_s)\longrightarrow  M/M_s \overset{a^{n-s}}\longrightarrow  a^{n-s}M/a^{n-s}M_s\longrightarrow 0$$ and the equality $0:_Ma^t=W$ for $t\gg 0,$ we get the conclusion.
\end{proof}

\begin{corollary} \label{cor} With the same notation as in the above Proposition, if $$e_1(\mathbb{M})-e_1(\mathbb{N})= s \, e_0(\mathbb{M})-\la (M/M_s)+\la(M_s+W/M_s),$$ then $M_{s+1}\subseteq aM_s+W.$
\end{corollary}
\begin{proof}

We simply notice  that we have an injective map
$$\left (M_{s+1}+W\right )/\left (aM_s+W\right ) \overset {a^{n-s-1}}\longrightarrow M_n/a^{n-s}M_s.$$
\end{proof}

\vskip 2mm   The converse does not hold, as the following example shows. Let $A=k[[t^3,t^4,t^5]]$ and let us
 consider the following $\m$-filtration $\mathbb{M}$ on $A:$ $$M=A, \ \ M_1=\m, \ \ M_2=\m^2, \ \ M_3=\m^2,\ \ M_j=\m^{j-1}$$ for $j\ge 4.$  
It is clear that  $t^3$
 is an $\mathbb{M}$-superficial element for $\m$  and $A$ is Cohen-Macaulay so that $W=0$. We have
$$P_{\mathbb{M}}(z)=\frac{1+2z-3z^2+3z^3}{1-z}, \ \ \ P_{\mathbb{N}}(z)=\frac{3}{1-z}$$ so that $e_0(\mathbb{M})=3,$
$e_1(\mathbb{M})=5$, $e_1(\mathbb{N})=0$ and the equality $e_1(\mathbb{M})-e_1(\mathbb{N})= e_0(\mathbb{M})-\la
(M/M_1)$ does not hold even if $M_2=t^3M_1$.
\vskip 2mm

The following  result  was proved in \cite{GN}, Lemma 2.1.   in the case $M=A$ and  $s=1.$

\begin{corollary} \label{corr} Let $\mathbb{M}=\{\q^jM\}_{j\ge 0}$ the $\q$-adic filtration on $M$ of dimension one.  If $a\in \q$ is an $\mathbb{M}$-superficial  element for $\q$ and $s\ge 1$ a given  integer.   Then   $$e_1(\mathbb{M})-e_1(\mathbb{N})= s \, e_0(\mathbb{M})-\la (M/M_s)  $$ if and only if  $M_{s+1}\subseteq aM_s+W $ and $W \subseteq M_s.$
\end{corollary}
\begin{proof}
By using Corollary \ref{cor} and Proposition \ref{d=1}, it is enough 
   to prove that $M_{s+1}\subseteq aM_s+W$ implies $ M_n=a^{n-s}M_s$ for $n\gg 0.$ 
     
We have  $M_{s+1} \subseteq aM_s+W$  and by multiplication by $\q^{n-s} $ the result follows since    $M_{n+1} \subseteq a M_n+\q^{n-s}W=a M_n $ for every $n\gg 0.$  
\end{proof}

Under the assumption dim $M=1, $  we recover    Theorem 1.3. in  \cite{GN}  and we give a positive answer to a question raised by
Corso in \cite{C}. 

\begin{theorem}\label{Buchsb}  Let $\mathbb{M}=\{\m^jM\}_{j\ge 0}$ be the $\m$-adic filtration of a Buchsbaum  module $M$ of dimension one.  Assume either  
\par \vskip 2mm $i)\  e_1(\mathbb{M})-e_1(\mathbb{N})=   \, e_0(\mathbb{M})-\la (M/M_1)$
\par  \noindent or  
\par ii) \ 
$ e_1(\mathbb{M})-e_1(\mathbb{N})= 2 \, e_0(\mathbb{M})-\la (M/M_2). $ 
  
\par \vskip 2mm \noindent   Then gr$_{\mathbb{M}}(M) $ is Buchsbaum. 
\end{theorem} 
\begin{proof} By the above corollary,     if either i) or ii) holds, then $M_{n+1} \subseteq  a M_n $ + W for every $n \ge 2.$
Then by Valabrega-Valla criterion applied to $ \mathbb{M}^{sat}= \mathbb{M}/W, $ it follows  that    $gr_{ \mathbb{M}^{sat}}(M^{sat}) $ is Cohen-Macaulay.   Denote  by $Q $ the graded maximal ideal of $ $gr$_{\m}(A). $    By using the fact that   gr$_{ \mathbb{M}^{sat}}(M^{sat}) $ is Cohen-Macaulay and $\m W=0,  $ it easy to see that   $Q \  H^0_{Q}(gr_{\mathbb{M}}(M)  )=0 $ and the result follows by \cite{SV}, Proposition 2.12.
\end{proof}

\vskip 3mm

For further applications,  we need to consider another filtration 
related to a superficial sequence. Given  a good $\q$-filtration $\mathbb{M}=\{M_j\}_{j\ge 0}$  of the $r$-dimensional module $M,$ let  $a_1,\cdots,a_r $ be  an  $\mathbb{M}$-superficial sequence  for $\q$ and  let $J:=(a_1,\cdots,a_r).$ We define the following filtration $$\mathbb{E}:\ \ \ \ M\supseteq M_1\supseteq JM_1 \supseteq J^2M_1 \supseteq \cdots \supseteq J^jM_1 \supseteq J^{j+1}M_1\supseteq \cdots $$ 
This filtration is $J$-good and we want to compare it with   $\mathbb{N},$ the $J$-adic filtration on $M.$ We need to remark that we have $$e_0(\mathbb{E})=e_0(\mathbb{M})=e_0(\mathbb{N}).$$

\begin{proposition}\label{EeN} Let $\mathbb{M}=\{M_j\}_{j\ge 0}$ be a good $\q$-filtration on $M$ of dimension one.   If $a $ is an $\mathbb{M}$-superficial element for $\q$,  then  we have 
$$e_1(\mathbb{E})-e_1(\mathbb{N})=  e_0(\mathbb{M})-h_0(\mathbb{M})+\la(M_1+W/M_1).$$
In particular $e_1(\mathbb{E})-e_1(\mathbb{N})=  e_0(\mathbb{M})-h_0(\mathbb{M})$ if and only if $ M_1 \supseteq  W.$
\end{proposition}
\begin{proof} As in the above proposition, we have for  every $n\gg 0$ the following equalities:
$$\la(M/a^nM_1)=p^1_{\mathbb{E}}(n)=e_0(\mathbb{M})(n+1)-e_1(\mathbb{E})$$ 
$$\la(M/a^{n}M)=p^1_{\mathbb{N}}(n-1)=e_0(\mathbb{M})n-e_1(\mathbb{N}).$$ We get   
\begin{equation*}
\begin{split}e_1(\mathbb{E})-e_1(\mathbb{N})&=e_0(\mathbb{M})-\la(M/a^nM_1)+\la(M/a^{n}M)\\ &=e_0(\mathbb{M})-\la(a^nM/a^nM_1)\\ &= e_0(\mathbb{M})-\la (M/M_1)+\la(M_1+W/M_1)\end{split}
\end{equation*} where the last equality follows from the exact sequence 
$$0\to (M_1+W)/M_1 \to M/M_1 \overset{a^n}\to a^nM/a^nM_1 \to 0.$$
\end{proof}

We conclude this   section about  the one dimensional case, with a notable extension of a bound on $e_1$ proved in the classical case by D. Kirby and extended to an $\m$-primary ideal by M.E. Rossi, G. Valla and W. Vasconcelos.  We assume $M$ is a Cohen-Macaulay module of dimension one.  We recall that 
  the integers $u_j$ are non negative because   we have $$u_j(\mathbb{M})=\la(M_{j+1}/aM_j).$$ 
In particular for  every $j \ge 0$ we have $$ H_{\mathbb{M}}(j) =e_0(\mathbb{M}) - \la(M_{j+1}/aM_j)  \le e_0(\mathbb{M}) $$
where $a$ is an  $\mathbb{M}$-superficial element for $\q.$
We define  
\begin{equation} \label{reduction} s(\mathbb{M}):= \min \{ j   :  H_{\mathbb{M}}(j)= e_0(\mathbb{M})\}= \min \{ j   :   M_{j+1}= a M_j \}.\end{equation}
 This  integer is the  {\it{  reduction number of $\mathbb{M}$ }} and   from the above equality it clear that it does not depend on $(a).$   Moreover  (\ref{PZ})  shows  that  $s(\mathbb{M}) $ is also the degree of the $h$-polynomial of ${\mathbb{M}}.$
\vskip 3mm
If $ A$ is a Cohen-Macaulay local ring of dimension one and we consider  the classical $\m$-adic filtration on $A,$ Sally proved  that $$s(\m)  \le   e_0(\m) -1$$ (see \cite{SaV} and \cite{S4}).  This result easily follows by a lower bound on the minimal number of generators proved by Herzog and Waldi (\cite{HW},Theorem 2.1). Herzog and Waldi's result can be generalized to modules in the case of the $\q$-adic filtration on $M.$ The proof    is a straightforward adaptation ion of the classical case and we give here a proof for completeness.  In the following $\mu(\cdot  ) $ denotes
the minimal number of generators of a module on a local ring and we assume $ s(\mathbb{M}) \ge 1.$ If 
$ s(\mathbb{M}) =  1 $ the statements become trivial.

 \begin{theorem}\label{herzog-waldi} \rm{(Herzog-Waldi)}   Let $\mathbb{M} $ be the $\q$-adic filtration of a Cohen-Macaulay module $M$ of dimension one.  Then 
 \vskip 2mm
 1. $\mu( \q^n M) > n $ for every $ n \le s(\mathbb{M}). $   
 \vskip 2mm
 2. $\mu( \q^n M) =  \mu( \q^{n+1} M)$  for every  $ n \ge  s(\mathbb{M}). $   
  \end{theorem}  
  
  \begin{proof} We prove 1. Let $a$ be an $\mathbb{M}$-superficial element for $\q, $ since  $ s=s(\mathbb{M}) \ge 1$
  and $\q^s M \not = a \q^{s-1}M, $ then there exist $x_1, \dots, x_s \in \q $ and $m \in M$ such that 
  $$x_1 \cdots x_s m \in \q^sM,    \mbox{ \  \ but \  \  }  x_1 \cdots x_s m \not \in  a \q^{s-1}M + \m \q^sM.$$
  For every $i=0, \dots, s $ we consider  $$y_i:= a^i x_{i+1} \cdots x_s m $$  and we claim that $\{ y_0=    x_{1} \cdots x_s m, \ y_1= a x_{2} \cdots x_s m, \dots \dots, y_s= a^s m\} $ is  part of a minimal system of generators of $\q^sM.$   Assume    $ r_0y_0 + r_1 y_1 + \dots + r_sy_s \in   \m \q^sM  $ with  $r_i \in A$ and we conclude  $r_i \in \m $ by arguing step by step on $ i=0, \dots, s.$ First $r_0 \in \m $  otherwise $y_0   \in  a \q^{s-1}M + \m \q^sM, $ Assume $i >0 $ and $r_0, \dots, r_{i-1} \in \m, $ hence  $ r_iy_i + r_{i+1} y_{i+1} + \dots + r_sy_s \in   \m \q^sM. $ Multiplying the sum with $x_1 \cdots x_i $ we obtain
   $$ a^i( r_iy_0 + r_{i+1} y_{i+1} x_1 \cdots x_i + \dots + r_sy_s x_1 \cdots x_i) \in   \m \q^{s+i} M. $$
   Since $\q^{s+i} M = a^i \q^s M $ and $ a $ is regular on $M, $  we get  $ r_iy_0 + r_{i+1} y_{i+1} x_1 \cdots x_i + \dots + r_sy_s x_1 \cdots x_i \in a \q^{s-1} M + \m \q^s M, $ therefore $r_i \in \m.$
   
   If $n \le s $ then it is easy to see that    $ y_{i,n}=   a^i  x_{i+1} \cdots x_{ n } m, $ $ i=0, \dots, n $  is  part of a minimal system of generators of $\q^nM.$   In fact  multiplying by $x_{n+1} \cdots x_s $ we   map the elements $ \{ y_{0,n}, \dots, y_{n,n} \} $  of $\q^nM $   onto  $\{y_0, \dots, y_n\} $ which is part of a minimal set of generators of 
 $\q^s M   $  and 1. is proved. We remark now that 2. is a trivial consequence of the definition of  $s(\mathbb{M}) $ and of  the fact that $a$ is $M$-regular. \end{proof}

  Here we  extend  Sally's result to our general setting.  

 \begin{proposition}\label{HW} Let $\mathbb{M} $ be the $\q$-adic filtration of a Cohen-Macaulay module $M$ of dimension one. Let $p $ be an integer such that $\q M \subseteq \m^p M$. Then 
 \vskip 2mm
 1. $H_{\mathbb{M}}(n) \ge  n +p  $ for every $ n \le s(\mathbb{M}) $ 
 \vskip 2mm  2. $  s(\mathbb{M}) \le  e_0(\mathbb{M}) - p.$  
  \end{proposition} 
\begin{proof}   Since 
$$ \q^nM \supseteq \m \q^nM  \supseteq \m^2 \q^nM \supseteq \dots \supseteq \m^p \q^n M \supseteq \q^{n+1}M $$
we get  $ \ \ \  H_{\mathbb{M}}(n)=  \la(\q^n M/ \q^{n+1}M) = \mu(\q^n M) + \mu(\m \q^nM)+ \dots + \mu(\m^{p-1} \q^nM) +  \la( \m^p \q^n M/ \q^{n+1}M ). $    \par \noindent  Now,  by   Theorem \ref{herzog-waldi},    we get   $H_{\mathbb{M}}(n) \ge  n +p  $  if $ n \le s(\mathbb{M}). $ 

Since   $H_{\mathbb{M}}(n) \le e_0(\mathbb{M}), $   as consequence of 1. we get   $H_{\mathbb{M}}(e_0(\mathbb{M}) - p)= e_0(\mathbb{M}).   $ Hence   
$$  s(\mathbb{M}) \le  e_0(\mathbb{M}) - p.$$  \end{proof}
 
 \vskip 3mm
\noindent From the proof of the above result we remark that  if $  \m^p \q^n M \neq \q^{n+1}M, $  then
 \vskip 2mm  (a)  $H_{\mathbb{M}}(n) >  n +p $ if $ n \le s(\mathbb{M}) $ 
\vskip 2mm (b)   $  s(\mathbb{M}) \le  e_0(\mathbb{M}) - p -1.$ 
  
\vskip 5mm

As a consequence of Proposition \ref{HW} and Proposition \ref{e1(sat)},  we obtain the following result which  refines a classical result proved by Kirby  for the maximal ideal  in  \cite{K}  and  extended to the  $\m$-primary  ideals in \cite{RVV}.   
 \vskip 2mm
 
\begin{proposition} \label{kirby} Let $\mathbb{M} $ be the $\q$-adic filtration of  a module  $M$ of dimension one.   Let $p$ be an integer such that $\q M \subseteq \m^p M$. Then
$$ e_1(\mathbb{M}) -e_1(\mathbb{N})  \le {{e_0(\mathbb{M}) -p +1}\choose 2}.$$
If $e_0(\mathbb{M}) \neq e_0(\m M), $ then
$$ e_1(\mathbb{M}) -e_1(\mathbb{N})  \le {{e_0(\mathbb{M}) -p  }\choose 2}.$$
 \end{proposition} 
\begin{proof} We recall that  $ e_0(\mathbb{M}) =e_0(\mathbb{M}^{sat}) $ and, if $\q M \subseteq \m^pM, $ clearly 
$\q M^{sat} \subseteq \m^p  M^{sat}.$ Then, by Proposition \ref{e1(sat)}, we may assume $M$ is a  Cohen-Macaulay module and $s({\mathbb{M}}) \ge 1.$ \par \noindent Since $e_1(\mathbb{M})= \sum_{j \ge0} (e_0(\mathbb{M})-H_{\mathbb{M}}(j)), $ by Proposition  \ref{HW},   it follows  that $$e_1(\mathbb{M})= \sum_{j =0}^{e_0(\mathbb{M})-p } (e_0(\mathbb{M})-H_{\mathbb{M}}(j) ) 
\le  \sum_{j =0}^{e_0(\mathbb{M})-p } (e_0(\mathbb{M})- j -p ) =  {{e_0(\mathbb{M}) -p +1}\choose 2}.$$
The last assertion is a consequence of  the remark  (a)  following Proposition  \ref{HW}. 
\end{proof}

\section{The higher dimensional  case}

We come now to the {\it{higher dimensional}} case. Here the strategy is  to lower dimension by using superficial elements. We do not assume that $M$ is Cohen-Macaulay, so  we will get formulas containing a correction term which vanishes in the Cohen-Macaulay case.  

If  $J:=(a_1,\cdots,a_r)$ is an  $\mathbb{M}$-superficial sequence  for $\q, $ let  $$\mathbb{N}:=\{J^jM\} $$  be the $J$-adic filtration on $M.$  It is not difficult to prove that  $e_0(\mathbb{M})=e_0(\mathbb{N}).$
 
If  $M$ is Cohen-Macaulay, then $e_i(\mathbb{N})=0$ for every $i\ge 1$, so that these integers are good candidates for being correction terms  when the Cohen-Macaulay assumption does not hold.
 Concerning $e_1(\mathbb{M}),$ we can say that  $ e_1(\mathbb{N}) $  is  the  penality for the lack of that condition.    This  character was  already used by Goto in studying the Buchsbaum case.  If  $M$ is a generalized Cohen-Macaulay
module,  then   following  \cite{GN},  we have  
$$ e_1(\mathbb{N})  \geq - \sum_{i=1}^{r-1}{r-2\choose
i-1}\lambda(H^i_{\mathfrak{m}}(M)) $$ 
with equality if $M$ is Buchsbaum. Hence if $M$ is Buchsbaum, then $ e_1(\mathbb{N})  $ is independent of $J.$ Very recently, under suitable assumptions and when $M=A,$  Goto and Ozeki proved that if  $ e_1(\mathbb{N})  $ is independent of $J, $  then  $A $ is Buchsbaum   (see \cite{GO1}).

\vskip 2mm

\noindent 
The following Lemma is the key to our investigation.  It is due to David Conti.

\begin{lemma}\label{key} Let $M_1,\dots,M_d $ be $ A$-modules of dimension $r,$  let $\q$ be an  $\m$-primary ideal of $A$ and $\M_1=\{M_{1,j}\},\dots,\M_d=\{M_{d,j}\}$ be good $ \q$-filtrations of $M_1,\dots,M_d$ respectively. Then we can find elements $a_1,\dots,a_r $ which are $\M_i$-superficial for $\q$ for every $i=1,\dots,d. $

\vskip 2mm \noindent If $r \ge 2$ and  $M_1=\dots=M_d =M,$  then  for every $1\le i\le j\le d$ we have $$e_1(\M_i)-e_1(\M_i/a_1M)=e_1(\M_j)-e_1(\M_j/a_1M).$$
\end{lemma}
\begin{proof} We have a filtration of the module $\oplus^d_{i=1}\M_i$
$$  \  \oplus_{i=1}^d M_i\supseteq \oplus_{i=1}^d M_{i,1}\supseteq \oplus_{i=1}^d M_{i,2}\supseteq  \dots \supseteq \oplus_{i=1}^d M_{i,j}\supseteq  \dots$$ which we denote by 
$\oplus^d_{i=1}\M_i. $  It is clear that this is a good $\q$-filtration on $\oplus_{i=1}^d M_i.$ Let us choose a $\oplus^d_{i=1}\M_i$-superficial  sequence
$\{a_1,\dots,a_r\}$ for $\q$. Then it is easy to see that $\{a_1,\dots,a_r \}$ is a sequence of $M_i$-superficial  elements
for $I $ for every $i=1,\dots,d.$ This proves the first assertion. As for the second one,  we know that $$e_1(\M_i)-e_1(\M_i/a_1M)=\begin{cases} \la(0:_Ma_1) & \text{if}\ \ r =2\\
0 & \text{if} \ \ r \ge 3\end{cases}$$ from which the conclusion follows.
\end{proof}

\vskip 4mm 
We start with the extension of 
Proposition \ref{e1<=}  to the higher dimensional case. To this end, given  a good $\q$-filtration $\mathbb{M}=\{M_j\}_{j\ge 0}$ of $M$ and      an ideal $J$ generated by a maximal sequence of $\mathbb{M}$-superficial  elements 
for $\q,$  we let for every $j\ge 0,$ \begin{equation}\label{vj}v_j(\mathbb{M}):=\la (M_{j+1}/JM_j).\end{equation}
When  $M$ is one-dimensional and Cohen-Macaulay one has $$v_j(\mathbb{M})=u_j(\mathbb{M})$$ where the $u_j$'s are defined as in 
(\ref{uj}).

\begin{proposition}\label{Vi} Let $\mathbb{M} $ be a good $\q$-filtration of a module $M$ of dimension $r$ and    $J$  an ideal generated by a maximal sequence of $\mathbb{M}$-superficial  elements 
for $\q;$ then we have $$e_1(\mathbb{M})-e_1(\mathbb{N})\le \sum_{j\ge 0}v_j(\mathbb{M}).$$
\end{proposition}
\begin{proof} If $\dim(M)=1,$ then we can apply Proposition  \ref{e1<=}. Let  $\dim(M)\ge 2;$  by Lemma \ref{key},  we can find a minimal system of generators  $\{ a_1,\cdots,a_r\}$ of $J$ such that $a_1$ is $\mathbb{N}$-superficial for $J$ and $\{a_1,\cdots,a_r\}$ is a sequence of $\mathbb{M}$-superficial  elements 
for $\q.$ The module $M/a_1M$ has dimension $r-1$ and $\mathbb{M}/a_1M$ is a good $\q$-filtration on it. Further it is clear that  $a_2,\cdots , a_r,$  is a maximal sequence of $\mathbb{M}/a_1M$-superficial  elements 
for $\q.$ Hence, if we let $K$ be the ideal generated by $a_2,\cdots , a_r,$  then by using induction on $\dim(M),$  we get
\begin{equation*}
\begin{split} 
e_1(\mathbb{M})-e_1(\mathbb{N})&=e_1(\mathbb{M}/a_1M)-e_1(\mathbb{N}/a_1M)
\\ &\le \sum_{j\ge 0}\la\left ( (M_{j+1}+a_1M)/(KM_j+a_1M)\right )
\\ & =\sum_{j\ge 0}\la\left ( (M_{j+1}+a_1M)/(JM_j+a_1M)\right )
\\ &=\sum_{j\ge 0}\la\left ( M_{j+1}/(JM_j+(a_1M\cap M_{j+1})\right )
\\ & \le \sum_{j\ge 0}\la\left ( M_{j+1}/JM_j\right )\\
\end{split} 
\end{equation*}
\end{proof}

In the classical case, when $M=A,$ $\mathbb{M}=\{\q^j\}$ with $\q$ an $\m$-primary ideal of the $r$-dimensional Cohen-Macaulay local ring $A,$ the above inequality is due to  S. Huckaba (see \cite{H1}). Huckaba also proved that equality holds if and only if the associated graded ring has depth at least $r-1.$ We will extend  this result in the next section (see Theorem \ref{e1=CM} and Theorem  \ref{th2} ).

\vskip 2mm

We move now to the extension of Proposition \ref{kirby}.
If $(A, \m) $ is a Cohen-Macaulay local ring,   we  recall that   Kirby (see \cite{K}) proved  
$$ e_1(\m) \leq { {e_0(\m)}\choose 2}$$
If $\q $ is an $\m-$primary ideal, the result has  been extended in \cite{RVV}.   In particular if $e_0(\q) \neq e_0(\m), $ then
$$ e_1(\q) \leq { {e_0(\q)-1 }\choose 2}$$
 
We improve the above results   by using the machinery already introduced.    
 
\begin{proposition} \label{kirbyr} Let $\mathbb{M} $ be the $\q$-adic filtration on $M.$  Let $p$ be an integer such that $\q M \subseteq \m^p M$. Then
$$ e_1(\mathbb{M}) -e_1(\mathbb{N})  \le {{e_0(\mathbb{M}) -p +1}\choose 2}.$$
 \end{proposition} 
 \begin{proof} We proceed by induction on $r=\dim M. $ If $r=1$ the result follows by Proposition  \ref{kirby}. If $r\ge 2, $ as before,  we can find a minimal system of generators  $\{ a_1,\cdots,a_r\}$ of $J$ such that $a_1$ is $\mathbb{N}$-superficial for $J$ and $\{a_1,\cdots,a_r\}$ is a sequence of $\mathbb{M}$-superficial  elements 
for $\q.$ The module $M/a_1M$ has dimension $r-1$ and $\mathbb{M}/a_1M$ is a good $\q$-filtration on it. 
Now $ e_1(\mathbb{M})-e_1(\mathbb{N}) =e_1(\mathbb{M}/ a_1M)-e_1(\mathbb{N}/a_1M), $ $e_0(\mathbb{M})=
e_0(\mathbb{M}/a_1M ) $ and $\q M/a_1M  \subseteq \m^p M/a_1M$. Hence the result follows by the inductive assumption. 
\end{proof}
For completeness we recall that,  by using a deeper investigation, for local Cohen-Macaulay rings of embedding dimension $b$, Elias in [E1], Theorem 1.6, proved 
 $$e_1(\m) \le  {{e_0(\m) }\choose 2} - {{b-1 }\choose 2}.$$
 An easier approach was presented by the authors in \cite{RV2}  where the result was  proved    for any $\m$-primary ideal $\q.$   
  
\begin{theorem}\label{ic}  Let $(A,\m)$ be a Cohen-Macaulay
local ring of dimension $r $  and  let $\q$  be  an $\m$-primary ideal in 
$A.$ Then
  $$e_1(\q) \le {e_0(\q)  \choose 2}- {\mu(\q)-r \choose
2}-\la(A/\q)+1.$$
\end{theorem}
 
Notice that in the particular case of an $\m-$primary ideal $\q \subseteq \m^2 $ a nice proof was  produced by Elias in \cite{E3}.

\vskip 3mm

We move now to the higher dimensional case   of Proposition \ref{d=1}.     This  improves  {\bf{Northcott's inequalty}} 
to filtrations  of  a module which is not necessarily Cohen-Macaulay.

\begin{theorem}\label{MeN} Let $\mathbb{M}=\{M_j\}_{j\ge 0}$ be a good $\q$-filtration of a module $M$ of dimension $r. $ If $s\ge 1$ is an integer and   $J$ is an ideal generated by a maximal  $\mathbb{M}$-superficial sequence  
for $\q,$ then we have:
$$e_1(\mathbb{M})-e_1(\mathbb{N})\ge    s\, e_0(\mathbb{M})-\la (M/M_{s-1})-\la(M/M_s+JM).$$
\end{theorem}
\begin{proof} If $r=1,$   by Proposition \ref{d=1}  we get $$e_1(\mathbb{M})-e_1(\mathbb{N})\ge  s\, e_0(\mathbb{M})-\la (M/M_s).$$ Hence we must prove that $$s e_0(\mathbb{M})-\la (M/M_s)\ge s \, e_0(\mathbb{M})-\la (M/M_{s-1})-\la(M/M_s+JM).$$ This is equivalent to proving  $$\la (M/M_{s-1})\ge \la (M_s+JM/M_s)=\la(JM/JM\cap M_s).$$ Since $J=(a)$ is a principal ideal and $aM_{s-1}\subseteq aM\cap M_s,$ we have a surjection $$M/M_{s-1}\overset {a} \to aM/aM\cap M_s$$ and the conclusion follows in this case.

Let $r\ge 2;$ by using the above remark, we can find a minimal system of generators  $\{a_1,\cdots,a_r\}$ of $J$ such that $a_1$ is $\mathbb{N}$-superficial for $J$ and $\{a_1,\cdots,a_r\}$ is a sequence of $\mathbb{M}$-superficial  elements 
for $\q.$
 
 The module $M/a_1M$ has dimension $r-1$ and $\mathbb{M}/a_1M$ is a good $\q$-~filtration on it. Furthermore  it is clear that  $a_2,\cdots , a_r $  is a maximal sequence of $\mathbb{M}/a_1M$-superficial  elements 
for $\q.$ Hence, if we let $K$ be the ideal generated by $a_2,\cdots , a_r $ and $\mathbb{K}$ the $K$-adic filtration on  $M/a_1M$, then 
 by  induction, and after a little standard computation,  we get 
\begin{equation*}
e_1(\mathbb{M}/a_1M)-e_1(\mathbb{K})\ge   s\, e_0(\mathbb{M}/a_1M)-\la (M/M_{s-1}+a_1M)-\la(M/M_s+JM).
\end{equation*} Since $e_0(\mathbb{M}/a_1M)=e_0(\mathbb{M}),$  $\mathbb{N}/a_1M=\mathbb{K}$ and 
$$e_1(\mathbb{M})-e_1(\mathbb{M}/aM)=e_1(\mathbb{N})-e_1(\mathbb{N}/aM),$$ we finally  get 
\begin{equation*}
\begin{split} 
e_1(\mathbb{M})-e_1(\mathbb{N})& \ge   s\, e_0(\mathbb{M})-\la (M/M_{s-1}+a_1M)-\la(M/M_s+JM) \\ 
&\ge s\, e_0(\mathbb{M})-\la (M/M_{s-1})-\la(M/M_s+JM)
\end{split} 
\end{equation*} 
which gives the conclusion.
 \end{proof}

\vskip 2mm
\begin{remark} {\rm {Let us apply our theorem to  the very particular case when $M_j=\q^j$ for every $j\ge 0$    and  $\q$ is a primary ideal of $A.$  It is well known that if $J$ is any minimal reduction of $\q,$ then we have an injection $$J/J\m\to \q/\q\m$$ which proves that any minimal system of generators of $J$ is part of a minimal system of generators of $\q$. Furthermore  $J$ can be minimally generated by a maximal sequence of $\mathbb{M}$-superficial  elements for $\q.$ Hence we can apply the above Theorem to  get: 
\vskip 2mm 
$\bullet \ \ \ s=1$ $$e_1(\q)-e_1(J)\ge e_0(\q)-\la(A/\q)$$ which is exactly Theorem 3.1 in \cite{GN}.
\vskip 2mm
$\bullet \ \ \ s=2$

$$e_1(\q)-e_1(J)\ge 2 e_0(\q)-\la(A/\q)-\la(A/\q^2+J).$$  This means $$2 e_0(\q)-e_1(\q)+e_1(J)\le 
2\la(A/\q)+\la(\q/\q^2+J).$$
Since $r=\la(J/J\m)$, if we let $t:=\la (\q/\q\m)-r$, we can find element   $x_1,\cdots,x_t\in \q$ such that $\q=J+(x_1,\cdots,x_t).$ Hence the  canonical map $$\varphi:(A/\q)^t\to \q/(\q^2+J)$$ given by  $\varphi(\overline{a_1},\cdots,\overline{a_t})=\sum\overline{a_ix_i}$ is surjective and we get 
$$2 e_0(\q)-e_1(\q)+e_1(J)\le(t+2)\la(A/\q)$$ which is Proposition 3.7 in \cite{C}.
\vskip 2mm

When  $M=A$ is Cohen-Macaulay and $\mathbb{M}$ is a Hilbert filtration on $A$, which means that $M_j=I_j$ with  $I_j$  ideals in $A,$ $I_0=A,$ $I_1$ is $\m$-primary and $\mathbb{M}$ is $I_1$-good, Guerrieri and Rossi proved in \cite{GR1} the following formula:
$$e_1(\mathbb{M})\ge 2e_0(\mathbb{M})-\big (\la(I_1/I_2)-r\la(A/I_1)+\la((I_2\cap J)/JI_1)+2\la(A/I_1)\big).$$ If we apply the above theorem in this situation, we get
$$e_1(\mathbb{M})-e_1(\mathbb{N})\ge 2e_0(\mathbb{M})-\la(A/I_1)-\la(A/I_2+J).$$

Since $A$ is Cohen-Macaulay, every superficial sequence is a regular sequence in $A$ and thus  $e_1(\mathbb{N})=0$ and $r\la(A/I_1)=\la (J/JI_1).$ Then, by easy computation, we can see that the two bounds coincide.}} \end{remark}
\vskip 4mm
 We now want   to extend  Proposition \ref{EeN} to the higher dimensional case. We recall that
given  a good $\q$-filtration $\mathbb{M}$  of the $r$-dimensional module $M,$ we can consider the ideal $J$ generated by a   maximal   $\mathbb{M}$-superficial sequence   for $\q,$ and we are interested in the study of two related  filtrations on $M$:  the $J$-adic filtration $\mathbb{N}:=\{J^jM\}$ already defined 
  and the filtration $\mathbb{E}$ given by 
 $$\mathbb{E}:\ \ \ \ M\supseteq M_1\supseteq JM_1 \supseteq J^2M_1 \supseteq \cdots \supseteq J^jM_1 \supseteq J^{j+1}M_1\supseteq \cdots $$ In the following for any  good $\q$-filtration $\mathbb{M}$ and for every  ideal $J$ generated by a   maximal   $\mathbb{M}$-superficial sequence   for $\q,$ we may associate the two good $J$-filtrations $\mathbb{N} $ and $\mathbb{E}.  $

 As in the remark before Theorem \ref{MeN}, by
 using Proposition \ref{ei},  we can easily see that $$e_1(\mathbb{E})-e_1(\mathbb{E}/aM)=e_1(\mathbb{N})-e_1(\mathbb{N}/aM).$$
 Notice that in general $e_1(\mathbb{N}) $ and $e_1(\mathbb{E}) $  depend on $J $  (see \cite{V1}, \cite{GO1}, \cite{GHV}, \cite{GGHOV}).

\begin{proposition}\label{E} Let $\mathbb{M}$ be a good $\q$-filtration of   a module $M$ of dimension $r $  and let $J$ be an ideal generated by a   maximal   $\mathbb{M}$-superficial sequence   for $\q. $ Then  we have 

$$e_1(\mathbb{E})-e_1(\mathbb{N})\ge   e_0(\mathbb{M})-h_0 (\mathbb{M})+\la(M_1+H^0(M)/M_1).$$
\end{proposition}
\begin{proof} If $r=1$ we apply Proposition \ref{EeN}. Let $r\ge 2;$  as before  we can find an element $a\in J$  which is superficial for $\mathbb{N}$ and $\mathbb{E}.$ Then we have {\small{ \begin{multline*}e_1(\mathbb{E})-e_1(\mathbb{N})=e_1(\mathbb{E}/aM)-e_1(\mathbb{N}/aM))\\ 
\ge e_0(\mathbb{M}/aM)-h_0 (\mathbb{M}/aM)+\la\left ((M_1+aM/aM)+H^0(M/aM))/(M_1+aM/aM)\right )\\
=e_0(\mathbb{M})-h_0 (\mathbb{M})+\la  \left (aM:_M\m^n+M_1/M_1\right)\ge e_0(\mathbb{M})-h_0 (\mathbb{M})+\la(M_1+H^0(M)/M_1).\end{multline*}}}
\end{proof}

\section{The border cases} 

The aim of  this section is the study of the  extremal cases with respect to the inequalities proved in the above section. With few exceptions we assume $M$ is a filtered Cohen-Macaulay module.
\vskip 2mm

 The following result, first   proved in the classical case by S. Huckaba in \cite{H1}, has been reconsidered  and extended in a series of recent papers (see \cite{HM2}, \cite{VP2}, \cite{JSV}, \cite{Co}), where, unfortunately,  the original heavy homological background was still essential. Recently   Verma in an expository paper on the Hilbert coefficients presents results  proved  by   Huckaba and  Marley \cite{HM2} by using our approach.  We remark that the statements involving the Hilbert coefficients $e_j$ with $j\ge 2$ are new even in the classical case, except  for the bound
on $e_2$ which had been already proved in \cite{CPR}.

We recall that,  given a good $\q$-filtration $\mathbb{M}$ of the module $M$  and     an ideal $J$ generated by a maximal sequence of $\mathbb{M}$-superficial  elements 
for $\q,$ we denote by $v_j(\mathbb{M})$ the non negative integers
$$v_j(\mathbb{M}):=\la(M_{j+1}/JM_j).$$ In Lemma \ref{eu} we proved that, if $M$ is one-dimensional and Cohen-Macaulay, then for every $j\ge 0$ 
$$e_i(\mathbb{M})=\sum_{j\ge i-1}\binom{j}{i-1}v_j(\mathbb{M}),$$ while, in Proposition \ref{Vi}, we proved that if $M$ is Cohen-Macaulay then 
$$e_1(\mathbb{M})\le \sum_{j\ge 0}v_j(\mathbb{M}).$$

 \begin{theorem}\label {e1=CM} Let $\mathbb{M}$ be a good $\q$-filtration of   a module $M$ of dimension $r $  and let $J$ be an ideal generated by a   maximal   $\mathbb{M}$-superficial sequence   for $\q. $ Then we have 
\vskip 2mm \noindent 
a)  $e_1(\mathbb{M})\le \sum_{j\ge 0}v_j(\mathbb{M})$ 
\vskip 2mm  \noindent 
b) $e_2(\mathbb{M})\le \sum_{j\ge 0}jv_j(\mathbb{M}).$
\vskip 2mm \noindent 
c) The following conditions are equivalent
\vskip 2mm
 
1. $\depth \ gr_{\mathbb{M}}(M)\ge r-1.$

\vskip 1mm
 
2. $e_i(\mathbb{M})=\sum_{j\ge i-1}\binom{j}{i-1}v_j(\mathbb{M})$ for every $i\ge 1.$

\vskip 1mm 
 
3. $e_1(\mathbb{M})=\sum_{j\ge 0}v_j(\mathbb{M}).$

\vskip 1mm 
 
 4. $e_2(\mathbb{M})=\sum_{j\ge 0}jv_j(\mathbb{M}).$

 \vskip 1mm 
 
 5. $P_{\mathbb{M}}(z)= {{\la(M/M_1) + \sum_{j\ge 0} (v_j(\mathbb{M}) - v_{j+1}(\mathbb{M}) ) z^{j+1}}\over {(1-z)^r}}$
 
    \end{theorem}
   \begin{proof} Let $J=(a_1,\cdots,a_r)$ and $\frak{a}=(a_1,\cdots,a_{r-1});$ we first remark that, by Lemma \ref{sup} and Theorem \ref{vv}, $\text{depth} \ gr_{\mathbb{M}}(M)\ge r-1$ if and only if $M_{j+1}\cap \frak{a}M  = \frak{a}M_j$ for every $j\ge 0.$ Further  we have $$v_j(\mathbb{M})=v_j(\mathbb{M}/\frak{a}M)+\la(M_{j+1}\cap \frak{a}M+JM_j/JM_j),$$ hence $v_j(\mathbb{M})\ge v_j(\mathbb{M}/\frak{a}M)$ and equality holds if and only if 
   $M_{j+1}\cap \frak{a}M   \subseteq JM_j.$ This  is certainly the case when $M_{j+1}\cap \frak{a}M=\frak{a}M_j. $ \par \noindent Hence, if $\text{depth} \ gr_{\mathbb{M}}(M)  \ge~ r-1,$ then $v_j(\mathbb{M})=v_j(\mathbb{M}/\frak{a}M)$ for every $j\ge 0.$ By induction on $j,$ we can prove that the converse holds. Namely $M_1 \cap \frak{a}M=\frak{a}M$ and, if $j\ge 1,$  then we have
   \begin{equation*}\begin{split} M_{j+1}\cap \frak{a}M & \subseteq JM_j\cap \frak{a}M=(\frak{a}M_j+a_rM_j)\cap \frak{a}M\\
   &=\frak{a}M_j+(a_rM_j\cap \frak{a}M)\subseteq \frak{a}M_j+a_r(M_j\cap \frak{a}M)\\
   &=\frak{a}M_j+a_r\frak{a}M_{j-1}=\frak{a}M_j.
   \end{split}
   \end{equation*}  where $a_rM_j\cap \frak{a}M\subseteq a_r(M_j\cap \frak{a}M)$ because $a_r$ is regular modulo $\frak{a}M,$ while  $M_j\cap \frak{a}M=\frak{a}M_{j-1}$ follows by induction.
   
   Since $M/\frak{a}M$ is Cohen-Macaulay of dimension one, we get
   $$e_1(\mathbb{M})=e_1(\mathbb{M}/\frak{a}M)=\sum_{j\ge 0}v_j(\mathbb{M}/\frak{a}M)\le \sum_{j\ge 0}v_j(\mathbb{M}).$$ Equality holds if and only if   $\text{depth} \ gr_{\mathbb{M}}(M)\ge r-1.$ This proves a) once more and moreover gives the equivalence between 1. and 3. in c). By using (\ref{PZ})  and Proposition \ref{ei} this also gives the equivalence between  1. and 5.  in c).
   
   Now, if $\frak{b}$ is the ideal generated by $a_1,\cdots,a_{r-2},$ then, as before, we get
   $$e_2(\mathbb{M})=e_2(\mathbb{M}/\frak{b}M)\le e_2(\mathbb{M}/\frak{a}M) =\sum_{j\ge 1}jv_j(\mathbb{M}/\frak{a}M)\le \sum_{j\ge 1}jv_j(\mathbb{M}).$$ This proves b) and 4 $\Longrightarrow$1. To complete the proof of the Theorem, we need only to show that
   1 $\Longrightarrow$ 2. If $\text{depth} \ gr_{\mathbb{M}}(M)\ge r-1,$ then $\mathbb{M}$ and $\mathbb{M}/\frak{a}M$ have the same $h$-polynomial; this implies that for every $i\ge  1$ we have $$e_i(\mathbb{M})=e_i(\mathbb{M}/\frak{a}M)=\sum_{j\ge i-1}\binom{j}{i-1}v_j(\mathbb{M}/\frak{a}M)=\sum_{j\ge i-1}\binom{j}{i-1}v_j(\mathbb{M}).$$
   \end{proof}

\vskip 2mm
In the above result the equality in b) does not force $gr_{\mathbb{M}}(M)$ to be  Cohen-Macaulay.  In fact we will see later that in a two dimensional local Cohen-Macaulay ring, $e_2$ can be zero,  but 
depth $gr_{\mathbb{M}}(M) = 0 $ (see Example \ref{1}). 
\vskip 3mm

We recall that  Proposition \ref{Vi}  extends  Huckaba's inequality without assuming the Cohen-Macaulayness of $M.$  In particular we proved that 
$$e_1(\mathbb{M})-e_1(\mathbb{N})\le \sum_{j\ge 0}v_j(\mathbb{M}) $$
where $\mathbb{N} $ is the   $J$-adic filtration on $M.$  

 If we do not assume that $M$ is Cohen-Macaulay, we are able to handle the equality  only for the $\m$-adic filtration on $A.$
 Surprisingly   in \cite{RV5}, Theorem 2.13,  the authors proved that  the equality in Proposition \ref{Vi} forces the ring  $A$  itself to be Cohen-Macaulay and hence, by Theorem \ref{e1=CM}, $gr_{\m}(A)$ to have almost maximal depth. The result is the following.

\begin{theorem}\label{th2} Let $(A,\m)$ be a local ring of dimension $r \ge 1$  and let $J$ be the ideal generated by a maximal  $\m$-superficial sequence. The following conditions are equivalent:

\vskip 2mm 
1. $e_1(\m)-e_1(J)=\sum_{j\ge 0}v_j(\m).$

\vskip 2mm
2. $A$ is Cohen-Macaulay and $ \depth  gr_{\m}(A)\ge r-1.$
\end{theorem}
\begin{proof} If $A$ is Cohen-Macaulay, then $e_1(J)=0$ and,  by the above result, we find  that 2) implies 1).

\noindent We prove now that 1) implies 2) by induction on $r.$ If $r=1,$ the result follows by Proposition \ref{e1<=}  since $W\subseteq \m.$ Let $r\ge 2$; by Lemma \ref{key} we can find a minimal basis $\{a_1,\dots,a_r\}$ of $J$ such that $a_1$ is $J$-superficial,  $\{a_1,\dots,a_r\}$ is an $\m$-superficial sequence  and $e_1(\m)-e_1(J)=e_1(\m/(a_1))-e_1(J/(a_1)).$
 Now $A/(a_1)$ is a local ring of dimension $d-1$ and  $J/(a_1)$ is generated by a maximal $\m/(a_1)$-superficial sequence.
 We can then apply  Proposition  \ref{Vi} to get  $ \ \sum_{j\ge 0}v_j(\m)=e_1(\m)-e_1(J)=e_1(\m/(a_1))-e_1(J/(a_1))\le \sum_{j\ge 0}v_j(\m/(a_1))\le \sum_{j\ge 0}v_j(\m).$ 
 
 This implies $$e_1(\m/(a_1))-e_1(J/(a_1))= \sum_{j\ge 0}v_j(\m/(a_1))$$  which, by the inductive assumption,  implies that $A/(a_1)$ is Cohen-Macaulay. By Lemma \ref{nuovo}  $A$ is Cohen-Macaulay so that $e_1(J)=0$ and then $e_1(\m)=\sum_{j\ge 0}v_j(\m);$ this implies  $\depth  gr_{\m}(A)\ge r -1$ and the result  is proved.
\end{proof}

 \vskip 2mm
\begin{remark} {\rm{As the reader can see, the   above result  had been presented for the $\m$-adic filtration. Actually a more general statement holds. We need a filtration $ \mathbb{M} $ such that $W   \subseteq M_1 $    in order to apply Proposition \ref{e1<=}   going down of dimension by superficial sequences.   }}
\end{remark}
\vskip 2mm

As a trivial consequence of the above result we have the following
\begin{corollary} Let $(A,\m)$ be a local ring of dimension $r\ge 1 $  and let $J$ be the ideal generated by a maximal  $\m$-superficial sequence. If $e_1(J)\le 0,$ then $$e_1(\m)\le \sum_{j\ge 0} v_j(\m).$$  Moreover, the following conditions are equivalent:

\vskip 2mm \noindent
1. $e_1(\m)=\sum_{j\ge 0} v_j(\m).$

\vskip 2mm \noindent 2. $A$ is Cohen-Macaulay and
$\depth  gr_{\m}(A)\ge r-1.$ 
\end{corollary}

 \begin{remark} {\rm{ Notice that the condition $ e_1(J) \le 0 $ is satisfied for example if $A$ is Buchsbaum
(see \cite{SV}, Proposition 2.7) or if $ \depth  A \ge r-1.  $ In fact if $ \depth  A \ge r-1 $ and   $ a_1, \dots, a_{r-1} $ is a superficial sequence in $J$, then $ e_1( J) =  e_1(J/(a_1, \dots, a_{r-1}) ) $ and $ 
 e_1(J/(a_1, \dots, a_{r-1}) ) \le 0 $ by Lemma \ref{W}. 
 
\noindent Hence  the previous result is an interesting extension of Huckaba-Marley's  result where  the Cohen-Macaulayness of $A $ is assumed.
 \vskip 2mm
   If $A$  is an unmixed,
equidimensional local ring that is a homomorphic image of a Cohen-Macaulay local
ring, then Vasconcelos conjectured that for any ideal $ J $ generated by a system of parameters, the Chern  coefficient  $ e_1(J) < 0 $ is equivalent to $A$  being non Cohen-Macaulay.

We remark that if  $\dim A= 1$, the property $e_1(J) = 0 $  is characteristic of
the Cohen-Macaulayness.
For $ \dim A \ge 2, $  the situation is somewhat different.  Consider the non  Cohen-Macaulay  ring
$ A = k[x, y, z]/(z(x,y, z))  $ and $J=(x,y)A. $  Despite the lack of the  Cohen-Macaulyness, 
 we have $e_1(J) = e_1(S) = 0  $ where $ S = k[x, y]\simeq  A/H^0_{\m}(A).  $  
  Recently a  very nice result was proved  by Ghezzi, Hong, Vasconcelos   who established the conjecture if $A$  is a homomorphic image of a Gorenstein ring, and for all universally catenary integral domains containing  fields, see  \cite{GHV}.  A remarkable   extension to the reduced case was  obtained in \cite{GGHOV}.  }}
\end{remark}
\vskip 4mm
In the Cohen-Macaulay case we describe now another set of numerical characters of the filtered module $M,$ which are important in  the study of the Hilbert coefficients. 

Let $\mathbb{M}=\{M_j\}_{j\ge 0}$ be a good $\q$-filtration of $M$ and    $J$  an ideal generated by   an $\mathbb{M}$-superficial sequence  
for $\q;$ then, for every $j\ge 0,$ we let  
\begin{equation}w_j(\mathbb{M}):=\la(M_{j+1}+JM/JM)=\la(M_{j+1}/M_{j+1}\cap JM).\end{equation}
Since $$JM_j\subseteq JM\cap M_{j+1}\subseteq M_{j+1},$$
we get \begin{equation}
v_j(\mathbb{M})=w_j(\mathbb{M})+\la(JM\cap M_{j+1}/JM_j)\end{equation}

The length of the abelian group $JM\cap M_{j+1}/JM_j$ will be denoted by $vv_j(\mathbb{M})$ since these groups  are the homogeneous components of the   Valabrega-Valla module $$VV(\mathbb{M}):=\bigoplus_{j\ge 0} (JM\cap M_{j+1}/JM_j)$$ of $\mathbb{M}$ with respect to $J$, as defined in \cite{VW2}, chapter 5. For example one has $$vv_0(\mathbb{M})=0, \ \ \ \
vv_1(\mathbb{M})=\la(JM\cap M_2/JM_1)$$ and since $M_{j+1}=JM_j$ for large $j,$ $vv_j(\mathbb{M})=0$ for $j \gg 0.$  It follows that, in the case of the $\m$-adic filtration $\{M_j=\m^j\}$ of the Cohen-Macaulay ring $A,$ one has 
$vv_1(\mathbb{M})=0$ by the analytic independence of a maximal regular sequence.

The significance of $VV(\mathbb{M})$ lies in the fact that, if $M$ is Cohen-Macaulay, by Valabrega-Valla, $gr_{\mathbb{M}}(M)$ is Cohen-Macaulay if and only if $VV(\mathbb{M})=0.$

As a consequence we have that $gr_{\mathbb{M}}(M)$ is Cohen-Macaulay if and only if $v_j(\mathbb{M})=w_j(\mathbb{M})$ for every $j\ge 0.$

In the one-dimensional Cohen-Macaulay case, we have $$e_1(\mathbb{M})=\sum_{j\ge  0}v_j(\mathbb{M})=\sum_{j\ge 0}w_j(\mathbb{M})+\sum_{j\ge 0}vv_j(\mathbb{M}),$$ so that $\sum_{j\ge 0}w_j(\mathbb{M})\le e_1(\mathbb{M})$ and equality holds if and only if 
$gr_{\mathbb{M}}(M)$ is Cohen-Macaulay.

This result can be extended to higher dimensions; we first need to remark that the integers $w_j$ do not change upon reduction  a superficial sequence. Namely  if $a\in J,$  then $$w_j(\mathbb{M}/aM)=\la \big ((M_{j+1}+aM/aM)+(JM/aM)\big /(JM/aM)\big )
=$$ $$=\la \big ( ( M_{j+1}+JM/aM)/(J(M/aM) \big )= \la(M_{j+1}+JM/JM)
=w_j(\mathbb{M}).$$

\vskip 2mm
\begin{theorem}\label {e1=wj} Let $\mathbb{M} $ be a good $\q$-filtration of the Cohen-Macaulay module $M$ of dimension $r\ge 1$ and let    $J$ be  an ideal generated by  a maximal  $\mathbb{M}$-superficial sequence  
for $\q. $ Then we have 
\vskip 2mm 
a)  $e_1(\mathbb{M})\ge \sum_{j\ge 0}w_j(\mathbb{M})$ 
\vskip 2mm 
b) $e_1(\mathbb{M})= \sum_{j\ge 0}w_j(\mathbb{M})$ if and only if $gr_{\mathbb{M}}(M)$ is Cohen-Macaulay.
\end{theorem}
\begin{proof} We may assume $r\ge 2$ and we let $J=(a_1,\cdots,a_r). $ Denote  $\frak{a}=(a_1,\cdots,a_{r-1}), $ then we have $$e_1(\mathbb{M})=e_1(\mathbb{M}/\frak{a}M)=\sum_{j\ge 0}v_j(\mathbb{M}/\frak{a}M)=\sum_{j\ge 0}w_j(\mathbb{M}/\frak{a}M)+\sum_{j\ge 0}vv_j(\mathbb{M}/\frak{a}M)=$$
$$=\sum_{j\ge 0}w_j(\mathbb{M})+\sum_{j\ge 0}vv_j(\mathbb{M}/\frak{a}M).$$
This proves a) and  also implies that $e_1(\mathbb{M})= \sum_{j\ge 0}w_j(\mathbb{M})$ if and only if $VV(\mathbb{M}/\frak{a}M)=0,$ if and only if $gr_{\mathbb{M}/\frak{a}M}(M/\frak{a}M)$ is Cohen-Macaulay. Hence b)  follows by Sally's machine.  

\end{proof}

\vskip 3mm

As a corollary of Theorem \ref {e1=wj},  we get the following achievement   which extends to our general  setting the main result  of A. Guerrieri in \cite{Gu1}. Here the proof is quite simple and is due to Cortadellas (see \cite{Co}).

\begin{corollary}  \label{GU} Let $\mathbb{M} $ be a good $\q$-filtration of the Cohen-Macaulay module $M$ of dimension $r\ge 1$ and  let  $J$  be an ideal generated by a maximal  $\mathbb{M}$-superficial sequence  
for $\q. $  If $\la(VV(\mathbb{M}))=1, $ then $\depth\  gr_{\mathbb{M}}(M)=r-1.$
\end{corollary}
\begin{proof} Since $\la(VV(\mathbb{M}))=1, $ then  $$\sum_{j\ge 0}w_j(\mathbb{M})=\sum_{j\ge 0}v_j(\mathbb{M})-1.$$ Since for some $n\ge 0$ we have $M_{n+1}\cap JM\not= JM_n$, we  have $\depth gr_{\mathbb{M}}(M)\le r-1,$ hence $e_1(\mathbb{M}) > \sum_{j\ge 0}w_j(\mathbb{M})$. We  get $$\sum_{j\ge 0}v_j(\mathbb{M})\ge e_1(\mathbb{M})> \sum_{j\ge 0}w_j(\mathbb{M})=\sum_{j\ge 0}v_j(\mathbb{M})-1.$$ This implies $\sum_{j\ge 0}v_j(\mathbb{M})= e_1(\mathbb{M})$ and the conclusion follows.\end{proof}

In the classical case, it was proved by Guerrieri in \cite{Gu2} that if $\la (\q^2\cap J/\q J)=2=\la(VV(\mathbb{M}))$, then $\depth  gr_{\mathbb{M}}(M)\ge r-2,$ a result which was extended by Wang in  \cite{W2}, where he proved that if $\la(VV(\mathbb{M}))=2$ then  $\depth \  gr_{\mathbb{M}}(M)\ge r-2.$   If  $ \la (\q^2\cap J/\q J)= \la(VV(\mathbb{M}))=3,$   Guerrieri and Rossi proved that $\depth  gr_{\mathbb{M}}(M)\ge r-3,$  provided $A$ is Gorenstein.  A conjecture of C. Huneke predicts that if $vv_j\le 1$ for every $j,$ then 
 $\depth \  gr_{\mathbb{M}}(M)\ge r-1.$ This is not true as shown by Wang. However Colom\'e and Elias proved that the condition $vv_j\le 1$ for every $j,$ implies $\depth \  gr_{\mathbb{M}}(M)\ge r-2.$ Concerning  this topic see also \cite{Gu1}, \cite{GR2}, \cite{W2}, \cite{Po}, \cite{CE}.

\vskip 2mm We want now to study the extremal case in Northcott's inequlity. First we need to recall a lower bound for  $e_1$ which was proved by Elias and Valla in \cite{EV}.
Given  a Cohen-Macaulay local ring $(A,\m),$ one has  \begin{equation}\label{ev}
e_1(\m)\ge 2e_0(\m)-h-2,\end{equation}
where $h=\mu(\m)-\dim(A)$ is the embedding codimension of $A.$ 
Equality holds above  if and only if  the $h$-polynomial  is short enough.
\vskip 2mm 
In our general setting we have the inequality given by Theorem \ref{MeN}, namely
$$e_1(\mathbb{M})-e_1(\mathbb{N})\ge    s\, e_0(\mathbb{M})-\la (M/M_{s-1})-\la(M/M_s+JM).$$
When  $M$ is Cohen-Macaulay, we have $e_1(\mathbb{N})=0$ so that, if $s=2,$ we get
$$e_1(\mathbb{M})\ge    2\, e_0(\mathbb{M})-\la (M/M_{1})-\la(M/M_2+JM).$$If we let 
\begin{equation} \label{h} h(\mathbb{M}):=\la(M_1/JM+M_2), \end{equation}  then we have 
\begin{equation}\label{8}
\begin{split} 
e_1(\mathbb{M}) &\ge 2\, e_0(\mathbb{M})-\la (M/M_{1})-\la(M/M_2+JM)
\\ &=2\, e_0(\mathbb{M})-h(\mathbb{M}) -2\, h_0(\mathbb{M}) 
\\
\end{split} 
\end{equation} a formula which extends (\ref{ev}) because $\la(\m/J+\m^2)=\mu(\m)-\dim(A)=h.$

\vskip 3mm

\vskip 2mm We remark that the integer $h(\mathbb{M})$ coincides with the embedding codimension in the case of the $\m$-adic filtration. Further we have \begin{equation}\label{h=}  h(\mathbb{M})=h(\mathbb{M}/JM)=h_1(\mathbb{M}/JM) \end{equation}  and also \begin{equation}\label{h=h1+} h(\mathbb{M})=h_1(\mathbb{M})+\la (M_2\cap JM/JM_1).\end{equation}

The proof of the following theorem is exactly the same as the original given in \cite{EV} and \cite{GR1}. We reproduce it here because is a  typical    example of  the  strategy to reduce dimension by using superficial elements and the Sally machine. 

We recall that, when  $M$ is Cohen-Macaulay and $\dim (M) =1, $ we   introduced the reduction number of the filtration  $\mathbb{M} $ as the integer  $s(\mathbb{M})= \min \{j : H_{\mathbb{M}} (j) =e_0(\mathbb{M}) \} $ and it turns out that it is also the degree of the $h$-polynomial of $M., $ see (\ref{reduction}).
\vskip 2mm Accordingly with this case, in 
the following $s(\mathbb{M})$ will denote the {\bf degree of the $h$-polynomial of $M$} (see Section 1.3. for the definition).

\begin{theorem}\label{elv} Let $\mathbb{M} $ be a good $\q$-filtration of a Cohen-Macaulay module $M$ and   let  $J$  be an ideal generated by   a maximal $\mathbb{M}$-superficial sequence  
for $\q.$  The following conditions are equivalent:
\vskip 2mm 
a) $e_1(\mathbb{M}) =2\, e_0(\mathbb{M})-2\, h_0(\mathbb{M})-h(\mathbb{M})$
\vskip 2mm 
b) $s(\mathbb{M})\le 2$ and $M_2\cap JM=JM_1.$
\vskip 2mm \noindent 
If either of the above conditions holds, then $gr_{\mathbb{M}}(M)$ is Cohen-Macaulay.
\end{theorem}
\begin{proof} We prove that b) implies a). Since $M_2\cap JM=JM_1$ and $M$ is Cohen-Macaulay,  we get $h(\mathbb{M})=h_1(\mathbb{M}).$  Since $s(\mathbb{M})\le 2$, we have 
 $$e_0(\mathbb{M})=h_0(\mathbb{M})+h_1(\mathbb{M})+h_2(\mathbb{M}) \  \  \  \ \ \ \ \ e_1(\mathbb{M})=h_1(\mathbb{M})+2\, h_2(\mathbb{M}).$$ Hence 
 \begin{equation*}
\begin{split} 2\, e_0(\mathbb{M})-2\, h_0(\mathbb{M})-h(\mathbb{M})& =2(h_0(\mathbb{M})+h_1(\mathbb{M})+h_2(\mathbb{M}))-2 h_0(\mathbb{M})-h_1(\mathbb{M})\\ 
&=h_1(\mathbb{M})+2 h_2(\mathbb{M})
=e_1(\mathbb{M}).
\\
\end{split} 
\end{equation*}
Let us prove the converse by induction on $r:=\dim(M).$ If $r=0,$ then we have $P_{\mathbb{M}}(z)=\sum_{i=0}^{s}h_i(\mathbb{M})$ with $h_i(\mathbb{M})\ge 0$ and
where we let $s:=s(\mathbb{M}).$ Since  $h(\mathbb{M})=h_1(\mathbb{M}),$ it is clear that $e_1(\mathbb{M}) \ge 2\, e_0(\mathbb{M})-2\, h_0(\mathbb{M})-h(\mathbb{M})$ 
and  if we have equality, then  $$ h_3(\mathbb{M})+2 h_4(\mathbb{M})+\cdots +(s-2)h_s(\mathbb{M})=0,$$ which implies $s\le 2.$

If $r\ge 1$,  let $J=(a_1,\cdots,a_r)$ and $K:=(a_1,\cdots,a_{r-1}).$ Then we have
\begin{equation*}
\begin{split}
2\, e_0(\mathbb{M})-2\, h_0(\mathbb{M})-h(\mathbb{M})&=e_1(\mathbb{M})=e_1(\mathbb{M}/KM)\ge e_1(\mathbb{M}/JM)\\
&\ge 
2\, e_0(\mathbb{M}/JM)-2\, h_0(\mathbb{M}/JM)-h(\mathbb{M}/JM)\\
&=2\, e_0(\mathbb{M})-2\, h_0(\mathbb{M})-h(\mathbb{M})
\\
\end{split}
\end{equation*} where we used several times Proposition \ref{ei}.
This gives $$e_1(\mathbb{M}/KM)= e_1(\mathbb{M}/JM)$$ which, again by Proposition \ref{ei}, implies $\text{depth} \  gr_{\mathbb{M}}(M/KM)=1.$ By Sally's  machine,  $gr_{\mathbb{M}}(M)$ is Cohen-Macaulay so that $s(\mathbb{M})=s(\mathbb{M}/JM)\le 2.$ Finally,  by  Valabrega and Valla,  $M_2\cap JM=JM_1,$ as wanted.
\end{proof}

We collect in the following formula some of the results we proved   in the case $M$ is Cohen-Macaulay and $J$  an ideal generated by a maximal sequence of $\mathbb{M}$-superficial  elements  for $\q.$ \begin{equation}\label{north}\begin{split} e_1(\mathbb{M})&\ge 2 e_0(\mathbb{M})-h(\mathbb{M}) -2 h_0(\mathbb{M}) \\
&=e_0(\mathbb{M})- h_0(\mathbb{M})+\la(JM+M_2/JM)\\
&= h_1(\mathbb{M})+\la(M_2/JM_1)+\la(M_2/JM\cap M_2).\\
\end{split}\end{equation} Here  the first inequality comes from Theorem \ref{MeN} with $s=1,$  the first equality comes from the identities $$e_0(\mathbb{M})=\la(M/JM),\ \ \  h_0(\mathbb{M})=\la(M/M_1),\ \ \ h(\mathbb{M})=\la(M_1/JM+M_2)$$ and, finally, the last equality  is a consequence of Proposition \ref{ab} which says that \begin{equation}\label{AB}e_0(\mathbb{M})=h_0(\mathbb{M})+ h_1(\mathbb{M})+\la(M_2/JM_1).\end{equation}

We are ready now to study the Hilbert function in the extremal case of Northcott's inequality and, at the same time, in the case of minimal multiplicity.

\begin{theorem}\label{nor} Let $\mathbb{M}=\{M_j\}_{j\ge 0}$ be a good $\q$-filtration of the Cohen-Macaulay module $M$ of dimension $r$ and   let $J$ be  an ideal generated by a maximal  $\mathbb{M}$-superficial sequence  
for $\q. $  Let us consider the following conditions:

\vskip 1mm
1. $s(\mathbb{M})\le 1 $ or, equivalently, $P_{\mathbb{M}}(z)=\frac{h_0(\mathbb{M}) + h_1(\mathbb{M})z}{(1-z)^r}$

\vskip 1mm
2. $e_1(\mathbb{M})=h_1(\mathbb{M}).$

\vskip 1mm 
3. $e_1(\mathbb{M})=e_0(\mathbb{M})-h_0(\mathbb{M}).$

\vskip 1mm 
4. $e_0(\mathbb{M})=h_0(\mathbb{M})+h_1(\mathbb{M}).$

\vskip 1mm
5. $M_2=JM_1.$

\vskip 1mm \noindent Then we have $$ \begin{matrix}
     1 &   \Longrightarrow & 2                  &     \Longrightarrow   & 4                        \\
    \Updownarrow     &                             & \Downarrow &                                 & \Updownarrow   \\
     1 &  \Longleftarrow    & 3                   &                                 & 5 
\end{matrix}$$
If any  of the  first three equivalent conditions holds, then $gr_{\mathbb{M}}(M)$ is Cohen-Macaulay.

\end{theorem}
\begin{proof}  It is clear that $1  \Longrightarrow 2,$  while, by using  (\ref{AB}), we get $4 \Longleftrightarrow 5.$  By looking at (\ref{north}), it is clear that $e_1(\mathbb{M})=h_1(\mathbb{M})$    implies $M_2=JM_1$ and $e_1(\mathbb{M})=e_0(\mathbb{M})-h_0(\mathbb{M})$  so that  $2 \Longrightarrow  3\ \text{and}  \  4.$   We need only to prove that $3 \Longrightarrow 1.$ 

If $e_1(\mathbb{M})=e_0(\mathbb{M})-h_0(\mathbb{M}),$ then $M_2\subseteq JM$, which implies $h_2(\mathbb{M}/JM)=0,$ and equality holds in (\ref{north}). By Theorem \ref{elv}, we get  $gr_{\mathbb{M}}(M)$ is Cohen-Macaulay and $s(\mathbb{M})\le 2.$ Hence $s(\mathbb{M})=s(\mathbb{M}/JM)\le 2;$ but we have seen that $h_2(\mathbb{M}/JM)=0,$ hence  $s(\mathbb{M})=s(\mathbb{M}/JM)\le 1,$ as required.\end{proof}

Notice that the example given after Corollary \ref{cor} shows that in the above theorem the condition $M_2=JM_1$ does not imply $s(\mathbb{M})\le 1.$ In order to have this  implication, we need to put some restriction on the filtration.

\vskip 3mm
 If  $L$ is any submodule of the given  Cohen-Macaulay module $M,$ and $\q$ an $\m$-primary ideal of $A$  such that $\q M\subseteq L,$ let us consider  the filtration \begin{equation}\label{find} \mathbb{M}_L: \{M_0=M,  M_{j+1}=\q^{j}L\}\end{equation}  for every $j\ge 0.$  It is clear that $\mathbb{M}_L$ is a good $\q$-filtration. $ \mathbb{M}_L $  will be called the \textbf{filtration induced} by $L $ and  it will appear in most of the results from now on. \par \noindent   It is important to remark that,  by   definition,   one has $M_{j+1} = \q M_j $ for every $j \ge 1.$ For example this property allows us  to conclude that  the condition $M_2=JM_1$  implies  $s(\mathbb{M})\le 1.$
  
 Let us compare the filtration   $\mathbb{M}_L $ with the already introduced filtration 
  \begin{equation*}  \mathbb{M}^I: \  (  \mathbb{M}^I)_0= M_0=M  \ \ \  \text{and}\ \ \    (  \mathbb{M}^I)_{j+1}=I M_j \ \ \ \text { for every  } j \ge 0. \end{equation*}  
  where $I $ is an ideal containing $\q.$

     Notice that if $\mathbb {M}= \q^nM $ is the $\q$-adic filtration, then $ \mathbb{M}^I = \mathbb{M}_L $  with $L=I M. $   
   
\noindent    When  $ \mathbb{M}  $ is the $\q$-adic filtration on $A$ and $I $ = $\m, $ then  $ \mathbb{M}^{\m}$ coincides with the   filtration 
 \begin{equation*}   \mathbb{M}_{\m}: \   A \supseteq \m \supseteq \m \q \supseteq \m \q^2 \supseteq \dots \supseteq
 \m\q^n \dots  \end{equation*}    already introduced in (\ref{find}).

 \vskip 3mm
 \begin{corollary} \label{ind} Let $M$ be a  given  Cohen-Macaulay module, $L$ a submodule and $\mathbb{M}=\mathbb{M}_L$ the good $\q$-filtration on $M$ induced by $L.$ If  $J$  is an ideal generated by a maximal  $\mathbb{M}$-superficial sequence  
for $\q, $  then all the conditions of the above theorem are equivalent.
 \end{corollary}
 \begin{proof} We prove that, for the filtration $\mathbb{M}$, we have $5  \Longrightarrow 1.$ If  $M_2=JM_1$, then 
 $\q L=JL$ so that $M_{j+1}=JM_j$ for every $j\ge 1.$ By Valabrega-Valla, this implies   $ gr_{\mathbb{M}}(M)$ is Cohen-Macaulay and $s(\mathbb{M}/JM)\le 1.$   Hence $$s(\mathbb{M})= s(\mathbb{M}/JM)\le 1,$$ as wanted.
 \end{proof}

 In the following theorem we study the  equality $e_1(\mathbb{M})=e_0(\mathbb{M})-h_0(\mathbb{M})+1.$ It turns out that we need some extra assumptions  in order to get a complete description of this case. Nevertheless, the theorem extends results of Elias-Valla (see\cite{EV}), Guerrieri-Rossi (see \cite{GR1}),  Itoh (see \cite{I1}), Sally (see \cite{S4}) and Puthenpurakal (see \cite{P2}).

 \begin{theorem}\label{nor+1} Let $M$ be a   Cohen-Macaulay module of dimension $r,$  $L$ a submodule of $M$ and $\mathbb{M}=\mathbb{M}_L$ the good $\q$-filtration on $M$ induced by $L.$ If  $J$ is   an ideal generated by   a maximal  $\mathbb{M}$-superficial sequence  
for $\q $   and we assume $M_2\cap JM=JM_1,$ then the following  conditions  are equivalent.

\vskip 1mm
1. $e_1(\mathbb{M})=e_0(\mathbb{M})-h_0(\mathbb{M})+1.$

\vskip 1mm
2. $e_1(\mathbb{M})=h_1(\mathbb{M})+2.$

\vskip 1mm 
3. $e_0(\mathbb{M})=h_0(\mathbb{M})+h_1(\mathbb{M})+1$ and $gr_{\mathbb{M}}(M)$ is Cohen-Macaulay.

\vskip 1mm 
4. $P_{\mathbb{M}}(z)=\frac{h_0(\mathbb{M})+h_1(\mathbb{M})z+z^2}{(1-z)^r}.$

\end{theorem}
\begin{proof}  It is clear that $4.$ implies $1.$  By (\ref{h=h1+}), the assumption $M_2\cap JM=JM_1$ gives the equality $h_1(\mathbb{M})=h(\mathbb{M}),$ hence, if $e_1(\mathbb{M})=e_0(\mathbb{M})-h_0(\mathbb{M})+1,$ we get by (\ref{north})
\begin{equation} \label{eq} \begin{split} 
e_1(\mathbb{M})&= e_0(\mathbb{M})- h_0(\mathbb{M})+1\\ 
&\ge 2e_0(\mathbb{M})- 2h_0(\mathbb{M})-h_1(\mathbb{M})\\
&=e_0(\mathbb{M})- h_0(\mathbb{M})+\la(M_2/JM_1)\\
&= h_1(\mathbb{M})+2\la(M_2/JM_1).\\
\end{split}
\end{equation}
Now $M_2\not=JM_1$, otherwise, by Corollary \ref{ind}, we have $e_1(\mathbb{M})=e_0(\mathbb{M})-h_0(\mathbb{M}).$  Hence $\la(M_2/JM_1)=1$
and we have equality above, so that $e_1(\mathbb{M})=h_1(\mathbb{M})+2.$ This proves that 1. implies 2. 

Let us assume that $e_1(\mathbb{M})=h_1(\mathbb{M})+2.$ By Corollary \ref{ind}, we have 
$M_2\not=JM_1$, so that equality holds in (\ref{eq})  and  $gr_{\mathbb{M}}(M)$ is Cohen-Macaulay by Theorem \ref{elv}. 

We need only to prove that  $3.$  implies $4.$ Since $gr_{\mathbb{M}}(M)$ is Cohen-Macaulay, $\mathbb{M}$ and $\mathbb{M}/JM$ have the same $h$-polynomial so that $$P_{\mathbb{M}}(z)=\frac{h_{\mathbb{M}/JM}(z)}{(1-z)^r }$$  But we have $e_0(\mathbb{M})=e_0(\mathbb{M}/JM),$ $ h_0(\mathbb{M})=h_0(\mathbb{M}/JM)$ and also
$h_1(\mathbb{M})=h(\mathbb{M})=h_1(\mathbb{M}/JM).$ Under the assumption  $e_0(\mathbb{M})=h_0(\mathbb{M})+h_1(\mathbb{M})+1$ this implies $$\sum_{j\ge 2}h_j(\mathbb{M}/JM)=1.$$ Now if $h_2(\mathbb{M}/JM)=0$, then $\q L\subseteq \q^2 L+JM$ which  implies $h_t(\mathbb{M}/JM)=0$ for every $t\ge 2.$ Hence $h_2(\mathbb{M}/JM)\not=0$ and since $h_j(\mathbb{M}/JM)\ge 0$ for every $j,$ we get 
$h_2(\mathbb{M}/JM)=1$ and $h_j(\mathbb{M}/JM)=0$ for $j\ge 3.$  The proof of the theorem is now complete.
\end{proof}
 
 As already shown in \cite{GR1}, the assumption $M_2\cap JM=JM_1$ is essential. The Cohen-Macaulay local ring $A=k[[t^4,t^5,t^6,t^7]]$ and the primary ideal $\q=(t^4,t^5,t^6)$ give, with $M=A$ and $L=\q$,  an example where  $e_1(\mathbb{M})=e_0(\mathbb{M})-h_0(\mathbb{M})+1$ but $gr_{\mathbb{M}}(M)$ is not Cohen-Macaulay.
 
 As in \cite{S4} we remark that we can apply  the theorem in the case $\mathbb{M}$ is the $\q$-adic filtration and $\q$ is integrally closed. Namely    Itoh has shown that if $\q$ is any $\m$-primary ideal,  then $$J\cap\overline {\q^2}=J\overline{\q}.$$ Hence, if $\q$ is integrally closed (e.g. $\q=\m$), then $$\q^2\cap J\subseteq J\cap\overline {\q^2}=J\overline{\q}=J\q,$$
so $J\q=\q^2\cap J.$ 
 
 In order  to get  rid of the assumption $M_2\cap JM=JM_1,$  we need one more ingredient,  the study of    the Ratliff-Rush filtration.

\vskip 2mm
In the next section, after discussing the basic properties of this filtration, we apply the results concerning the second Hilbert coefficient $e_2(\mathbb{M}) $ thus completing the study of the equality $$e_1(\mathbb{M})= e_0(\mathbb{M}) -h_0(\mathbb{M}) +1 $$ (see Theorem \ref{sal}).

 \chapter{Bounds for $e_2(\mathbb{M})$ }

If $M $ is a Cohen-Macaulay module, then   in the previous chapter  we showed  that $e_0(\mathbb{M}) $ and $e_1(\mathbb{M}) $ are positive integers. 

In the classical case of an   $\m$-primary ideal $\q $ of a Cohen-Macaulay local ring $A$, as far as
the higher Hilbert coefficients  are concerned, it is
a famous result of M. Narita that $e_2(\q)  \geq 0$ \cite{Na}. This result is  extended here to the case of modules.  In the same paper,
Narita  also showed that if $ \dim A=2,$ then $e_2(\q)=0$ if and only if $\q^n$ has
reduction number one for some large $n.$ Consequently, ${ 
gr}_{\q^n}(A)$ is Cohen-Macaulay. There are examples which show that the result
cannot be extended to higher dimensions.    Very recently Puthenpurakal presented some new results concerning  this problem, see \cite{Puth}. Interesting results on $e_2(\q) $ can   also be   found  in \cite{CPR} which investigates  the interplay between the integrality, or even the normality, of the ideal $\q$ and $e_2(\q).  $  
Classical bounds for  $e_2(\mathbb{M}) $ can be improved and reformulated in  our general setting by using    special  good $\q$-filtrations  on the module $M  $  described in the first  section.  We will introduce the use of  the Ratliff-Rush filtration in studying Hilbert coefficients. This is  a device  which  will be crucial also in the next chapter.

Unfortunately, the positivity does not extend to the higher Hilbert coefficients.  Indeed, in \cite{Na} M. Narita
showed that it is possible for $e_3(\mathbb{M}) $ to be negative. However, a
remarkable result of S. Itoh says that if $\q$ is a normal ideal
then $e_3(\q)  \geq 0$ \cite{I}. A nice  proof of this result was
also given  by S. Huckaba and C. Huneke in \cite{HH}.  In general, it
seems that the integral closedness $($or the normality$)$ of the
ideal $\q $ has  non trivial consequences for  the Hilbert
coefficients of $I$ and, ultimately, for ${ \depth}\, gr_\q(A). $

\section{The Ratliff-Rush filtration}

  Given a good $\q$-filtration on the module $M$, we shall  introduce a new filtration which was constructed  by Ratliff and Rush in \cite{RR}. Here we extend the construction to the general case of a filtered module by following  the definition given by  W.  Heinzer et al. in \cite{HLS}, Section 6.  A further generalization was studied by T.J. Puthenpurakal  and F. Zulfeqarr in \cite{PZ}. 
   \vskip 3mm
Let $\q$ be an $\m$-primary ideal in $A $ and let $\mathbb{M}$ be a good $\q$-filtration on the module $M. $ We define the filtration $\widetilde{\mathbb{M}}$ on $M$ by letting $$
   \widetilde{M}_n:=\bigcup_{k\ge 1}(M_{n+k}:_M\q^k).$$ If  there is no confusion, we will omit  the subscript  $M$ in the   colon.
   It is clear that $\widetilde{M}_0=\widetilde{M}=M$ and, for every $n\ge 0,$ $M_n\subseteq
   \widetilde{M}_n.$   
   
   Further, since $M$ is Noetherian,  there is a positive integer $t$, depending on $n$ such that $$\widetilde{M}_n=M_{n+k}:\q^k  \ \ \ \ \forall  k\ge t.$$     
   \noindent The filtration  $\widetilde{\mathbb{M}}$ is called the \textbf{Ratliff-Rush} filtration associated to $\mathbb{M}.$ 
   \vskip 2mm
   If $\mathbb{M}$ is the $\q$-adic filtration of $M=A, $ then for every integer $n $ there exists an integer $k$ such that 
   $$\widetilde{M}_n=  \widetilde{\q^n}= \q^{n+k} : \q^k.$$
   
   The most important properties of $\widetilde{\mathbb{M}}$ are collected in the following lemma.
   \begin{lemma}\label{rr} Let   $\mathbb{M}$ be  a good $\q$-filtration on the $r$-dimensional module $M, $ such that $\depth_{\q}(M) \ge 1.$  Then we have:
   
   1. There exists an integer $ n_0 $ such that $M_n= \widetilde{M}_n$ for all $n\ge n_0.$
  \vskip 1mm
  2. $\widetilde{\mathbb{M}}$ is a good $\q$-filtration on $M.$
    \vskip 1mm
   3. If $a$ is  $\mathbb{M}$-superficial for $\q$, then it is also $\widetilde{\mathbb{M}}$-superficial for $\q.$
    \vskip 1mm
   4. $\widetilde{\mathbb{M}}$ and $\mathbb{M}$ share the same Hilbert-Samuel  polynomial so that $e_i(\widetilde{\mathbb{M}})=e_i(\mathbb{M})$ for every $i=0,\cdots,r.$
    \vskip 1mm
   5. If $a$   is $\mathbb{M}$-superficial for $\q,$ then  $ \widetilde{M}_{j+1}:a= \widetilde{M}_j$ for every $j\ge 0$, so that $\depth  \ gr_{\widetilde{\mathbb{M}}}(M)\ge 1.$
   \vskip 1mm
   6. $\depth \ gr_{\mathbb{M}}(M)\ge 1$ if and only if $M_n= \widetilde{M}_n$ for every $n.$
    \end{lemma}
   \begin{proof} Let $a$ be  $\mathbb{M}$-superficial for $\q.$ Since  $\depth_{\q}(M) \ge 1, $  $a$ is a regular element for $m$ and, by Theorem \ref{super},  there exists an integer $n$ such that $M_{j+1}:a=M_j$ for every $j\ge n.$ We have 
   \begin{equation*}
   \begin{split}
   \widetilde{M}_n &=M_{n+k}:\q^k\subseteq M_{n+k}:a^k=(M_{n+k}:a):_Ma^{k-1}\\ &=M_{n+k-1}:a^{k-1}=\cdots =M_{n+1}:a=M_n.\\ 
   \end{split}
   \end{equation*}
   which proves the first assertion. Also we have $$\q   \widetilde{M}_n=\q (M_{n+k}:\q^k)\subseteq \q M_{n+k}:\q^k\subseteq M_{n+k+1}:\q^k\subseteq \widetilde{M}_{n+1}$$ which proves that $\widetilde{\mathbb{M}}$ is a $\q$-filtration. Further, since $M_n= \widetilde{M}_n$ for $n\gg 0,$  $\widetilde{\mathbb{M}}$ is $\q$-good, $a$ is $\widetilde{\mathbb{M}}$-superficial for $\q$
   and  $\widetilde{\mathbb{M}}$ and $\mathbb{M}$ share the same Hilbert-Samuel polynomial.  This proves  that $e_i(\widetilde{\mathbb{M}})=e_i(\mathbb{M})$ for every $i=0,\cdots,r.$
   
   We  prove now that $\widetilde{M}_{j+1}:a=\widetilde{M}_j$ for every $j\ge 0.$ It is clear that we can find an integer $k$ such that $M_{j+1+k}:a=M_{j+k}$ and $\widetilde{M}_{j+1}=M_{j+1+k}:\q^k.$ Then we get $$\widetilde{M}_{j+1}:a=(M_{j+1+k}:\q^k):a=(M_{j+1+k}:a):\q^k=M_{j+k}:\q^k\subseteq \widetilde{M}_j.$$
   Finally we must prove that $\text{depth} \ gr_{\mathbb{M}}(M)\ge 1$ implies $\widetilde{M}_n =M_n$ for every $n.$ But we have $$\widetilde{M}_n=M_{n+k}:\q^k\subseteq M_{n+k}:a^k=M_n$$ because $a^*$ is a regular element on $gr_{\mathbb{M}}(M).$
   \end{proof}
\vskip 2mm

The filtration $\widetilde{\mathbb{M}}$ on $M$ is a good $\q$-filtration because,  by the definition,  we have $\q \widetilde{M}_n \subseteq  \widetilde{M}_{n+1} $ and, by Lemma \ref{rr}, for large $n$ we have $\widetilde{M}_{n+1}= {M}_{n+1}=
\q {M}_{n }= \q \widetilde{M}_{n}.$ 
\vskip 2mm
It's worth  recalling  that T.J. Puthenpurakal  and F. Zulfeqarr  in \cite{PZ} and in \cite{Puth}  further analyzed   the case when $\q$ is not a regular ideal, i.e. it does not  contain an $M$-regular element. 

\section{Bounds for  $e_2(\mathbb{M}) $.}

 \vskip 3mm
  The Ratliff-Rush filtration is  a very useful tool in proving results on the Hilbert coefficients,  
but in general it does not behave well in  inductive arguments, except  few, for instance for the  integrally closed ideals. \par \noindent 

Denote by $\overline {\q} $ the integral closure of the ideal $\q $ in $A.$ It is easy to see that
$$ \q \subseteq  \widetilde{\q} \subseteq \overline {\q}. $$
Hence, if $\q$ is an integrally closed  ideal, then $ \widetilde{\q} = \q $  (say that $\q$ is Ratliff-Rush closed).  In this case  the Ratliff-Rush closure commutes with the quotient by a superficial element. In fact S.  Itoh proved that if $\q$ is an integrally closed  ideal,  there exists a superficial element $a \in \q $  such that      $ \q B $ is integrally closed in  $B=A/aA,  $  hence $ \q B $  is Ratliff-Rush closed in $B. $ In particular $$ {\widetilde{\q} }B = {\widetilde{\q B }  }.  $$   This is not true in general. \par
Superficial elements do not behave well,  even if we consider  Ratliff-Rush closed ideals (non integrally closed).  Consider $\q =(x^l, xy^{l-1}, y^l) , l>2,  $  in $ A=k[[x,y]] $ (see \cite{RS}); in this case all the powers of $\q $ are Ratliff-Rush closed,  nevertheless  there is no  superficial element $a \in \q$ for which  $\q/(a) $ is not Ratliff-Rush closed in $B=A/aA, $ hence $ {\widetilde{\q} }B \neq  {\widetilde{\q B }  }.  $   In \cite{Puth}, Theorem 3.3. and Theorem 5.5.,  T. Puthenpurakal gives    a complete characterization of the existence of a superficial element $a\in \q$ for which   $$  {\widetilde{\q^i}B } = {\widetilde{{\q^i }B} }  $$  for every integer $i.$

  \vskip 3mm
     An important fact proved by Huckaba and Marley in  \cite{HM1},  Corollary 4.13 is an  easy consequence of our approach. 
  
 Let us assume  $\mathbb{M}$ is  a good $\q$-filtration on the $2$-dimensional Cohen-Macaulay module $M, $    then by Theorem \ref{e1=CM} c) and  Lemma \ref{rr}, (4) and  (5)  we have  
\begin{equation} \label{HM} e_1(\mathbb{M}) = e_1(\widetilde{\mathbb{M}})=   \sum_{j\ge 0} v_j(\widetilde{\mathbb{M}}) \\ \ \ \ 
   e_2(\mathbb{M}) = e_2(\widetilde{\mathbb{M}}) =  \sum_{j\ge 1} j v_j(\widetilde{\mathbb{M}})  \end{equation} 
   We recall that   $v_j(\widetilde{\mathbb{M}})= \lambda(\widetilde{M}_{j+1}/J \widetilde{M}_{j}) $ where $J$ is a maximal $\mathbb{M}$-superficial sequence for $\q, $ hence by Lemma \ref{rr} (3),  a maximal $\widetilde{\mathbb{M}}$-superficial sequence  for $\q. $
 As a  first application of the previous formula we obtain a short   proof of the non negativity of $e_2(\mathbb{M}).$

 \begin{proposition} \label{e2>0} Let $\mathbb{M}$ be a good $\q$-filtration of theCohen-Macaulay module $M$ of dimension $r.$ Then $$  e_2(\mathbb{M}) \ge 0.$$ 
 \end{proposition}
 \begin{proof} If $r=1, $ it is clear by Lemma \ref{eu}. Let $r \ge 2,  $ by Proposition \ref{ei} we may assume $r=2.$    Hence by (\ref{HM}) 
 \begin{equation*}\begin{split} e_2(\mathbb{M})&=e_2(\widetilde{\mathbb{M}})=\sum_{j\ge 1}jv_j(\widetilde{\mathbb{M}}) \ge 0.
\end{split}
\end{equation*}

\end{proof} 
\vskip 3mm
 The
following example given in \cite{CPR} shows that  $e_2(\q)=0 $ does not imply the
Cohen-Macaulayness of ${ gr}_\q(A)$.

\begin{example}\label{e2=0}{\rm
Let $A$ be the   regular local ring
$k[\![x,y,z]\!]$, with $k$ a field and $x,y,z$ indeterminates and consider 
 $\q = (x^2-y^2, y^2-z^2, xy, xz, yz), $ then 
\[
P_\q(z)=\frac{5 + 6z^2 - 4z^3 +z^4}{(1-z)^3}.
\]
In particular, $e_2(\q)=0 $ and we prove that $ gr_{\q}(A)$ has depth zero. In fact   we can find a superficial element  for $\q$ whose   initial form is not regular on $ gr_{\q}(A).   $   
Computing by CoCoA  the Hilbert coefficients $e_i (\q/(xy))$ we can see that  $xy$ is a
superficial element for $\q  $ (see Remark \ref{criterion}), but  its   initial form is a zero-divisor
on  $ gr_{\q}(A)$ since $P_{\q}(z) \neq P^1_{\q/(xy)}(z)  $ (see Proposition \ref{ei}). }
\end{example}

\vskip 2mm

In the two-dimensional case, one can prove that if $e_2(\mathbb{M})=0,$ then $gr_{\widetilde{\mathbb{M}}}(M)$ is Cohen-Macaulay. We prove this in the next theorem, where    we extend  also results by Sally and Narita (see \cite{S5} and \cite {Na}).

\begin{theorem} \label{e2} Let $\mathbb{M}=\{M_j\}_{j\ge 0}$ be a good $\q$-filtration of the Cohen-Macaulay module $M$ of dimension $2$ and   let  $J$  be an ideal generated by   a maximal  $\mathbb{M}$-superficial sequence  
for $\q. $  Then 
\vskip 2mm
1. $e_2(\mathbb{M})\ge e_1(\mathbb{M})-e_0(\mathbb{M}) +\la(M/\widetilde{M}_1)\ge 0$
\vskip 2mm
2. If $ e_2(\mathbb{M})=0$  and $M_1=\w{M}_1,$ then $e_1(\mathbb{M})=e_0(\mathbb{M})-h_0(\mathbb{M})$ so that $s(\mathbb{M})\le 1$ and $gr_{\mathbb{M}}(M)$ is Cohen-Macaulay.

\vskip 2mm 
3. $gr_{\widetilde{\mathbb{M}}}(M)$ is Cohen-Macaulay if at least one of the following conditions holds:

\vskip 2mm 
  \ \ \ \ \ \  a) $ e_2(\mathbb{M})=0,$

\vskip 2mm
 \ \ \ \ \ \  b) $e_2(\mathbb{M})= e_1(\mathbb{M})-e_0(\mathbb{M}) +\la(M/\widetilde{M}_1)$ and $\widetilde{M}_2\cap JM=J\widetilde{M}_1,$

\end{theorem}
\begin{proof} Since $\depth  \ gr_{\widetilde{\mathbb{M}}}(M)\ge 1=r-1$,  we have $e_2(\mathbb{M})=e_2(\widetilde{\mathbb{M}})=\sum_{j\ge 1}jv_j(\widetilde{\mathbb{M}})$ 
and $e_1(\mathbb{M})=e_1(\widetilde{\mathbb{M}})=\sum_{j\ge 0}v_j(\widetilde{\mathbb{M}}).$ Hence we get
\begin{equation*}\begin{split} e_2(\mathbb{M})&=e_1(\mathbb{M})-v_0(\widetilde{\mathbb{M}})+\sum_{j\ge 2}(j-1)v_j(\widetilde{\mathbb{M}})\\
&= e_1(\mathbb{M})-\la(\widetilde{M}_1/JM)+\sum_{j\ge 2}(j-1)v_j(\widetilde{\mathbb{M}})\\
&= e_1(\mathbb{M})-e_0(\mathbb{M})+\la(M/\widetilde{M}_1)+\sum_{j\ge 2}(j-1)v_j(\widetilde{\mathbb{M}})\\
&\ge e_1(\mathbb{M})-e_0(\mathbb{M})+\la(M/\widetilde{M}_1)\ge 0
\end{split}
\end{equation*} where the last inequality follows by  (\ref{north}) applied to the Ratliff-Rush filtration $\widetilde{M}$. This proves 1. which trivially gives  2.  
As for 3.,  if $e_2(\mathbb{M})=0,$ then $e_1(\mathbb{M})=e_0(\mathbb{M})-\la(M/\widetilde{M}_1)$ and $gr_{\widetilde{\mathbb{M}}}(M)$ is Cohen-Macaulay by Theorem \ref{nor}.

If $e_2(\mathbb{M})= e_1(\mathbb{M})-e_0(\mathbb{M}) +\la(M/\widetilde{M}_1),$ then
$v_j(\widetilde{\mathbb{M}})=0$ for every $j\ge 2.$ This means that $\widetilde{M}_{j+1}=J\widetilde{M}_j$ for every $j\ge 2$, and since $\widetilde{M}_2\cap JM=J\widetilde{M}_1,$ $gr_{\widetilde{\mathbb{M}}}(M)$ is Cohen-Macaulay by Valabrega-Valla.
\end{proof}

The inequality $e_2(\mathbb{M})\ge e_1(\mathbb{M})-e_0(\mathbb{M}) +\la(M/\widetilde{M}_1)\ge 0$ was proved by Sally in \cite{S5}, Corollary 2.5, in the special case $M_j=\q^j.$ The methods there involve the local cohomology of the Rees ring. 
\vskip 3mm

  As a  consequence of the previous theorem  we obtain  a classical result by Narita (see \cite{Na}).      

\begin{corollary} \label{NS} Let $\q$ be an $\m$-primary ideal  which is integrally closed in a Cohen-Macaulay local ring $(A, \m) $ of dimension $r \ge 2.$  Then
 \vskip 2mm
1.  $e_2(\q) \ge e_1(\q) - e_0(\q) + \la(A/\q) $ 
 \vskip 2mm
2.  If $e_2(\q)=0 $ then $gr_\q(A) $ is Cohen-Macaulay and $e_i(\q)=0 $ for  $i \ge 2.$
  
 \end{corollary}
  
In the above corollary we get rid  of the assumption $\dim (M)=2 $ of Theorem \ref{e2} since $\q$ is integrally closed and,  by Proposition  \ref{Itoh},  we can find a superficial sequence $a_1, \dots, a_{r-2} $ in $\q$  such that $ \q/(a_1, \dots, a_{r-2}) $ is integrally closed in $A/(a_1, \dots, a_{r-2}) $ which is a local Cohen-Macaulay ring of dimension two. Moreover,  by Proposition \ref{ei}, the   numerical invariants involved  are preserved going modulo  $(a_1, \dots, a_{r-2}) $ and we may apply Theorem \ref{e2} in order to prove $1. $  The second assertion comes from $1. $ and Theorem \ref{nor}.  
\vskip 2mm
Example \ref{e2=0} shows  that in Corollary \ref{NS}  the assumption  that 
the ideal $\q $ is   integrally closed  cannot be weakened.  In the example  $e_2(\q)=0, $ but $gr_{\q}(A) $ is not Cohen-Macaulay.

Notice that  inequality $1.$ was already proved by Itoh in \cite[12]{I}. Assertion 2.  of Corollary \ref{NS} was proved by Puthenpurakal in \cite{P2}.    \vskip 3mm
 Later, it
was  conjectured by Valla \cite[6.20]{Val} that if the
equality $e_2= e_1 - e_0 + \lambda(A/\q)$ holds when  $\q $ is the maximal ideal ${\mathfrak m}$ of $A, $ then the
associated graded ring ${\rm gr}_{\mathfrak m}(A)$ is
Cohen-Macaulay. Unfortunately, the following example given by Wang
shows that the conjecture is false.

\begin {example} \label{Wang2}  {\rm
Let $A $ be the two dimensional local Cohen-Macaulay ring
\[
k[\![x,y,t,u,v]\!]/(t^2,tu,tv,uv,yt-u^3,xt-v^3),
\]
with $k$ a field and $x,y,z,u,v $ indeterminates. Let ${\mathfrak m}$   be the
maximal ideal  of $A$. One has that the associated
graded ring ${\rm gr}_{{\mathfrak m}}(A )$ has depth zero and
\[
P_{ {\mathfrak m}}(z)= \frac{1 + 3z + 3z^3 -z^4}{(1-z)^2}.
\]
In particular, one has $e_2 = e_1 - e_0 + 1$, that is, $e_2$ is
minimal according to Itoh's bound. }
\end{example}

Hence the  condition   $\lambda(A/\q )=e_0(\q)-e_1(\q)+e_2(\q)$ does not imply 
  that $  gr_{\q}(A)$ is Cohen-Macaulay
even for an integrally closed ideal $\q$. However,  
 Corso, Polini and Rossi in \cite{CPR}  proved  that the conjecture is true if $\q$ is normal (i.e. $\q^n$ is integrally closed for every $n$).  
 
 \begin{theorem} Let   $\q$ be a normal  $\m$-primary ideal    in a Cohen-Macaulay local ring $(A, \m) $ of dimension $r \ge 2.$  Then
 \vskip 2mm
1.  $e_2(\q) \ge e_1(\q) - e_0(\q) + \la(A/\q). $ 
\vskip 2mm
2. If $e_2(\q) =  e_1(\q) - e_0(\q) + \la(A/\q) $  then $gr_{\q}(A) $ is Cohen-Macaulay.
 
\end{theorem}

 \vskip 2mm
In the above result 1. follows by Corollary \ref{NS}. The crucial point in the proof of 2. is  to prove  that  the reduction number of $\q$ is at most two.  Here a  result  by   Itoh in \cite{I1} had been fundamental. Another context where the normality plays an important role in  the reduction number of $\q$ is in pseudo-rational two-dimensional normal local ring, see a notable result by  Lipman and  Teissier \cite{LT}, Corollary 5.4.

\vskip 2mm
We present  here a
short proof of a result of Narita for modules which  characterizes  $e_2(\mathbb{M})=0 $  
when $M$ is a Cohen-Macaulay module of dimension two and $\mathbb{M} $ is the $\q$-adic filtration on $M.$   We will write $e_i(\q^n M) $ when we consider the Hilbert coefficients of the $\q^n$-adic filtration on $M $ with $n$ a fixed integer.

\begin{proposition}\label{Narita}
  Let $\q$ be an $\m$-primary ideal and let $ M$ be     a Cohen-Macaulay module  of dimension $two.  $ 
 
\noindent Then $e_2(\q M )=0 $  if and only if $\q^nM $ has reduction number one for
some positive integer $n$. Under these circumstances   ${\rm
gr}_{\q^nM}(M)$ is Cohen-Macaulay.
\end{proposition}
\begin{proof} We first recall that $e_2(\q M) =e_2(\q^m M )$ for every positive
integer $m$. Assume $e_2(\q M) =0$ and let $n$ be an integer such that
$\widetilde{\q^nM }=\q^n M$. Hence by Theorem  \ref{e2}, $0= e_2(\q^n M) \ge
e_1(\q^n M) -e_0(\q^n M) + \lambda(M/\widetilde{\q^n M}) = e_1(\q^nM)
-e_0(\q^n M) + \lambda(M/ {\q^n M})$. Hence $e_1(\q^n M) - e_0(\q^n M) +
\lambda(M/ {\q^n M})=0$ because it cannot be negative by
Northcott's inequality. The result follows now by Theorem \ref{nor}.  For the converse, if
$\q^nM $ has reduction number one for some $n$, then $e_2(\q^n M) = 0$ and  ${ 
gr}_{\q^nM}(M)$ is Cohen-Macaulay
again by Theorem \ref{nor}.  In particular $e_2(\q M )=e_2(\q^n M)=0$.
  \end{proof} 

\medskip

We remark that Narita's result cannot be extended to a local
Cohen-Macaulay ring of dimension $>2$ without changing the statement. The ideal $\q $ presented in
Example~\ref{e2=0} satisfies $e_2(\q) =0$, however $\q^n$ has  
reduction number greater than one for every $n$. In fact, it is enough to
remark that $\q $ does not have reduction number one $({  gr}_{\q}(A)$ is
not Cohen-Macaulay$)$ and $\q^n=(x,y,z)^{2n}$ for $n >1$ which has
reduction number two. An extension of Narita's result to the the higher dimensional case is given
by Puthenpurakal in \cite{Puth}.

\medskip

We can prove now the following result which is, at the same time, an extension and a completion of a deep  theorem proved by Sally in \cite{S4}. This result is new even in the classical case and it completes the study of the equality   $e_1(\mathbb{M})=e_0(\mathbb{M})-h_0(\mathbb{M})+1 $ (see Chapter 2).   

\begin{theorem} \label{sal} Let $M$ be a     Cohen-Macaulay module of dimension $r\ge 2$, $L$ a submodule of $M$ and $\mathbb{M}=\mathbb{M}_L$ the good $\q$-filtration on $M$ induced by $L. $ Assume that  $e_1(\mathbb{M})=e_0(\mathbb{M})-h_0(\mathbb{M})+1.$ Then the following conditions are equivalent:
\vskip 2mm 
1. $e_2(\mathbb{M})\not=0$
\vskip 2mm 
2. $e_2(\mathbb{M})=1$
\vskip 2mm 
3. $\depth  \ gr_{\mathbb{M}}(M)\ge r-1.$
\vskip 2mm 
4.  $P_{\mathbb{M}}(z)=\frac{h_0(\mathbb{M})+h_1(\mathbb{M})z+z^2}{(1-z)^r}.$

\end{theorem}
\begin{proof} We recall that we are considering the filtration
$$ \mathbb{M}_L:  \  \  M \supseteq L \supseteq \q L \supseteq  \dots  \supseteq \q^j L \supseteq  \dots $$
First we prove that 1., 2. and 3. are equivalent. As usual, $J$  is an ideal generated by a maximal sequence of $\mathbb{M}$-superficial  elements 
for $\q.$ Let us first consider the case $r=2.$ We have \begin{equation*}v_0(\widetilde{\mathbb{M}}) =\la(\widetilde{M}_1/JM)
=e_0(\mathbb{M})-h_0(\mathbb{M})+\la(\widetilde{M}_1/M_1).\end{equation*} and  $\depth  \ gr_{\widetilde{\mathbb{M}}}(M)\ge 1=r-1.$ This implies  $$e_2(\mathbb{M})=e_2(\widetilde{\mathbb{M}})=\sum_{j\ge 1}jv_j(\widetilde{\mathbb{M}})$$ 
\begin{equation*}
e_0(\mathbb{M})-h_0(\mathbb{M})+1=e_1(\mathbb{M})=e_1(\widetilde{\mathbb{M}})=\sum_{j\ge 0}v_j(\widetilde{\mathbb{M}})\end{equation*}
 
so that  \begin{equation}\label{bla}
\sum_{j\ge 1}v_j(\widetilde{\mathbb{M}})=\sum_{j\ge 0}v_j(\widetilde{\mathbb{M}})-v_0(\widetilde{\mathbb{M}})=1-\la(\widetilde{M}_1/M_1).
\end{equation}
Let us assume that 1. holds,  then $\sum_{j\ge 1}v_j(\widetilde{\mathbb{M}})>0,$ so that $\sum_{j\ge 1}v_j(\widetilde{\mathbb{M}})=1$ and $\la(\widetilde{M}_1/M_1)=0.$
But if $\widetilde{M}_1=M_1,$ we cannot have  $v_1(\widetilde{\mathbb{M}})=0,$ otherwise $$M_2\subseteq \widetilde{M}_2=J\widetilde{M}_1=JM_1,$$  and, by 
Theorem \ref{nor}, $e_1(\mathbb{M})=e_0(\mathbb{M})-h_0(\mathbb{M}),$ a contradiction.
Hence $v_1(\widetilde{\mathbb{M}})=1$ and $v_j(\widetilde{\mathbb{M}})=0$ for every $j\ge 2,$
which implies $e_2(\mathbb{M})=1.  $  This proves that 1. implies 2. 

Let now assume that $e_2(\mathbb{M})=1.$ Then we must have $v_1(\widetilde{\mathbb{M}})=1$ and $v_j(\widetilde{\mathbb{M}})=0$ for every $j\ge 2,$ so that, by (\ref{bla}), $\widetilde{M}_1=M_1.$ Thus 
$$1=\la( \widetilde{M}_2/JM_1)\ge \la(M_2/JM_1)\ge 1$$ which implies $\widetilde{M}_2=M_2.$ Now if $j\ge 2$ and $\widetilde{M}_j=M_j$, then we have \begin{equation*}\label {i} M_{j+1}\subseteq \widetilde{M}_{j+1}=J \widetilde{M}_{j}=JM_j\subseteq M_{j+1}. \end{equation*} Hence, by induction, we get $\widetilde{M}_{t}=M_t$ for every $t\ge 1.$  By the above Lemma, this  implies $\text{depth} \ gr_{\mathbb{M}}(M)>0,$ thus proving that 2. implies 3. 

Finally, the condition $ {depth} \ gr_{\mathbb{M}}(M)>0$ implies $\widetilde{M}_1=M_1$, hence $\sum_{j\ge 1}v_j(\widetilde{\mathbb{M}})= 1$ and $e_2(\mathbb{M})\not=0.$ This completes  the proof of the equivalence of 1., 2. and 3. in the case $r=2.$ 

Let us now consider the general case, when $r\ge 3.$ Let $\frak{a}$ be an ideal generated by an $\mathbb{M}$-superficial sequence for $\q$ of length $r-2.$ Then we have $e_i(\mathbb{M})=e_i(\mathbb{M}/\frak{a}M)$ for $i=0,1,2$ and $h_0(\mathbb{M})=h_0(\mathbb{M}/\frak{a}M).$ Hence   the assumption holds for the $2$-dimensional Cohen-Macaulay module $M/\frak{a}M.$ The conclusion follows because, by Sally machine, $\text{depth} \ gr_{\mathbb{M}}(M)\ge r-1$ if and only if $\text{depth} \ gr_{\mathbb{M}}(M/\frak{a}M)\ge 1.$

We end the proof of the theorem by proving that 2. is equivalent to 4. We notice that if $e_2(\mathbb{M})=1$
and $\depth \ gr_{\mathbb{M}}(M)\ge r-1,$
then, by Theorem \ref{e1=CM}, we get $$1=e_2(\mathbb{M})=\sum_{j\ge 1}jv_j(\mathbb{M}).$$ Hence $v_1(\mathbb{M})=1$ and $v_j(\mathbb{M})=0$ for $j\ge 2$. Since $e_i(\mathbb{M})=\sum_{j\ge i-1}\binom{j}{i-1}v_j(\mathbb{M}),$ we also get $e_j(\mathbb{M})=0$ for $j\ge 3.$ These values of the $e_i$' s  give the required Hilbert series and conversely.

\end{proof}

\vskip 2mm The following example from  \cite{S5} shows that in the above result the assumption  $e_2(\mathbb{M})\not=0$ is essential. We remark that the $\q$-adic filtration of $A$ is a filtration of the type $\mathbb{M}_L$  induced on $A$ by $L=\q$ itself. 
  

\vskip 2mm 
\begin{example} {\label{1}}{\rm {Consider the ideal $q=(x^4,x^3y,xy^3,y^4) \subseteq  A=k[[x,y]]. $  The ideal $\q$ is not integrally closed and if we consider on $A$ the $\q$-adic filtration,   we have $$P_{\q}(z)=\frac{11+3z+3z^2-z^3}{(1-z)^2}.$$  This gives $e_0=16,$ $e_1=6,$ $e_2=0,$ $h_0(\q)=11,$  so that 
$e_1(\q)=e_0(\q)-h_0(\q)+1.$ It is clear that $x^2y^2\notin \q$ while $x^2y^2 \q\subseteq \q^2$ so that $\w{\q}\not= \q$ and 
$gr_{\q}(A)$ has depth zero by Lemma \ref{rr}.  }}
\end{example} 
\vskip 2mm  
 Very little is known about  the Hilbert Function of the filtered module $M$ when $e_2(\mathbb{M})=0 $ and $ M_1\not =\w{M}_1. $ As completion of Theorem \ref{sal}   and by using the machinery introduced in the next section,  we will give an  answer,  provided       $e_1(\mathbb{M})=e_0(\mathbb{M})-h_0(\mathbb{M})+1.$

\vskip 2mm The following example shows that in the above theorem the assumption $e_1(\mathbb{M})=e_0(\mathbb{M})-h_0(\mathbb{M})+1$ and $e_2(\mathbb{M})=1$ does not imply $gr_{\mathbb{M}}(A)$ is Cohen-Macaulay.
\vskip 1mm
\begin{example} {\label{2}}{\rm{Consider the ideal   $\q=(x^6,x^5y^3,x^4y^7,x^3y^8,  x^2y^{10},xy^{11}, $ $ y^{22})   $  in $A=k[[x,y]. $ We have 
$$P_{\q}(z)=\frac{61+26z+z^2}{(1-z)^2}.$$ This gives  $e_0(\q)=88,$ $e_1(\q)=28,$ $e_2(\q)=1,$ $h_0(\q)=61,$  so that we have 
$e_1(\q)=e_0(\q)-h_0(\q)+1.$ However $gr_{\q}(A)$ is not Cohen-Macaulay because $\q$ is an $\m$-primary  ideal in a regular ring of dimension two and $s(\q)=2>1$ (see \cite{HM1} Theorem A and Proposition 2.6, \cite{CDJR} Proposition 2.9).}}
\end{example} 
\vskip 2mm
The previous example underlines  the difference between the case of  the $\m$-adic filtration  and the more general case of an $\m$-primary ideal $\q. $ In the first case \cite{EV} Elias and Valla  proved that if the degree of the $h$-polynomial is less than or equal to two, then the associated graded ring is Cohen-Macaulay.  The above example shows that this is not the case when $\q$ is not maximal, even if $A$ is regular and $h_2=1.$

\vskip 2mm

  The second statement in Theorem 2 of  \cite{P2} says that if $M=A$ is Cohen-Macaulay, $\q=\m,$ $\dim A=2$ and $e_1(A)=2e_0(A)-\mu(\m)+1,$ then either $gr_{\m^n}(A)$ is Cohen-Macaulay for $n\gg 0$ or $\text{depth} \  gr_{\m}(A)\ge 1.$
 
 We notice that we have $$2e_0(A)-\mu(\m)+1=2e_0(A)-2h_0(A)-h(A)+1,$$ so that this last result will be a consequence of the following theorem which is a further step after Theorem \ref{elv}.
 
\begin{theorem} \label{r=2} Let $M$ be a    Cohen-Macaulay module of dimension $r\ge 2$, $L$ a submodule of $M$ and $\mathbb{M}=\mathbb{M}_L$ the good  $\q$-filtration on $M$ induced by $L.$ Let $J$  be  generated by   a maximal  $\mathbb{M}$-superficial sequence  
for $\q $  and assume that  $e_1(\mathbb{M})=2e_0(\mathbb{M})-2h_0(\mathbb{M})-h(\mathbb{M})+1,$  $\widetilde{M}_1=M_1$ and $\widetilde{M}_2\cap JM=JM_1.$   Then we have 
  \vskip 2mm 
 1. If  $\widetilde{M}_2\not=M_2,$ then $ gr_{\widetilde{\mathbb{M}}}(M)$ is Cohen-Macaulay.
 \vskip 2mm
 2. If $r=2$, then 
  \vskip 2mm   a) $e_2(\mathbb{M})=e_0(\mathbb{M})-h_0(\mathbb{M})-h_1(\mathbb{M})+1$ if and only if $\depth  \ gr_{\mathbb{M}}(M)=0.$
\vskip 2mm   b) $e_2(\mathbb{M})=e_0(\mathbb{M})-h_0(\mathbb{M})-h_1(\mathbb{M})+2$ if and only if $\depth \ gr_{\mathbb{M}}(M)\ge 1.$
 \vskip 2mm  Further, in case a), $ gr_{\widetilde{\mathbb{M}}}(M)$ is Cohen-Macaulay; in case b), 
$$P_{\mathbb{M}}(z)=\frac{h_0(\mathbb{M})+h_1(\mathbb{M})z+h_2(\mathbb{M})z^2+z^3}{(1-z)^2}.$$
 \end{theorem}
\begin{proof} We have $e_i(\mathbb{M})=e_i(\widetilde{\mathbb{M}})$ for $i=0,1,2$ $$h_0(\widetilde{\mathbb{M}})=\la(M/\widetilde{M}_1)=\la(M/M_1)=h_0(M)$$ and 
\begin{equation*}
\begin{split}
h(\widetilde{\mathbb{M}})&=\la(\widetilde{M}_1/JM+\widetilde{M}_2)=\la(M_1/JM+\widetilde{M}_2 )\\
&=\la(M_1/JM+M_2 )-\la(JM+\widetilde{M}_2/JM+M_2)\\
&=h(\mathbb{M})-\la\left (\widetilde{M}_2/M_2+(\widetilde{M}_2\cap JM)\right)\\
&=h(\mathbb{M})-\la (\widetilde{M}_2/M_2).
\end{split}
\end{equation*}

Since $M_2\cap JM\subseteq \widetilde{M}_2\cap JM=JM_1,$ we also have $M_2\cap JM=JM_1,$ which implies by (\ref{h=h1+}) $$h(\mathbb{M})=h_1(\mathbb{M}).$$ Further 
\begin{equation*}
\begin{split} 2e_0(\mathbb{M})-2h_0(\mathbb{M})-h(\mathbb{M})+1&=e_1(\mathbb{M})=e_1(\widetilde{\mathbb{M}})\\
&\ge 2e_0(\widetilde{\mathbb{M}})-2h_0(\widetilde{\mathbb{M}})-h(\widetilde{\mathbb{M}})\\
&=2e_0(\mathbb{M})-2h_0(\mathbb{M})-h(\mathbb{M})+\la (\widetilde{M}_2/M_2)
\end{split}
\end{equation*} so that $0\le \la (\widetilde{M}_2/M_2)\le 1.$
\vskip 1mm
If $\widetilde{M}_2\not = M_2,$ then $\la (\widetilde{M}_2/M_2)= 1$ and $$e_1(\widetilde{\mathbb{M}})
= 2e_0(\widetilde{\mathbb{M}})-2h_0(\widetilde{\mathbb{M}})-h(\widetilde{\mathbb{M}})$$
so that $gr_{\widetilde{\mathbb{M}}}(M)$ is Cohen-Macaulay by Theorem \ref{elv}. This proves 1. 

Let us prove 2. We have $r=2$ and $\depth \ gr_{\widetilde{\mathbb{M}}}(M)\ge 1=r-1$ so that, by Theorem \ref{e1=CM}, $$e_1(\mathbb{M})=e_1(\widetilde{\mathbb{M}})=\sum_{j\ge 0}v_j(\widetilde{\mathbb{M}}),\ \ \ \ \ e_2(\mathbb{M})=e_2(\widetilde{\mathbb{M}})=\sum_{j\ge 1}jv_j(\widetilde{\mathbb{M}}).$$
Now $$v_0(\widetilde{\mathbb{M}})=\la(\widetilde{M}_1/JM)=\la(M_1/JM)=e_0(\mathbb{M})-h_0(\mathbb{M})$$

\begin{equation*}
\begin{split}
v_1(\widetilde{\mathbb{M}})&=\la(\widetilde{M}_2/J\widetilde{M}_1)=\la(\widetilde{M}_2/JM_1)=\la(\widetilde{M}_2/M_2)+\la(M_2/JM_1)\\
&=\la(\widetilde{M}_2/M_2)+e_0(\mathbb{M})-h_0(\mathbb{M})-h_1(\mathbb{M})\\
&=\la(\widetilde{M}_2/M_2)+e_0(\mathbb{M})-h_0(\mathbb{M})-h(\mathbb{M}),
\end{split}
\end{equation*} where we used the equality $e_0(\mathbb{M})=h_0(\mathbb{M})+h_1(\mathbb{M})+\la(M_2/JM_1)$ proved in Proposition \ref{ab}. 

\noindent This implies 
{\small{\begin{equation*}
\begin{split} 2e_0(\mathbb{M})-2h_0(\mathbb{M})-h(\mathbb{M})+1&= e_1(\mathbb{M})=
e_1(\widetilde{\mathbb{M}})\\
&=v_0(\widetilde{\mathbb{M}})+v_1(\widetilde{\mathbb{M}})+\sum_{j\ge 2}v_j(\widetilde{\mathbb{M}})\\
&=\la(\widetilde{M}_2/M_2)+2e_0(\mathbb{M})-2h_0(\mathbb{M})-h(\mathbb{M})+\sum_{j\ge 2}v_j(\widetilde{\mathbb{M}})
\end{split}
\end{equation*} }} so that $$\la(\widetilde{M}_2/M_2)+\sum_{j\ge 2}v_j(\widetilde{\mathbb{M}})=1.$$
In the case $\sum_{j\ge 2}v_j(\widetilde{\mathbb{M}})=1,$ we have $\widetilde{M}_2=M_2$ and 
$e_1(\mathbb{M})=v_0(\widetilde{\mathbb{M}})+v_1(\widetilde{\mathbb{M}})+1.$ We claim that this implies $M_3\not=JM_2$ and $v_2(\widetilde{\mathbb{M}})=1.$ Namely, if  $M_3=JM_2,$ then $\q^2L=J\q L$ so that $M_{j+1}=JM_j$ for every $j\ge 2$. Since $M_2\cap JM=JM_1,$ by Valabrega-Valla $gr_{\mathbb{M}}(M)$ is Cohen-Macaulay with $v_j(\mathbb{M})=0$ for every $j\ge 2.$ But then $e_1(\mathbb{M})
=v_0(\mathbb{M})+v_1(\mathbb{M}),$  a  contradiction.

Hence $M_3\not=JM_2,$ so that $$J\widetilde{M}_2=JM_2\subset M_3\subseteq \widetilde{M}_3$$ and $v_2(\widetilde{\mathbb{M}})=1.$  This proves the claim.
\vskip 1mm Now we can write 
\begin{equation*}
\begin{split} e_2(\mathbb{M})=
e_2(\widetilde{\mathbb{M}})&=v_1(\widetilde{\mathbb{M}})+\sum_{j\ge 2}jv_j(\widetilde{\mathbb{M}})\\
&=\la(\widetilde{M}_2/M_2)+e_0(\mathbb{M})-h_0(\mathbb{M})-h_1(\mathbb{M})+\sum_{j\ge 2}jv_j(\widetilde{\mathbb{M}})\\
&=e_0(\mathbb{M})-h_0(\mathbb{M})-h_1(\mathbb{M})+1+\sum_{j\ge 2}(j-1)v_j(\widetilde{\mathbb{M}}).
\end{split}
\end{equation*}
Hence  we have only two possibilities for $e_2(\mathbb{M}),$ namely 

\begin{equation}\label{for} e_2(\mathbb{M})=
\begin{cases}
  e_0(\mathbb{M})-h_0(\mathbb{M})-h_1(\mathbb{M})+1   & \text{if} \ \ \ \ \sum_{j\ge 2}v_j(\widetilde{\mathbb{M}})=0, \\
 e_0(\mathbb{M})-h_0(\mathbb{M})-h_1(\mathbb{M})+2     & \text{otherwise}.
\end{cases}\end{equation}

Now, if $\depth  \ gr_{\mathbb{M}}(M)\ge 1,$ then $ \widetilde{M}_2=M_2,$ hence 
$\sum_{j\ge 2}v_j(\widetilde{\mathbb{M}})=1$ and we have $e_2(\mathbb{M})=e_0(\mathbb{M})-h_0(\mathbb{M})-h_1(\mathbb{M})+2.$

Conversely, if $e_2(\mathbb{M})=e_0(\mathbb{M})-h_0(\mathbb{M})-h_1(\mathbb{M})+2,$
then $v_2(\widetilde{\mathbb{M}})=1$ and $v_j(\widetilde{\mathbb{M}})=0$ for every $ j\ge 3.$ This implies  
$$1=\la(\widetilde{M}_3/J\widetilde{M}_2)=\la(\widetilde{M}_3/JM_2)\ge \la(M_3/JM_2)\ge 1,$$ so that $\widetilde{M}_3=M_3.$ Hence, since $v_j(\widetilde{\mathbb{M}})=0$ for every $ j\ge 3,$ we get $\widetilde{M}_j=M_j$ for  every $ j\ge 0,$ which is equivalent to $\depth  \ gr_{\mathbb{M}}(M)\ge 1.$

This proves a); as for b), it follows by a) and (\ref{for}). We come now to the last assertions of the theorem.

In case a), we have $\widetilde{M}_2\not=M_2,$ so that $ gr_{\widetilde{\mathbb{M}}}(M)$ is Cohen-Macaulay by 1. In case b), $\mathbb{M}=\widetilde{\mathbb{M}}$  so that $v_2(\mathbb{M})=v_2(\widetilde{\mathbb{M}})=1$ and $v_j(\mathbb{M})=v_j(\widetilde{\mathbb{M}})=0$ for every $j\ge 3.$ Since by Theorem \ref{e1=CM} we have 
$e_i(\mathbb{M})=\sum_{j\ge i-1}\binom{j}{i-1}v_j(\mathbb{M}),$ we get 
$e_3(\mathbb{M})=1$ and $e_j(\mathbb{M})=0$ for every $j\ge 4.$
These values of the $e_i$'s give the required Hilbert series.
\end{proof}
 
 The assumptions $\widetilde{M}_1=M_1$ and $\widetilde{M}_2\cap JM=JM_1$  in the above theorem seem very strong, but they are satisfied  by the $\q$-adic filtration of any  primary integrally closed ideal $\q$.  In particular by the $\m$-adic filtration.
 \begin{corollary} Let $\q$ be an $\m$-primary ideal in the Cohen-Macaulay local ring $A$ of dimension $r$ and $\mathbb{M}$ the $\q$-adic filtration on $A.$  If $\q$ is integrally closed and $e_1(\mathbb{M})=2e_0(\mathbb{M})-2h_0(\mathbb{M})-h(\mathbb{M})+1,$ the following conditions are equivalent and each implies $$P_{\mathbb{M}}(z)=\frac{h_0(\mathbb{M})+h_1(\mathbb{M})z+h_2(\mathbb{M})z^2+z^3}{(1-z)^r}.$$
  
 \ \ \ \ \ a) $\depth  \ gr_{\mathbb{M}}(M)\ge r-1.$
  \vskip 2mm 
 \ \ \ \ \ b) $e_2(\mathbb{M})=e_0(\mathbb{M})-h_0(\mathbb{M})-h_1(\mathbb{M})+2.$
 
 \vskip 2mm \noindent If this is not the case, then  $e_2(\mathbb{M})= e_0(\mathbb{M})-h_0(\mathbb{M})-h_1(\mathbb{M})+1.$ \end{corollary}
 \begin{proof}  Since $\q\subseteq \widetilde{\q}\subseteq \overline{\q},$ we have $\q=\widetilde{\q};$ on the other hand, if  $J$  is an ideal generated by a maximal sequence of $\mathbb{M}$-superficial  elements for $\q,$ by a result of Huneke and Itoh (see \cite{H} and \cite{I}), we have 
 $$\widetilde{\q^2}\cap J\subseteq \overline{\q^2}\cap J=J\overline{\q}=J\q$$ so that 
 $$\widetilde{\q^2}\cap J=J\q.$$ Hence the equivalence between a) and b) follows by the theorem  if $r=2.$ 
When  $r\ge 3, $ by a result of Itoh (see \cite{I1}), we can find  an ideal $\frak{a}$ generated by an $\mathbb{M}$-superficial sequence for $\q$ of length $r-2$ such that $\q/\frak{a}$ is integrally closed. Then we have $e_i(\mathbb{M})=e_i(\mathbb{M}/\frak{a})$ for $i=0,1,2,$  $$h_0(\mathbb{M})=\la (A/J)=h_0(\mathbb{M}/\frak{a})$$ and $$h_1(\mathbb{M})=h(\mathbb{M})=\la(A/J+\q^2)=h(\mathbb{M}/\frak{a})=h_1(\mathbb{M}/\frak{a}).$$ Hence   all the assumptions of the theorem  hold for the $2$-dimensional Cohen-Macaulay local ring $A/\frak{a}$ and the integrally closed primary ideal $\q/\frak{a}.$ The equivalence between a) and b) follows by the theorem because, by Sally machine, $\depth  \ gr_{\mathbb{M}}(A)\ge r-1$ if and only if $\depth  \ gr_{\mathbb{M}}(A/\frak{a})\ge 1.$
 
 As for the last assertion  if   $ \depth \ gr_{\mathbb{M}}(M) <  r-1,$ then by using Sally's machine, we deduce   $\depth  \ gr_{\mathbb{M}}(A/\frak{a})= 0.$ Since $\q/\frak{a} $ is integrally closed and $A/\frak{a} $ is a $2$-dimensional local Cohen-Macaulay ring, we may apply  Theorem \ref{r=2}  and it is easy to see  that $e_2(\mathbb{M})=e_0(\mathbb{M})-h_0(\mathbb{M})-h_1(\mathbb{M})+1$ since the integers involved do not
 change passing to $\mathbb{M}/\frak{a}.$ 
\end{proof}

 \chapter{Sally's conjecture and applications}

Let $\q $ be an $\m$-primary ideal  of $A $ and  let $M$ be  a  Cohen-Macaulay $A$-module of dimension $ r. $ Consider the good $\q$-filtration $\mathbb{M} =\mathbb{M}_L$ induced by a  submodule $L$ of $M$ (see (\ref{find})).  
 In Theorem \ref{nor} we proved that, if $\mathbb{M} $ has  minimal multiplicity, namely $e_0(\mathbb{M})=h_0(\mathbb{M})+h_1(\mathbb{M}), $ then $$P_{\mathbb{M}}(z)=\frac{h_0(\mathbb{M})+h_1(\mathbb{M})z}{(1-z)^r} $$ and gr$_{\mathbb{M}}(M) $ is Cohen-Macaulay.
 
\noindent Furthermore, minimal multiplicity is equivalent to have  $e_1(\mathbb{M})=e_0(\mathbb{M})-h_0(\mathbb{M}),$ the minimal value  with respect to   Northcott's  bound, and  Hilbert series
   $$P_{\mathbb{M}}(z)=\frac{h_0+h_1z }{(1-z)^r }. $$ 
   \vskip 2mm
   In the next case, when $e_0(\mathbb{M})= h_0(\mathbb{M})+h_1(\mathbb{M})+1,$ we say that $\mathbb{M} $ has almost  minimal multiplicity. 
 Almost minimal multiplicity is much more difficult to handle, even for   the $\m$-adic filtration on a Cohen-Macaulay local ring. 
  In this particular case $h_0(\m)=1 $ and $h_1(\m) =  \mu(\m) -r  = h,  $   the embedding codimension.  Then $$ e_0(\m) = h +2 $$
For example the Cohen-Macaulay one-dimensional local ring $A=k[[t^4,t^5, t^{11}]]$ has almost minimal multiplicity (with respect the $\m$-adic filtration) and its Hilbert series is $$P_{\m}(z)=\frac{1+ h z+z^3}{(1-z)}, $$ but the associated graded ring is not Cohen-Macaulay. 
 
 It was conjectured by  Sally in \cite {S3} that, for an  $r$-dimensional Cohen-Macaulay local ring, in the classical case of the $\m$-adic filtration, almost minimal multiplicity forces the depth of the associated graded ring to be at least $ r-1.$
After $13$ years,  the conjecture was proved by  Wang in \cite{W1} and at the same time by  Rossi and   Valla in \cite{RV1}. In particular  it was  proved that an  $r$-dimensional Cohen-Macaulay local ring $A$ has  almost minimal multiplicity if and only if 
 $$P_A(z)=\frac{1+h z+z^s}{(1-z)^r }$$ for some integer $s$ such that $2\le s\le e_0 -1.$ 
 
\noindent   Later the conjecture was  stated  for any $\m$-primary ideal of a Cohen-Macaulay ring and an extended version  was proved  in \cite{H2}, \cite{CPVa}, \cite{E2} and \cite{R2} by following Rossi and Valla's proof.    
 
 \vskip 3mm In this chapter  we present a proof of  this result in the general case of a module endowed with the filtration induced by $L.$ The crucial point of  this result  is a bound on the reduction number of $\mathbb{M}.$

\noindent   As we have already seen in Chapter 1, if $J$ is an ideal generated by a maximal  $\mathbb{M}$-superficial sequence  
for $\q, $ then $\mathbb{M}$ is a  good $J$-filtration and hence for large $n$ we have $$ M_{n+1}=J M_n.$$ 
In particular $J$ is a  minimal $\mathbb{M}$-reduction of   $\q.$
 Following  the classical theory of  reductions  of an ideal,  we denote by 
 $$ r_J(\mathbb{M}):= \text{min} \{ n \in {\bf N}\  |\   M_{j+1}=J M_j \ \ \text{for every} \ \ j\ge n \} $$
the {\it reduction number} of  $\mathbb{M} $ with respect to $J.$ Since  $\mathbb{M}=   \mathbb{M}_L$ for a given submodule $L$ of $M,$ 
 we clearly have $$ r_J(\mathbb{M}):= \text{min} \{ n \in {\bf N}\  |\   M_{n+1}=J M_n \}=\text{min} \{ n \in {\bf N}\ | \   v_n(\mathbb{M})=0 \}. $$
If $\mathbb{M}$ is the $\q$-adic filtration on the ring $A, $ then we write $r_J(\q)   $    instead of  
$r_J(\mathbb{M}).    $  In the $1$-dimensional case
$$  r_J(\mathbb{M})\le e_0(\mathbb{M}) -1.$$
This bound can be easily extended to higher dimension under the assumption    $ {\depth} \ gr_{\mathbb{M}}(M)\ge r-1. $  Moreover, as in the classical  case (see  \cite{T}) and \cite{S7}),    if $ M$ is Cohen-Macaulay and $  {\depth} \ gr_{\mathbb{M}}(M)\ge r-1, $  then $ r_J(\mathbb{M})$ is independent  of $J.  $   

\noindent In Theorem \ref{r(I)} we prove that if $\dim M$=$2 $ or, more in general,  if we assume $ \depth  \ gr_{\mathbb{M}}(M)\ge r-2, $ then 
$$  r_J(\mathbb{M})\le e_1(\mathbb{M})- e_0(\mathbb{M}) + h_0(\mathbb{M})+1.$$
By using this bound,    as a bonus,  we get  easy proofs of new  ``border cases"   theorems.  In particular we consider the filtrations  having $ e_1(\mathbb{M})=e_0(\mathbb{M})-h_0(\mathbb{M})+1$ or  $ e_1(\mathbb{M})=e_0(\mathbb{M})-h_0(\mathbb{M})+2. $  
  
\noindent  If  $e_1(\mathbb{M})=e_0(\mathbb{M})-h_0(\mathbb{M})+1,$ we have $$P_{\mathbb{M}}(z)=\frac{h_0(\mathbb{M})+h_1(\mathbb{M})z+z^2}{(1-z)^r}
 $$ if $M_2\cap JM=JM_1$ (Theorem \ref{nor+1}) or $e_2(\mathbb{M})\not= 0$ (Theorem \ref{sal}).   
In this  case $\mathbb{M} $ has  almost minimal multiplicity. The case when $M_2\cap JM\not=JM_1$ and 
 $e_2(\mathbb{M})= 0 $ is more difficult and, by using the usual approach,  
 we reprove in dimension two a nice  result due to   Goto,  Nishida and  Ozeki    (see \cite{GNO2}).

 \vskip 4mm  
 \section{A bound on the reduction number}
 
 Let $\q$ be an  $\m$-primary ideal of $A$ and  let  $\mathbb{M} = \{M_j \}_{ j\ge 0}$  be a  { good   $\frak{q}$-filtration} of a  finitely generated module
$M.$ 
  

We start with   a result on the reduction number which was  proved in \cite{R1},
Theorem 1.3. in the case of the $\q$-adic filtration on $A.$ As usual we denote by $J$ an ideal generated by a maximal  $\mathbb{M}$-superficial sequence for $\q.$       
 
 Assume $\depth_{\q} M \ge 1, $ by Lemma \ref{rr}, there exists an integer $p$ such that $M_j =\w{M}_j $ for every $j \ge p.$  This implies that  
 $$ \w{M}_{j+1} =J \w{M}_j +M_{j+1} $$  for every $j \ge p.$ Hence the module
 $$N:= \oplus_{j\ge 0}( \w {M}_{j+1}/ J\w{M}_j+M_{j+1})$$ has finite length. In the following we denote 
   by $\nu $ the    minimal number  of generators of $ N $ and remark that
  $$\nu =\dim_k N/\m N <  \lambda (N).$$

 \begin{theorem}   Let  $\mathbb{M}$ be   a good $\q$-filtration of  $M $ and  let $J$  be an ideal generated by a maximal  $\mathbb{M}$-superficial  sequence  
for $\q. $ Then 
 $$ \q^{\nu} \subseteq  J\q^{\nu -1} + (M_{\nu +n} : \w{M}_n )  $$
 for every positive integer $n.$
 \end{theorem}
 \begin{proof}  Let $p$ be the integer such that $  \w{M}_n  = J \w{M}_{n-1}  + M_n $ for all $n>p.$ For all $ n=1,\dots,
p$ we consider the elements $m_{1n},\dots, m_{\nu_n n} \ \in \w{M}_n $ such that the corresponding elements in $ N_n=  \w {M}_n/  J \w{M}_{n-1} +M_n $ form a  minimal  system of generators as an $ A-$module. In particular if we define  the submodules 
 $$ L_n:= <  m_{1n},\dots, m_{\nu_n n} > \subseteq \w{M}_n,$$
 then $  \w{M}_n = J \w{M}_{n-1} + L_n + M_n $ (with  $L_n=0 $ if $n>p$). It is easy to see that 
for every $n \ge 1 $ we can write 
$$  \w{M}_n = \sum_{j=0}^n J^{n-j} L_j + M_n. $$
  
\noindent   We have $\nu = \sum_{n=1}^p
\nu_n$ and so $|(in)|=\nu$ if  $i=1,\dots,\nu_n $ and $n=1,\dots,p.$  
 
Denote by $a_{in} $  an element of $ \q, $ then $a_{in} m_{in} \in \w {M}_{n+1}=\sum_{j=0}^{n+1} J^{n+1-j} L_j + M_{n+1}.$  Then 
  there exist $ c_{(in)(kj)}\in J^{n+1-j}$ and $\alpha _{i n+1} \in M_{n+1} $
such that
$$  a_{in}  m_{in}
=\sum_{j=1}^{n+1}\sum_{k=1}^{\nu_j}c_{(in)(kj)}  m_{kj} + \alpha _{i n+1} $$ with $m_{k\ p+1}=0$ for every $k.$

Thus if we consider the relations
$$\sum_{j=1}^{n+1}\sum_{k=1}^{\nu_j}c_{(in)(kj)}  m_{kj} -  a_{in}   m_{in} = \alpha _{i n+1} $$ we get a
system of $\nu$ linear equations in the $\nu$ variables $ {m_{kj}}$ where $j=1,\dots,p$ and
$k=1,\dots,\nu_j.$ 
  
Hence      we consider the linear system
  $$ B  \begin{pmatrix}    m_{1 1}  \\  \dots \\   m_{\nu_p p} \end{pmatrix} =   \begin{pmatrix}   \alpha_{1 2} \\  \dots \\  \alpha_{\nu_p p+1}  \end{pmatrix}  $$
 where  $B=(b_{(in)(kj)})  $ is the   matrix of the coefficients of the variables with
 $(in)$ and $(kj) $ running through $  \{ (1 1),\dots, (\nu_1 1), \dots \dots, (1 p), \dots, (\nu_p p) \}.$
 The matrix  $B$    has entries in $A$ and it has size $\nu \times \nu. $ We remark that  the  ${(in)(kj)}$-entry  
   is an element   in $ \q^{ n+1-j}$ which is zero  if $n+1< j. $   In particular, on the diagonal $(in)=(kj), $ we have    $ b_{(in)(in)} =c_{(in)(in)}-a_{in} \in  \q. $  

\noindent We remark now  that the determinant of $B$ is   an element of  $ \q^{\nu}, $ in particular
we can prove that  all  the cofactors of $B$  are   elements of  $ \q^{\nu}. $    
Indeed  a cofactor $ \prod b_{(in)(kj)} $ is a product of $\nu$ elements of $B$ involving one entry for each row $(in)$ and all the $\nu $ columns  $(kj) $ are reached.  Due to the fact that $b_{(in)(kj)} \in  \q^{ n+1-j}$  for
every  $n,j \in \{1,\dots,p\}$ and  $ i=1,\dots,\nu_n,\   k=1,\dots,\nu_j,$  it is easy to see that   
$$ \prod b_{(in)(kj)} \in \q^{\nu_1 +2 \nu_2 + \dots + p \nu_p + \nu -(\nu_1 +2 \nu_2 + \dots + p \nu_p)} = \q^{\nu}.$$

\noindent  Let $a= \prod a_{in} $ be where the product is over  $n=1,\dots,p$ and $i=1,\dots,\nu_n.$ \par By expanding out the determinant of $B$ and considering  the fact that for $n > j$ we have $b_{(in)(kj)} \in J \q^{n-j},$    
one can  see   that
$$\det(B)=(-1)^\nu ( a -\sigma)  $$ for a suitable $\sigma \in J
\q ^{\nu -1}.$ 

Now the Cayley-Hamilton theorem  asserts that
$$ \det(B) \begin{pmatrix}   m_{1 1} \\  \dots \\   m_{\nu_p p} \end{pmatrix} =  B^* \begin{pmatrix}   \alpha_{1 2} \\  \dots \\  \alpha_{\nu_p p+1}  \end{pmatrix} $$ 
where 
  $B^*$ denotes the adjoint matrix of $B.$ We claim $$\det(B) m_{kj} \in M_{\nu+j} $$
for every $j=1,\dots,p$ and
$k=1,\dots,\nu_j. $ 

We remark that  $B^* $  has size  $\nu  \times \nu  $ and, since   the ${(in)(kj)}$-entry  of $B$ 
   is an element   in $ \q^{ n+1-j}, $  as before in the computation of $\det(B), $    we can prove that the  $ (kj)(in)$-entry of $B^*, $ is an element  of  $ \q^{\nu -( n+1-j)}. $ Hence the product of the $(kj)$-row of $B^*$ and $
\begin{pmatrix}   \alpha_{1 2} \\  \dots \\  \alpha_{k j+1} \\ \dots \\ \alpha_{\nu_p p+1}  \end{pmatrix} $ is  an element of $\q^{\nu -( n+1-j)} M_{j+1} \subseteq M_{\nu +j}, $ as required.
 
 Since for every $n \ge 1 $ we can write 
 $  \w{M}_n = \sum_{j=0}^n J^{n-j} L_j + M_n  $    and $\det(B) L_j \subseteq M_{\nu +j}, $   we easily get $$\det (B)=  a -\sigma \in M_{n+\nu} :  \w{M}_n.$$  

We may repeat the same procedure for all monomial $a= \prod a_{in} $ in $\q^{\nu}$ and
the result follows.
 \end{proof}

 Since $M_{n+1}=JM_n =J \w{M_n} $ for large $n,$ we may define the integer
$$ k:=  \text{min} \{ t \ | \  M_{t+1} \subseteq J \w{M}_t \}. $$

 \begin{corollary}\label{nu}Let  $\mathbb{M}$ be   a good $\q$-filtration on $M $ and 
   let $J$  be an ideal generated by a maximal  $\mathbb{M}$-superficial  sequence  
for $\q. $  Then
$$\q^{\nu} M_{k+1} \subseteq J M_{k+\nu}.$$
\end{corollary}

\begin{proof} From the above theorem we get
$$\q^{\nu} M_{k+1}\subseteq (J\q^{\nu -1}+M_{\nu +k}:_A\w{M}_k)M_{k+1} \subseteq JM_{k+\nu} + (M_{\nu +k}:_A\w{M}_k)J\w{M}_{k } \subseteq J M_{k+\nu}.$$ 
  \end{proof}
  
 \begin{corollary}\label{r} Let    $L$ be a submodule of a   module $M$ and let $\mathbb{M}=\mathbb{M}_L$ be the good $\q$-filtration on $M$ induced by $L.$  Let  $J$  be  an ideal generated by a maximal  $\mathbb{M}$-superficial  sequence 
for $\q.$  With the above notations we have $$  r_J(\mathbb{M}) \le k+\nu.$$
 \end{corollary}
 \begin{proof}   We remark that,  by the definition of $\mathbb{M}_L, $ we have $M_{k+\nu +1} =\q^{\nu}M_{k+1}.$
 Then the result  follows by Corollary  \ref{nu}.   \end{proof}
 
 In the following  we denote by $$ S_J := \{ n \in {\mathbb N} \  | \   M_{j+1}\cap J\w{M}_j=JM_j \ \ \text{for every}  \ j,\ \  0\le j \le n\}.$$   We remark that $ S_J \neq \emptyset  $ since $0 \in S_J.$
 
 \vskip 2mm
 The following result extends  Theorem 1.3. in \cite{R1}  to modules. We include here a proof even if it is essentially a natural recasting 
   of the original  result proved in the classical case.  
 
\begin{theorem}
\label{1.3} Let  $L$ be a submodule of the   module $M$ and let $\mathbb{M}=\mathbb{M}_L$ be the good $\q$-filtration on $M$ induced by $L.$ Let   $J$  be  an ideal generated by a maximal  $\mathbb{M}$-superficial  sequence 
for $\q.$  If    $n \in S_J, $ then 
$$  r_J(\mathbb{M}) \le \sum_{i\ge 0} v_i(\w{\mathbb{M}}) + n+1 -  \sum_{i=0}^n v_i({\mathbb{M}}) $$ 
\end{theorem}
\begin{proof} By Corollary \ref{r} we have 
$$  r_J(\mathbb{M}) \le \nu + k= \sum_{i\ge 0} \nu_i +k \le \sum_{i\ge 0} \la(\w{M}_{i+1 }/J \w{M}_{i }+M_{i+1}) +k. $$
But it is clear that $$\la(\w{M}_{i+1 }/J \w{M}_{i }+M_{i+1})=v_i(\w{\mathbb{M}})-\la(M_{i+1}/J\w{M}_i\cap M_{i+1})\le v_i(\w{\mathbb{M}}),$$ so that 
\begin{equation*} \la(\w{M}_{i+1 }/J \w{M}_{i }+M_{i+1}) =v_i(\w{\mathbb{M}})-v_i({\mathbb{M}})   \ \ \  \text{if $0\le i \le n$} \end{equation*}
and $ \la(\w{M}_{i+1 }/J \w{M}_{i }+M_{i+1}) < v_i(\w{\mathbb{M}})\ \ \   \text{if $0\le i\le k-1.$}$ 
This implies 
$$ r_J(\mathbb{M}) \le \sum_{i=0}^n(v_i(\w{\mathbb{M}})-v_i({\mathbb{M}}))+\sum_{i\ge n+1}v_i(\w{\mathbb{M}})+k$$ so that the conclusion follows if $k\le n+1.$ If $k\ge n+2, $ then we have

\begin{equation*} \begin{split} r_J(I) & \le \sum_{i=0}^n ( v_i(\w{\mathbb{M}}) -  v_i( {\mathbb{M}}))  + \sum_{i=n+1}^{k-1} ( v_i(\w{\mathbb{M}})  - 1) +\sum_{i\ge k}  v_i(\w{\mathbb{M}})  +\ k \\  &=\sum_{i\ge 0} v_i(\w{\mathbb{M}}) +n +1 - \sum_{i= 0}^n v_i( {\mathbb{M}}).\\ \end{split}\end{equation*}

\end{proof} 
 
 A nice application of the above theorem is given by  Kinoshita,  Nishida,  Sakata,  Shinya in \cite{N}.
 \vskip 2mm
In some cases  we can get rid  of the $v_i$'s  in the  formula given in Theorem \ref{1.3}
.  
  \begin{theorem} \label{r(I)}  Let  $L$ be a submodule of the  $r$-dimensional Cohen-Macaulay module $M$ and let $\mathbb{M}=\mathbb{M}_L$ be the good $\q$-filtration on $M$ induced by $L.$  If    $ \depth  \ gr_{\mathbb{M}}(M)\ge r-2,  $   then
  $$ r_J(\mathbb{M})  \le e_1(\mathbb{M}) -e_0(\mathbb{M}) + h_0(\mathbb{M})+1 $$
  for every ideal $J$ generated by a maximal $\mathbb{M}$-superficial sequence for $\q. $
  \end{theorem}
 \begin{proof}  First of all we prove that we may reduce the problem to dimension $r \le 2.$ Let $r > 2$ and $J=(a_1,\dots,a_r)$   be  an ideal generated by a maximal  $\mathbb{M}$-superficial  sequence 
for $\q$ ; we consider the ideal   $K:=(a_1,\dots,a_{r-2}).   $    Since $\depth  \ gr_{\mathbb{M}}(M)\ge r-2, $ by Lemma \ref{sup} and the Valabrega-Valla criterion, we get $M_{n+1} \cap KM=KM_{n} $ for every $n,$ which easily implies $$v_j( \mathbb{M})=v_j(\mathbb{M}/KM)$$ for every $j.$ 
 Hence  $$r_J(\mathbb{M})= r_{\frak{a}}  (\mathbb{M}/KM)$$ where $ \frak{a}=(a_{r-1},a_r).$  Moreover, $e_1(\mathbb{M})=e_1(\mathbb{M}/KM),  $ $e_0(\mathbb{M})=e_0(\mathbb{M}/KM)  $ 
  and $ h_0(\mathbb{M})=h_0(\mathbb{M}/KM) $ since $KM \subseteq L$ by definition. 
  
  \noindent   Hence we may assume $r \le 2;$   by Theorem \ref{1.3} with $n=0,$  we have  
  $$  r_J(\mathbb{M}) \le \sum_{i\ge 0} v_i(\w{\mathbb{M}})  \ +1\ - \   v_0({\mathbb{M}})=\sum_{i\ge 0} v_i(\w{\mathbb{M}})  \ +1-e_0({\mathbb{M}})+h_0({\mathbb{M}}). $$ Since $\depth  \ gr_{\w{\mathbb{M}}}(M)\ge 1\ge r-1, $ by Theorem \ref{e1=CM}   we have
   $   \sum_{i\ge 0} v_i(\w{\mathbb{M}}) = e_1(\w{\mathbb{M}})= e_1(\mathbb{M}).$ The conclusion follows.
  \end{proof}

 \vskip 2mm  
 
 \begin{remark} \label{corollaries}{\rm{We note  that the above bound  is sharp and it gives easy proofs of some results presented in Chapter 2.   For example if we consider   $e_1(\mathbb{M}) = e_0(\mathbb{M}) -  h_0(\mathbb{M}), $ i.e. the minimum value of $e_1(\mathbb{M})  $ with respect to Northcott's bound, we immediately get  $  r_J(\mathbb{M}) \le 1, $ hence
 $M_2=JM_1 $ and obviously $gr_{\mathbb{M}}(M) $ is Cohen-Macaulay (see Theorem \ref{nor}).  
 
 If $e_1(\mathbb{M}) = e_0(\mathbb{M}) -  h_0(\mathbb{M}) +1,  $ then  $  r_J(\mathbb{M}) \le 2, $ hence
 $M_3=JM_2 $ and, if $M_2 \cap J= JM_1, $  by Valabrega-Valla criterion,  we conclude that
 $gr_{\mathbb{M}}(M) $ is Cohen-Macaulay (see Theorem \ref{nor+1}). 
 }}
 \end{remark}

  \vskip 4mm  
 \section{A generalization of Sally's conjecture}
 \vskip 2mm
As a consequence of Theorem \ref{r(I)},  we present now an extended version of the Sally conjecture.  We recall that in the case of the $\m$-adic filtration, the question raised by Sally was the following:  if  $A$ is Cohen-Macaulay of  dimension $r$ and  it has almost minimal multiplicity, that is $e_0(\m)= \mu(\m) -r +2, $  is it  true that $ \depth gr_\m(A) \ge r-1?$ 

Sally's conjecture  was proved by \cite{RV1} and independently  by \cite{W1}. The next  theorem is  a generalization of this  result.

Let $\mathbb{M}$ be a  good $\q$-filtration of $M$ and $J$ an ideal generated by a maximal  sequence of $\mathbb{M}$-superficial elements for $\q.$ 
For every integer $j\ge 0$, we have defined the integers
$$\begin{aligned} v_j( \mathbb{M}) &=\la(M_{j+1}/JM_j)\\ 
vv_j( \mathbb{M}) &=\la(M_{j+1}\cap JM/JM_j)\\
w_j( \mathbb{M})&=\la(M_{j+1}/M_{j+1}\cap JM).
\end{aligned}$$ so that we have the formula 
$$v_j( \mathbb{M})=w_j( \mathbb{M})+  vv_j( \mathbb{M}).$$ Further, 
if $x \in J$ is a superficial element, it is easy to see that $$\begin{aligned} vv_j( \mathbb{M})&\ge vv_j( \mathbb{M}/xM)\\  w_j( \mathbb{M})&=w_j( \mathbb{M}/xM)\\
   v_j( \mathbb{M})&\ge v_j( \mathbb{M}/xM)\end{aligned}
   $$ for every $j\ge 0.$
 \vskip 2mm
 We recall that if $vv_j( \mathbb{M})=0$ for every $j\ge  0, $ then by the Valabrega-Valla criterion, $gr_{\mathbb{M}}(M) $ is Cohen-Macaulay. 
 
  \begin{theorem} \label{p} Let  $L$ be a submodule of the  $r$-dimensional Cohen-Macaulay module $M$ and let $\mathbb{M}= \mathbb{M}_L $  be the good $\q$-filtration on $M$ induced by $L.$ Let   $J$  be  an ideal generated by a maximal  $\mathbb{M}$-superficial  sequence 
for $\q.$  
  If  there exists a positive integer $p$ such that
  
  1.    $vv_j( \mathbb{M})=0$ for every $j\le p-1$
  
 \noindent  and 
  
  2.  $v_p( \mathbb{M}) \le 1,$  
  
 \vskip 1mm  \noindent then     
   $\depth  \ gr_{\mathbb{M}}(M)\ge r-1. $
\end{theorem}
 
\begin{proof}  Conditions 1. and 2. are preserved modulo  superficial elements in $J,$ so that, by using the Sally machine, we may reduce the problem to the case $ r=2 $ and  we have to prove that $ {\depth} \ gr_{\mathbb{M}}(M)\ge 1. $    
  
\noindent   We may assume that  $v_p( \mathbb{M})=\lambda (M_{p+1}/JM_p) =1, $ otherwise, by the Valabrega-Valla criterion, we  immediately get that $gr_{\mathbb{M}}(M) $ is  Cohen-Macaulay.
 
\noindent    Since $M_{p+1}=\q M_p$ is generated over $A$ by the products $am$ with $a\in  \q$ and $m\in M_p,$ the condition $\la (M_{p+1}/JM_p) =1$  implies $M_{p+1}=JM_p+(a)m$ with 
   $a \in \q,$ $m\in M_p$  and $am\notin JM_p.$  Then  for every $n\ge p$ the multiplication by $a$  gives  a surjective map  from
$M_{n+1}/JM_n$ to $M_{n+2}/JM_{n+1};$    this implies
$$v_n(\mathbb{M})=\la(M_{n+1}/JM_n) \le 1$$ for every $ n\ge p.$ 

Let $J=(x,y)$ and $s:=r_{(y)}(\mathbb{M}/xM)$ the reduction number of $\mathbb{M}/xM$ with respect to the ideal $(y).$ This means that $$  v_j(\mathbb{M}/xM) = 0 \ \ \ \text{if} \ \ j\ge s, \ \ \ v_j(\mathbb{M}/xM)  > 0 \ \ \ \text {if}\ \  j< s.$$ It follows that, when  $s\le p,$ we get $vv_j(\mathbb{M}/xM)=0$ for every $j\ge 0,$ so that $gr_{\mathbb{M}/xM }(M/xM) $ is Cohen-Macaulay. By Sally's machine, this implies $gr_{\mathbb{M}}(M) $ is  Cohen-Macaulay as well.
   
  Hence we may assume  $s >p \ge 1$ and we  prove  that
\begin{equation} \label{claim} {v_j( \mathbb{M})= v_j( \mathbb{M}/xM)} \end{equation}   for $j=0,\dots,s-1$. If  $0\le j \le p-1,$ we have $vv_j(\mathbb{M})=0$ by assumption, so that $$v_j( \mathbb{M})=w_j(\mathbb{M})=w_j(\mathbb{M}/xM)\le v_j(\mathbb{M}/xM)\le v_j( \mathbb{M}).$$ On the other hand, if $p\le j \le s-1,$ we have $$0<  v_j(\mathbb{M}/xM)\le v_j( \mathbb{M})\le 1.$$ This proves (\ref{claim}) and also that $${  v_j( \mathbb{M})= v_j( \mathbb{M}/xM)=1}$$ for all $p\le j \le s-1.$

\noindent Now,  
    for every
$j\ge 0,$ we have $v_j(\mathbb{M}/xM)=\la(M_{j+1}/JM_j+x(M_{j+1}:x)),$ and hence  for ${\  j=0,\dots,s-1} $
$$ v_j(\mathbb{M})=v_j(\mathbb{M}/xM)=\la(M_{j+1}/JM_j+x(M_{j+1}:x)). $$ 
From the above equality and  the following  exact sequence
{\small{$$ 0\longrightarrow M_j:x/M_j:J \overset{y} \longrightarrow M_{j+1}:x/M_j \overset{x} \longrightarrow M_{j+1}/JM_j\longrightarrow
M_{j+1}/JM_j+x(M_{j+1}:x)\longrightarrow 0 $$}}we get $$\la(M_j:x/M_j:J)=\la(M_{j+1}:x/M_j)$$ for every $j=0,\dots,s-1.$ With $j=1$ this gives $M_2:x=M_1$, so that,
by induction on $j,$ we get ${  M_{j+1}:x=M_j}$ for ${\  j=0,\dots,s-1}.$

\noindent   We claim that   if   $ r_J(\mathbb{M})\le s, $ then ${  M_{j+1}:x=M_j }$
for every $j\ge 0$ and so depth $ gr_{\mathbb{M}}(M) > 0, $ as required. In fact, if we have $ r_J(\mathbb{M})\le s,$ then $M_{t+1}=JM_t$ for every $t\ge s.$ Let us  assume  by induction that $j\ge s$ and $M_j:x=M_{j-1} $ (we know that ${  M_{s}:x=M_{s-1}}$), then we get $$M_{j+1}:x=JM_j:x=(xM_j+yM_j):x\subseteq M_j+y(M_j:x)=M_j+yM_{j-1}=M_j$$ which proves the claim.
    
\noindent     It remains  to prove that $ r_J(\mathbb{M})\le s.$ We have $$e_1({\mathbb{M} })= e_1(\mathbb{M}/ xM) =\sum_{j\ge 0} v_j(\mathbb{M}/xM)= \sum_{j\le p-1} v_j( {\mathbb{M} }) + s-p.$$  Further, since $vv_j(\mathbb{M})=0$ for every $j\le p-1,$ we have $p-1\in S_J$, so that, by Theorem \ref{1.3}, we get  $$  r_J(\mathbb{M}) \le \sum_{j\ge 0} v_j(\w{\mathbb{M}}) +\ p \ - \ \sum_{j=0}^{p-1}  v_j({\mathbb{M}})= e_1(\mathbb{M}) + \ p \ - \ \sum_{j=0}^{p-1} v_j({\mathbb{M}})= s. $$ 
\end{proof}

\begin{corollary}   \label{sp}
With the  notation  of Theorem \ref{p}, assume 
there exists a positive integer $p$ such that
  
  1.    $vv_j( \mathbb{M})=0$ for every $j\le p-1$

  2.  $v_p( \mathbb{M})= 1.$   
  
\noindent   Then   $$P_ {\mathbb{M}}(z)= {{\  \sum_{n=0}^{p-1} \lambda (M_n/M_{n+1} + JM_{n-1})z^n \ + \  (\lambda(M_p/JM_{p-1})-1\  )z^p \ + \       z^s} \over (1-z)^r}$$ 
for some $s, $   $p+1 \le s \le e_0(\mathbb{M})-1 .$   

\noindent Furthermore if $M_{p+1}  \cap J =J M_{p},  $ then
   $\ gr_{\mathbb{M}}(M)$ is Cohen-Macaulay if and only if   $s=p+1. $ 
   \end{corollary}
   
   \begin{proof} Now, by assumption,  $M_{n} \cap JM=JM_{n-1}  $ for every $n \le p, $ then for $n <p $ we  have $v_{n-1}(\mathbb{M})-v_n(\mathbb{M}) =\lambda ( M_{n}+JM /JM)-
 \lambda ( M_{n+1}+JM /JM) 
 =   \lambda ( M_n +JM / M_{n+1} + JM )=\lambda (M_n/(M_{n+1} + JM_{n-1})).$  Further by using the information coming from the   proof of above theorem,  if   $p<s, $  then $ v_p(\mathbb{M})=\dots=v_{s-1}(\mathbb{M})=1.$ The Hilbert series follows  by Theorem  \ref{e1=CM} (5.)  and   by Theorem  \ref{p}.  
  
 It is clear that,  if $\ gr_{\mathbb{M}}(M)$ is Cohen-Macaulay,  then  $s=p+1. $ Conversely if  $s=p+1  $ and $M_{p+1}  \cap J =J M_{p},  $ we can prove that the $h$-polynomial of ${\mathbb{M}}$ coincides with that of ${\mathbb{M}}/JM$ and hence  $\ gr_{\mathbb{M}}(M)$ is Cohen-Macaulay. In fact $h_n({\mathbb{M}})=h_n({\mathbb{M}}/JM)=   \lambda (M_n/M_{n+1} + JM_{n-1}) $ for $n \le p-1. $ Further $h_p({\mathbb{M}}/JM)= \la(M_p+JM/ M_{p+1} +JM)=\la(M_p/M_{p+1} +J M_{p-1})= \la(M_p/JM_{p-1})-\la(M_{p+1}/M_{p+1} \cap J)= \la(M_p/JM_{p-1})-v_p({\mathbb{M}})= h_p({\mathbb{M}}).$ Finally $ 1=h_{p+1}({\mathbb{M} })= h_{p+1}({\mathbb{M}}/JM) $ since $ e_0({\mathbb{M}})=e_0({\mathbb{M}}/JM)= \sum_{i\ge 0}h_i({\mathbb{M}}/JM).$ 
 
     \end{proof}

\vskip 2mm


    
The assumptions of the above result are satisfied if we consider the $\m$-adic filtration on a local Cohen-Macaulay ring of initial degree $p-1. $ Hence Corollary \ref{sp} extends Theorem 3.1. \cite{RV0}.

\vskip 3mm
The next result is the promised  extension to modules of the classical Sally  conjecture. We point out that    the statements a) and b)  are   independent of $J. $ As in the following, results  where the effective role of the superficial sequences is limited to the proof and it disappears in the  statements,  will  be   highly appreciated. 

\begin{corollary}{\label{Sallymodule}}
\label{1.6}  Let  $L$ be a submodule of the  $r$-dimensional Cohen-Macaulay module $M$ and let $\mathbb{M}=\mathbb{M}_L$ be the good $\q$-filtration on $M$ induced by $L.$  The following conditions are equivalent:

a)  $e_0(\mathbb{M})= h_0(\mathbb{M})+ h_1(\mathbb{M})  + 1$ 

\vskip 2mm 
b) $P_{\mathbb{M}}(z)=\frac{ h_0(\mathbb{M})+ h_1(\mathbb{M})z+z^s}{(1-z)^r}$ for some integer $s\ge 2.$
\vskip 2mm 
\noindent Further, if either of the above conditions holds, then we have

 \vskip 2mm 
 c)  $ \depth  \ gr_{\mathbb{M}}(M) \ge r-1.$
 
 \vskip 2mm
 d) Let $J$ be the  ideal generated by a maximal $\mathbb{M}$-superficial sequence and assume  $M_2 \cap J= J M_1.$ Then  $gr_{\mathbb{M}}(M)$ is Cohen-Macaulay $\Longleftrightarrow$ $s=2$ 
 $\Longleftrightarrow$ $e_1(\mathbb{M})=e_0(\mathbb{M})-h_0(\mathbb{M})+1$  $\Longleftrightarrow$
 $e_1(\mathbb{M})=h_1(\mathbb{M})+2.$
\end{corollary}
\begin{proof} It is clear that b) implies a). By the Abhyankar-Valla  formula, if we have $e_0(\mathbb{M})=h_0(\mathbb{M})+h_1(\mathbb{M}) + 1,$ then 
$\la(M_2/JM_1)=1 $  for every maximal superficial sequence $J;$ hence  a) implies b) 
by the above corollary,  and a) implies c) by Theorem  \ref{p} with $p=1.$  Finally the first equivalence in d) follows by Corollary \ref{sp} and the remaining part is a trivial computation.  \end{proof}

  We remark that    the assumption $M_{2}  \cap J =J M_{1}  $ in Corollary \ref{Sallymodule}  d) is necessary.    In fact if   we consider Example \ref{huno}, the $\q$-adic filtration   has almost minimal multiplicity 
 since  $\la(\q^2/t^4 \q)=1, $   $P_{\q}(z)= {{ 2 +z+z^2}\over (1-z)}, $  but the associated graded ring  $\ gr_{\q}(A)$ is not Cohen-Macaulay.
 
 \vskip 2mm
By using the results  of this chapter we easily a   collection of results  proved in \cite{H}, \cite{O},
\cite{HM2}, \cite{I}
and \cite{GR1}
   and already discussed  here  (see Theorem \ref{nor}, Theorem \ref{nor+1}, Theorem \ref{nor+2}) by using easier methods. 
See also  Remark \ref{corollaries} and Corollary  1.9.  \cite{R1}.

\section{The case  $e_1(\mathbb{M})=e_0(\mathbb{M})-h_0(\mathbb{M})+1$ }

 \vskip2mm
 
We prove now some results on $ e_1(\mathbb{M}) $ for which we cannot avoid the hard   machinery introduced in this chapter.   
\vskip 2mm
The next result   completes  Theorem \ref{sal}. In  the classical case, in  \cite{Goto}  Theorem 1.2,    a different proof of the same result had been presented  which involved  the structure of the Sally module.   
\vskip 2mm

 \begin{theorem} \label{goto2} Let $M$ be a     Cohen-Macaulay module of dimension two, $L$ a submodule of $M$ and let $\mathbb{M}=\mathbb{M}_L$ be the good $\q$-filtration on $M$ induced by $L.$     Then  $e_1(\mathbb{M})=e_0(\mathbb{M})-h_0(\mathbb{M})+1$  if and only if either   $$P_{\mathbb{M}}(z)=\frac{h_0(\mathbb{M})+h_1(\mathbb{M})z+ z^2}{(1-z)^2} $$ 
or  
 $$P_{\mathbb{M}}(z)=\frac{h_0(\mathbb{M})+h_1(\mathbb{M})z+3 z^2-z^3}{(1-z)^2}.$$  In the first case $\depth  \ gr_{\mathbb{M}}(M) > 0, $ while in the second   $\depth    gr_{\mathbb{M}}(M)=~ 0. $
 
 \end{theorem}
 \begin{proof} First we remark that if  $M$ is  a     Cohen-Macaulay module of dimension one and $ e_1(\mathbb{M})=e_0(\mathbb{M})-h_0(\mathbb{M})+1, $ then $$P_{\mathbb{M}}(z)=\frac{h_0(\mathbb{M})+( e_0(\mathbb{M})-h_0(\mathbb{M}) -1 )z+ z^2}{(1-z) }.$$ In fact,  by Theorem \ref{e1=CM}, we have $$ e_1(\mathbb{M})= \sum_{i \ge 0} v_i(\mathbb{M}) = e_0(\mathbb{M})-h_0(\mathbb{M})+\sum_{i\ge 1} v_i(\mathbb{M}) = e_0(\mathbb{M})-h_0(\mathbb{M})+1. $$ Necessarily  $v_1=1 $ and $v_i =0 $ for every $i \ge 2, $ hence  the Hilbert  series follows. In particular we remark that $e_2=1.$
 
 Let now $M$ be a     Cohen-Macaulay module of dimension two and assume  $ e_1(\mathbb{M})=e_0(\mathbb{M})-h_0(\mathbb{M})+1  $; if $e_2(\mathbb{M}) \not = 0, $ then the result follows by Theorem \ref{sal}.  Hence we assume
 $e_2 (\mathbb{M})=0.$ Since $\depth \ gr_{\widetilde{\mathbb{M}}}(M)\ge 1 $,  we have $e_2(\mathbb{M})=e_2(\widetilde{\mathbb{M}})=\sum_{j\ge 1}jv_j(\widetilde{\mathbb{M}})$ 
and $e_1(\mathbb{M})=e_1(\widetilde{\mathbb{M}})=\sum_{j\ge 0}v_j(\widetilde{\mathbb{M}}).$ It follows  
$v_j(\widetilde{\mathbb{M}})=0 $ for $j \ge 1 $ and $e_1(\mathbb{M})=e_1(\widetilde{\mathbb{M}})=v_0(\widetilde{\mathbb{M}})= e_0(\widetilde{\mathbb{M}}) - h_0(\widetilde{\mathbb{M}}). $   By Theorem \ref{nor} it follows that  $gr_{\widetilde{\mathbb{M}}}(M)$ is Cohen-Macaulay  and 
\begin{equation} \label{Ptilde}   P_{\w{\mathbb{M}}}(z)=\frac{h_0({\w{ \mathbb{M}}})+( e_0(\mathbb{M})-h_0(\w{\mathbb{M}}) )z }{(1-z)^2}.
\end{equation}
 
\noindent  Since $e_1(\mathbb{M})=  v_0(\widetilde{\mathbb{M}})= v_0( {\mathbb{M}}) + \la(\w{M}_1/M_1)=e_0(\mathbb{M})-h_0(\mathbb{M})+ \la(\w{M}_1 /M_1), $   then $ \la(\w{M}_1/M_1)=1.$ 
 We prove   now that
 $\w{M_i}= M_i $ for every $i \ge 2.$
  
\noindent First we remark that $v_1(\mathbb{M}) >1. $ In fact if $v_1(\mathbb{M})= 1,  $ then   $e_0(\mathbb{M})= h_0(\mathbb{M})+ h_1(\mathbb{M})  + 1$ and  by Corollary \ref{Sallymodule}, we would have $e_2 \not = 0.$ Let  $x $ be a $\mathbb{M}$-superficial element and write $\overline {\mathbb{M}} = \mathbb{M}/x M.$ Now  $ M_2 : x \not = M_1$ since 
$ v_1(\mathbb{M}) > 1=v_1(\overline{ \mathbb{M}}). $

\noindent  Moreover we know that  $e_2(\mathbb{M})=0 $ and 
$ e_2(\overline {\mathbb{M}})=1 $ hence,  by Proposition \ref{ei}, 4), we get $  \sum_{i \ge 0} \lambda(M_{i+1}  : x / M_i)= 1$ and then  $ M_{i+1} : x = M_i $ for every $i \ge 2.$ As a consequence  it follows  that $\w{M}_i= M_i $ for every $i \ge 2,  $ as wanted. 
  Since  $$P_{\mathbb{M}}(z)= P_{\w{\mathbb{M}}}(z) + \sum_{i\ge 0} [ \la(\w{M}_{i+1}/ M_{i +1}) - \la(\w{M}_i/ M_i)] z^i, $$
 by (\ref{Ptilde})   and the previous fact,  we get
$$P_{\mathbb{M}}(z)  = \frac{h_0({\w{ \mathbb{M}}})+( e_0(\mathbb{M})-h_0(\w{\mathbb{M}}) )z }{(1-z)^2}
 + (1-z) =$$ $$= \frac{h_0(\mathbb{M})+( e_0(\mathbb{M})-h_0(\mathbb{M}) -2 )z+3 z^2-z^3}{(1-z)^2}.$$
 
\end{proof}
\vskip 3mm

 Example \ref{1} and Example \ref{2} show that in  Theorem \ref{goto2} both the Hilbert series can occur. 
 We remark that Theorem \ref{goto2}  cannot be extended to dimension $\ge 3.$  
 In fact if we consider in $R=k[[x,y,z]]$ the $\q$-adic filtration with  $\q =(x^2-y^2, x^2-z^2,xy,xz,yz)$, then 
$h_0(\q)= \la(R/\q)=5, $ $e_0(\q)=8, $ $e_1(\q)=4 = e_0(\q) - h_0(\q) +1, $ but 
$$ P_{\q}(z) =\frac {5 + 6z^2 - 4z^3 + z^4}{ (1-z)^3.}$$
In this case $\text{depth} \ gr_{\q}(R) =0    $  because $x^2 \not \in \q $ and $x^2 \in \q^2: \q.$

In the last section we will   characterize all the possible Hilbert series in higher dimension   by using a recent  result   by  Goto, Nishida, Ozeki on the structure of Sally modules of an $\m$-primary ideal $\q  $ satisfying the equality
$e_1(\q)= e_0(\q) - \lambda(A/\q) +1 $ (see \cite{GNO2}).

\vskip 2mm  It is clear that, by using  Sally's machine,  the above Theorem \ref{goto2} has a natural extension to higher dimensions.

  \begin{corollary} \label{goto} Let $M$ be a     Cohen-Macaulay module of dimension $r\ge 2$, $L$ a submodule of $M$ and let $\mathbb{M}=\mathbb{M}_L$ be the good  $\q$-filtration on $M$ induced by $L.$   Assume  $e_1(\mathbb{M})=e_0(\mathbb{M})-h_0(\mathbb{M})+1$ and   $\depth  \ gr_{\mathbb{M}}(M) \ge r-2,  $  then either 
 $$P_{\mathbb{M}}(z)=\frac{h_0(\mathbb{M})+h_1(\mathbb{M})z+ z^2}{(1-z)^r} $$ or 
 $$P_{\mathbb{M}}(z)=\frac{h_0(\mathbb{M})+h_1(\mathbb{M})z+3 z^2-z^3}{(1-z)^r}.$$
 
 \end{corollary}
 \vskip 4mm

\section{The case   $e_1(\mathbb{M})=e_0(\mathbb{M})-h_0(\mathbb{M})+2 $}

\vskip 3mm
In  the  case where $e_1(\mathbb{M})=e_0(\mathbb{M})-h_0(\mathbb{M})+2  $ we present  only  partial results.
The problem is open  if  $ M_2\cap JM\neq  JM_1  $ and $e_2 \neq 2.$
  \vskip 2mm
 
 \begin{theorem} \label{primo} Let $M$ be a     Cohen-Macaulay module of dimension $r\ge 2$, $L$ a submodule of $M$ and let $\mathbb{M}=\mathbb{M}_L$ be the  good  $\q$-filtration on $M$ induced by $L.$   Assume that  $e_1(\mathbb{M})=e_0(\mathbb{M})-h_0(\mathbb{M})+2$ and $M_2\cap JM=JM_1$ where $J$ is an ideal generated by a maximal  $\mathbb{M}$-superficial sequence  for $\q.$  Then either $ gr_{\mathbb{M}}(M)$ is Cohen-Macaulay and 
 $$P_{\mathbb{M}}(z)=\frac{h_0(\mathbb{M})+h_1(\mathbb{M})z+2z^2}{(1-z)^r},$$ or $\depth  \ gr_{\mathbb{M}}(M)= r-1$ and 
 $$P_{\mathbb{M}}(z)=\frac{h_0(\mathbb{M})+h_1(\mathbb{M})z+z^3}{(1-z)^r}.$$
 
 \end{theorem}
 \begin{proof}  Since the assumptions are preserved modulo superficial elements in $J$, we may assume  $r=2.$ By Theorem \ref{r(I)} we have $$r(\mathbb{M})\le 3;$$ if $r(\mathbb{M})\le 2$ then $v_j(\mathbb{M})=0$ for every $j\ge 2$ so that $vv_j(\mathbb{M})=0$ for the same values of $j.$ Since by assumption $vv_1(\mathbb{M})=0,$ the associated graded module $gr_{\mathbb{M}}(M) $ is Cohen-Macaulay. Thus, by Theorem \ref{e1=CM}, we get $$e_0(\mathbb{M})-h_0(\mathbb{M})+2=e_1(\mathbb{M})=v_0(\mathbb{M})+v_1(\mathbb{M}),$$ which implies $v_1(\mathbb{M})=2$, $e_2(\mathbb{M})=2$ and $e_j(\mathbb{M})=0$ for every $j\ge 3.$ These values of the $e_i$'s give the required Hilbert series.
 
 Let $r(\mathbb{M})= 3;$ we have $$M_2\cap J\w{M}_1\subseteq  M_2\cap JM=JM_1$$ so that $1\in S_J.$ Further, $\text{depth} \ gr_{\w{\mathbb{M}}}(M)\ge 1 =r-1,$ hence   $$e_1(\mathbb{M})=e_1(\w{\mathbb{M}})=\sum_{i\ge 0}v_i(\w{M}).$$ By Theorem \ref{1.3} we get \begin{equation*}\begin{split} 3=r(\mathbb{M})& \le  \sum_{i\ge 0}v_i(\w{M})+2-v_0(\mathbb{M})-v_1(\mathbb{M})\\ &= e_1(\mathbb{M})+2-v_0(\mathbb{M})-v_1(\mathbb{M})\\
 &= e_0(\mathbb{M})-h_0(\mathbb{M})+2+2-e_0(\mathbb{M})+h_0(\mathbb{M})-v_1(\mathbb{M})\\ &=4-v_1(\mathbb{M})\end{split}\end{equation*}
 which implies $v_1(\mathbb{M})\le 1.$ Since $r(\mathbb{M})= 3,$ we cannot have 
 $v_1(\mathbb{M})=0,$ hence $v_1(\mathbb{M})=1$ and, by Corollary \ref{Sallymodule}, we get 
 $$P_{\mathbb{M}}(z)=\frac{ h_0(\mathbb{M})+ h_1(\mathbb{M})z+z^s}{(1-z)^2}.$$
 This gives $$h_1(\mathbb{M})+s=e_1(\mathbb{M})=e_0(\mathbb{M})-h_0(\mathbb{M})+2=h_0(\mathbb{M})+h_1(\mathbb{M})+1-h_0(\mathbb{M})+2,$$ which implies $s=3$ and $\depth  \ gr_{\mathbb{M}}(M)= 1.$
 \end{proof}
\vskip 2mm Let us remark that, as already noted in \cite{ERV} for  the case of the $\m$-adic filtration, both the Hilbert functions given in the theorem are realizable.

 In the above theorem we can get rid of the assumption $M_2\cap JM=JM_1$, but then we need another strong requirement, the condition $e_2(\mathbb{M})=2.$ The following theorem has been proved in the case of $\q$-adic filtration by  Sally (see \cite{S5}).  This is  a very deep result which gives a new class of ideals for which there is equality of the Hilbert function $H_I(n)$ and the Hilbert polynomial $p_I(n)$ at $n=1.$ 
 
 \begin{theorem} \label{nor+2} Let $M$ be a     Cohen-Macaulay module of dimension $r\ge 2$, $L$ a submodule of $M$ and let $\mathbb{M}=\mathbb{M}_L$ be the $\q$-good filtration on $M$ induced by $L.$     Assume that  $e_1(\mathbb{M})=e_0(\mathbb{M})-h_0(\mathbb{M})+2$  and $e_2(\mathbb{M})=2.$ Then $ \depth  \ gr_{\mathbb{M}}(M)\ge r-1$ and $$P_{\mathbb{M}}(z)=\frac{h_0(\mathbb{M})+h_1(\mathbb{M})z+2z^2}{(1-z)^r}.$$
 \end{theorem}
 \begin{proof}   By  Sally's  machine, we may   reduce the problem to  the case $ {r=2}.$
As usual, $J$  is an ideal generated by a maximal sequence of $\mathbb{M}$-superficial  elements 
for $\q.$  First of all, we show that $M_1=\widetilde{M}_1.$ Let $$\widetilde{L}:=\widetilde{M}_1=M_{k+1}:\q^k=\q^kL:\q^k$$ and $\mathbb{N}$ the 
 $\q$-good filtration on $M$ induced by $\widetilde{L},$ so that 
 $$\widetilde{\mathbb{M}}=\{M\supseteq \widetilde{M}_1\supseteq  \widetilde{M}_2 \supseteq \cdots \supseteq\widetilde{M}_{j+1}\supseteq \cdots \}$$ where $\widetilde{M}_{j+1}=\bigcup_k(\q^{k+j}L:\q^k)$ and
 $$\mathbb{N}:=\{M\supseteq \widetilde{L}\supseteq  \q\widetilde{L}\supseteq \cdots \supseteq \q^j\widetilde{L}\supseteq \cdots \}.$$ If $a\in \q^j$ and $m\in \widetilde{L}$, then we have $$a m \q^k\subseteq a\q^kL\subseteq\q^{j+k}L.$$ 
 
\noindent  Hence $$M_{j+1}=\q^jL\subseteq N_{j+1}=\q^j\widetilde{L}\subseteq \q^{j+k}L:\q^k\subseteq\widetilde{M}_{j+1}.$$ This implies that $$e_i(\mathbb{M})=e_i(\mathbb{N})=e_i(\widetilde{\mathbb{M}})$$ for every $i=0,\cdots,r.$ 
 
\noindent  We apply (\ref{north}) to the filtration $\mathbb{N}$ and we get $$e_1(\mathbb{N})=e_1(\mathbb{M})\ge e_0(\mathbb{N})-h_0(\mathbb{N})=e_0(\mathbb{M})-h_0(\mathbb{N}).$$ By  Theorem \ref{nor} and Theorem \ref{sal}  and since   $e_2(\mathbb{M})=e_2(\mathbb{N})=2,$ we get $$e_1(\mathbb{M})\ge  e_0(\mathbb{M})-h_0(\mathbb{N})+2.$$ So we have $$\la(M/\widetilde{M}_1)=h_0(\mathbb{N})\ge e_0(\mathbb{M})-e_1(\mathbb{M})+2=h_0(\mathbb{M})=\la(M/M_1)$$ which implies $ {M_1=\widetilde{M}_1}.$
 
\noindent      Since $\text{depth} \ gr_{\widetilde{\mathbb{M}}}(M)\ge 1=r-1$,  we have $$e_2(\mathbb{M})=e_2(\widetilde{\mathbb{M}})=\sum_{j\ge 1}jv_j(\widetilde{\mathbb{M}})=2.$$
  If $v_1(\widetilde{\mathbb{M}})=0,$ then $M_2\subseteq \widetilde{M}_2=J\widetilde{M}_1=JM_1\subseteq  M_2$ and therefore, by Corollary \ref{ind},  $e_2(\mathbb{M})=0,$  a contradiction.

\noindent  Thus $ v_1(\widetilde{\mathbb{M}})=2,$ and $ v_j(\widetilde{\mathbb{M}})=0$ for every $j\ge 2.$   From this we get \begin{equation}\label{eq} 2=\la(\widetilde{M}_2/JM_1)=\la(\widetilde{M}_2/M_2)+\la(M_2/JM_1).\end{equation} We cannot have $\la(M_2/JM_1)=v_1(\mathbb{M})\le 1$, otherwise $e_2(\mathbb{M})\not =2$ by  corollaries \ref{ind} and  \ref{Sallymodule}
 
\noindent  Hence, we must have $v_1(\mathbb{M})=2,$  so that, by (\ref{eq}),  $\widetilde{M}_2=M_2.$  Since $ v_j(\widetilde{\mathbb{M}})=0$ for every $j\ge 2,$ we immediately get $\widetilde{M}_j=M_j$ for every $j\ge 0$ and $\text{depth} \ gr_{\mathbb{M}}(M)\ge 1.$ By Theorem \ref{e1=CM}, this implies $e_2(\mathbb{M})=\sum_{j\ge 1}jv_j(\mathbb{M})=2,$ hence $v_j(\mathbb{M})=0$ for every $j\ge 2$ and   $e_j(\mathbb{M})=0$ for every $j\ge 3$;  this  gives the required Hilbert series.
 \end{proof}
 \vskip 2mm The following example shows that the assumptions $e_1(\mathbb{M})=e_0(\mathbb{M})-h_0(\mathbb{M})+2$  and $e_2(\mathbb{M})=2$ in Theorem  \ref{nor+2}   
 do not imply that $gr_{\mathbb{M}}(M)$ is Cohen-Macaulay. 
 
\begin{example}{\rm  Let $A=k[[x,y]],$ $\q=(x^6,x^5y,x^4y^9,x^3y^{15},x^2y^{16},xy^{22},y^{24})$ and $\mathbb{M}$ the $\q$-adic filtration. Then we have $$P_{\mathbb{M}}(z)=\frac{87+37z+2z^2}{(1-z)^2}$$ so that $e_0(\mathbb{M})=126,$ $e_1(\mathbb{M})=41,$ $e_2(\mathbb{M})=2,$ and $h_0(\mathbb{M})=87.$ We have $41=126-87+2$ but the associated graded ring is not Cohen-Macaulay.}
\end{example}
 
 \vskip 2mm  
 
 The following example shows that,    in  Theorem \ref{nor+2}, the assumption $e_2=2$ is essential. 
 
 \begin{example} {\rm Let $A=k[[x,y]],$ $\q=(x^6,x^5y^2,x^4y^6,x^3y^8,x^2y^9,xy^{11},y^{13})$ and $\mathbb{M}$ the $\q$-adic filtration. Then we have $$P_{\mathbb{M}}(z)=\frac{49+25z+3z^2+z^3-z^4}{(1-z)^2}$$ so that $e_0(\mathbb{M})=77,$ $e_1(\mathbb{M})=30$  and $h_0(\mathbb{M})=49.$ We have $30=77-49+2.$ 
 Further $X^4Y^5\notin \q,$ $X^4Y^5 \in \q^3:\q^2,$ so that $M_1\not= \w{M}_1.$ This implies that the associated graded ring has depth zero. Of course,  $e_2(\mathbb{M})=0\not=2.$
 }
 \end{example}

 \chapter{Applications to the Fiber Cone}

 Let $(A, \m) $ be a commutative local ring and let $\q$ be an ideal of $A.$ As usual $\mathbb{M} $ denotes a good $\q$-filtration of a module $M $ of dimension $r $  and $gr_{\q}(A)=\oplus_{n \ge 0}\q^n/\q^{n+1}$  the associated graded ring to $\q.$  Given   an ideal $I  $ containing $\q, $ we define  the  graded   module on $gr_{\q}(A)$  $$F_I (\mathbb{M}):= \oplus_{n\ge 0} M_n/I M_n.  $$
  
\vskip 2mm
\noindent  $F_I (\mathbb{M}) $ is  called the Fiber cone of $\mathbb{M}$ with respect to $I.$ If $\mathbb{M}$ is the $\q$-adic filtration on $A $ and $I =\m, $ then we write $F_\m(\q)=\oplus_{n\ge 0} \q^n/\m \q^n $ which is the classical  definition of the  Fiber cone of $\q. $ It coincides with $ gr_{\m}(A) $ when  $ \q=\m.$ \par This graded object encodes a lot of  information about  $\q.$ For instance, its dimension gives  the minimal number of generators of any minimal reduction of $\q,$ that is the analytic spread of $\q $ and its Hilbert function  determines  the minimal number of generators of the powers of $\q.$  We remark  that  if 
  $F_\m(\q) $ is Cohen-Macaulay, then the reduction number of $\q$ can be read off directly from the Hilbert series of  $F_\m(\q)  $  (see Proposition 1.85 \cite{VW3}).

Usually the arithmetical properties of the Fiber cone  and those of the associated graded ring  were   studied via apparently   different  approaches. The literature concerning the associated graded rings is  much   richer, but  new and specific   techniques were  necessary in order to study the Fiber cone.   In spite of the fact  that    $F_I (\mathbb{M}) $  is not the  graded module associated  to a filtration,   the aim of this section is to show that  it is possible to deduce information about  $F_I(\mathbb{M})$ as a consequence of the   theory on    filtrations.   In particular we will obtain recent results on the fiber cone of an ideal as an easy consequence  of classical results on the  associated graded rings of certain special filtrations.  

In this chapter we prove  some results  
which were  recently  obtained   using  different devices, often very technical ones.   Of course,  we are not going to give a complete picture of the   literature on the Fiber Cone, but   we have selected  some results which illustrate well  the use  of this approach.  Our hope is that  this method will  be useful in the future to prove new results on this topic. 

\vskip 2mm

     \section{Depth of the Fiber Cone}
     \vskip 2mm
  
\noindent Cortadellas and Zarzuela     proved in   \cite{CZ} and  \cite{Co} the existence of an exact sequence  of  the homology  of modified  Koszul complexes  which relates   $F_I (\mathbb{M}) $  with the associated graded modules to the filtrations $ \mathbb{M}  $ and $ \mathbb{M}^I.  $ We recall that the filtration  $\mathbb{M}^I $ on $M$ is defined as follows:
\begin{equation} \label{mk} \mathbb{M}^I : \ \ M\supseteq IM\supseteq IM_1\supseteq\dots \supseteq I M_n\dots\supseteq \dots \end{equation} It is clear that $\mathbb{M}^I $ is a good $\q$-filtration on $M, $ thus  $e_0(\mathbb{M})=e_0(\mathbb{M}^I ).$

The  idea by   Cortadellas and   Zarzuela has not been exploited so deeper. Starting from their work, we are going to present several applications.  First we prove a  result  which relates    the   depth  of  $F_I(\mathbb{M}), $ with  the depths of $  \ gr_{\mathbb{M}}(M) $ and $gr_{  \mathbb{M}^I}(M).$  Since the involved objects  are graded modules on $gr_\q(A),$ the depths are always computed with respect to $Q=\oplus_{n>0} \q^n/\q^{n+1}.$
\vskip 2mm

\vskip 2mm
\begin{proposition} \label{depth} Let $\mathbb{M} $ be  a good $\q$-filtration on  a module $M$ and let $I$ be an ideal containing $\q$ such that $M_{n+1} \subseteq I M_n$ for every $n\ge 1.$ We have

 \vskip 2mm
1. $\depth \ F_{I}(\mathbb{M}) \ge \text{min}\{ \depth  \ gr_{\mathbb{M}}(M) +1, \depth \ gr_{  \mathbb{M}^I }(M)\} $
\vskip 2mm
2. $ \depth  \ gr_{  \mathbb{M}^I}(M)  \ge  \text{min}\{ \depth  \ gr_{\mathbb{M}}(M), \  \depth  \ F_{I}(\mathbb{M})\} \  $ 
\vskip 2mm
3. $\depth  \ gr_{\mathbb{M}}(M) \ge  \text{min}\{ \depth  \ gr_{  \mathbb{M}^I}(M), \ \depth  \ F_{I}(\mathbb{M})\}- 1.$  

\end{proposition} 
\begin{proof}  We have the following homogeneous exact sequences of
$gr_{\q}(A)$-graded modules:
$$ 0\to N \to  \ gr_{\mathbb{M}}(M) \longrightarrow \ F_{I}(\mathbb{M}) \to 0$$ 
$$ 0\to \ F_{I}(\mathbb{M})  \to  \ gr_{  \mathbb{M}^I}(M) \longrightarrow N(-1) \to 0$$ 
where $N = \oplus_{n \ge 0} I M_n/M_{n+1}.$   

 It is enough to remark that we have the following exact sequences of the corresponding homogeneous parts of degree $n:$
$$0\to I M_n/M_{n+1} \to  M_n/M_{n+1}  \to  M_n/I  M_n \to 0$$ 
$$0\to \ M_n/I M_n   \to  I M_{n-1}/I M_n \longrightarrow   I M_{n-1}/M_{n} \to 0.$$
The inequality  between the depths  follows  from  standard facts (see for example  depth's formula  in  \cite{BH}).
\end{proof}
 Several examples show that $F_\m(\q) $ can be Cohen-Macaulay even if gr$_\q(A)$ is not Cohen-Macaulay and conversely. The above proposition clarifies the intermediate role of  the   graded module  associated to the filtration ${\mathbb{M}_{\m}}= \{\m \q^n \}$.  
 
 It will be useful   to remember that,  by Remark \ref{conti},  it is possible to find a superficial sequence $a_1, \dots, a_r $ in $\q$ which is both $\mathbb{M}$-superficial and $ \mathbb{M}^I$-superficial for $\q.$

 As  a consequence of the above proposition,  we immediately  obtain    already known  results: for example Theorem 1, Theorem 2 in  \cite{SH},   Theorem 2 in \cite{S2}, Theorem 3.4. in  \cite{HS},  Proposition 4.1, Corollary 4.3., Proposition 4.4. in \cite{Co}.

We present  here a proof of  Theorem 1 in  \cite{SH}, to show an explicit application  of our approach. We prove  the original statement, but it could easily be extended  to  modules.

  \begin{theorem} \label{Shah} Let $\q $ be an ideal of a local  ring $(A, \m) $ and let  $J$ be   an ideal generated by a superficial regular sequence  for $\q$  such that $\q^2=J \q.$ Then $F_{\m}(\q)$ is Cohen-Macaulay.
  \end{theorem}
  \begin{proof} By using the assumption, we  get   that $ \q^{n+1} \cap J=J\q^n$ and $\m\q^{n+1} \cap J = J \m \q^n$ for every integer $n.$
  By the Valabrega-Valla criterion it follows  that the    filtrations $\{\q^n\} $ and $\{\m\q^n\}$ on $A$ have associated graded rings  of 
    depth at least $ \mu(J)= $ dim $F_{\m}(\q).$ The result follows now by Proposition \ref{depth} 1.
  \end{proof}
  \vskip 3mm
  
\noindent   We remark that in the above  case we can easily write the Hilbert series of  the standard graded $k$-algebra $F_{\m}(\q),  $   that is $P_{F_\m(\q)}(z)= \sum_{i\ge 0} \dim_k (\q^i/\m \q^{i})  z^i.$ In fact 
   $$P_{F_\m(\q)}(z)=\frac{ 1}{(1-z)^r}P_{F_\m(\q)/JF_\m(\q)}(z)=\frac{1}{(1-z)^r}\sum_{i\ge 0}^s \ \dim_k(\q^i/J \q^{i-1} + \m \q^i)z^i  $$
   Since $\q^2=J\q $ and $\dim_k (\q/J+\m \q)= \mu(\q) -r, $ one has 
     $$P_{F_\m(\q)}(z)=\frac{ 1+ (\mu(\q)-r )z}{(1-z)^r}.$$
   \vskip 5mm
   \section{The Hilbert function of the Fiber Cone}
   \vskip 2mm
 If  $ \lambda(M/I  M) $ is finite, then $M_n/IM_n$ has finite length   and    we may define for every integer $n$ the numerical function
 $$H_{F_{I}(\mathbb{M})}(n):= \lambda(M_n/IM_n)$$
which is  the Hilbert function of $F_{I}(\mathbb{M}).$ We denote by $P_{F_{I} }(z) $ the corresponding Hilbert series, that is $ \sum_{i\ge 0} H_{F_{I}(\mathbb{M})}(i) z^i.  $   \par   \vskip 2mm  From now on we shall assume  that $\la (M/\q M) $ is finite. In this case    $\dim F_{I}(\mathbb{M})=r= \dim M.$  We recall that $ H_{F_{I}(\mathbb{M})}(n)$ is a polynomial function and the corresponding polynomial $ p_{F_{I} }(X) $ has  degree $r-1.$      It is the  Hilbert polynomial of $F_{I}(\mathbb{M}) $ and, as usual, we can write   $$p_{F_{I}(\mathbb{M})}(X)=\sum_{i=0}^{r-1}(-1)^i f_i(\mathbb{M})\binom{X+r-i-1}{r-i-1}.$$ 
 The coefficients $f_i(\mathbb{M})$ are integers and they are  called  the Hilbert coefficients of  $F_{I}(\mathbb{M}).$ In particular $f_0(\mathbb{M})$ is the multiplicity of the fiber cone of $\mathbb{M}.$ 
 \vskip 2mm
We can relate the Hilbert coefficients of  $F_{I}(\mathbb{M}) $ to  those  of  the filtrations $\  \mathbb{M} $ and $   \mathbb{M}^I  $ in a natural way.
 We remark that, for every $n \ge 0,$ we have 
 \begin{equation} \label{HS}  \lambda(M/M_n) \ + \  \lambda(M_n/I M_n)\ = \  \lambda(M/I M_n)   \end{equation}
 Hence
 \begin{equation} \label{P} p^1_{ \mathbb{M} }(X-1) +   p_{F_{I} }(X) =\ p^1_{   \mathbb{M}^I }(X) \end{equation}
 and
 \begin{equation} \label{H} zP^1_{ \mathbb{M} }(z) +   P_{F_{I} }(z) =\ P^1_{   \mathbb{M}^I}(z). \end{equation}
 
\noindent Since 
$$ p^1_{\mathbb{M}}(X-1)=\sum_{i=0}^{r}(-1)^i e_i(\mathbb{M})\binom{X+r-i-1}{r-i} $$  and  $$ p^1_{   \mathbb{M}^I }(X)=\sum_{i=0}^{r}(-1)^i e_i(   \mathbb{M}^I)\binom{X+r-i}{r-i  }   $$    
from (\ref{P}),  it is possible to prove that 
 \vskip 4mm
 
\begin{eqnarray}  \label{ffi}   {\fbox{  $e_0(\mathbb{M})=e_0(  \mathbb{M}^I) \  \     \text{and} \ \ \ \  
f_{i-1}(\mathbb{M })= e_i(\mathbb{M}) + e_{i-1}(\mathbb{M}) - e_i(\mathbb{M}^I)  $ }} \end{eqnarray}

\vskip 4mm \noindent
for every $ i=1,\dots, r.$ 

\vskip 2mm
\noindent Hence the theory developed    in the previous  sections on the Hilbert coefficients of the   graded module associated to a good filtration on $M$ can be applied to  $ e_i(\mathbb{M}) $ and $ e_i(  \mathbb{M}^I) $ in order to get  information, via   (\ref{ffi}),  on the coefficients of the Fiber cone of  $\mathbb{M}. $   

\section{A version of  Sally's conjecture for the Fiber Cone}

 \vskip 2mm
 We present a short proof of     the  main result of \cite{JV1}, Theorem 4.4.,  which is   the analog   of Sally's conjecture in the case of the fiber cone. 
 
  \vskip 2mm
 \begin{theorem} {\label{G-almost}} Let  $\mathbb{M}$ be the  $\q$-adic filtration on a Cohen-Macaulay module $M $  of dimension $r$ and let  $I $ be  an ideal containing $\q. $  Assume $\mathbb{M}  $   has  almost Goto minimal multiplicity  with respect to $I$ and   $\depth \ gr_{\mathbb{M}}(M) \ge r-2. $ 
  
  \noindent  Then   $\depth F_I(\mathbb{M}) \ge r-1.$
 \end{theorem}
\begin{proof} We recall that $\mathbb{M}  $   has  almost Goto minimal multiplicity  with respect to $I $ if and only if $ \ \mathbb{M}^I \ $ has \ almost  minimal multiplicity if and only if  one has $ \ \lambda(I M_1/J I M)~=~1  $ for every ideal $J$ generated by a maximal superficial sequence for $\q, $ equivalently $\lambda(I\q M/J IM)=1.  $   Hence by Corollary {\ref{Sallymodule}}, we get $ \depth  gr_{\mathbb{M^I}}(M) \ge r-1  $ and the result follows now  by Proposition \ref{depth}.

  \end{proof}
 \vskip 3mm
We remark that, under the assumptions  of Theorem \ref{G-almost}, we are able to write the Hilbert series of $F_I(\mathbb{M}).$  In fact, by using (\ref{H}) and Theorem \ref{Sallymodule}, we get
$$P_{F_I(\mathbb{M})}(z)= \frac{ \la(M/IM) + [e_0((\mathbb{M})- \la(M/IM)-1]z  +z^s- z h_{\mathbb{M}}(z)}{(1-z)^{r+1}}$$ for some integer $s\ge 2.$
\vskip 3mm
The following  example shows that in Theorem \ref{G-almost},   the assumption  $\depth \ gr_{\mathbb{M}}(M) \ge r-2  $   is necessary.

\begin{example} {\rm{ Let $A= k [[x,y,z]] $ and $ \q= (y^2-x^2, z^2-y^2, xy, yz, zx). $  The ideal $J= (y^2-x^2, z^2-y^2, xy) $ is generated by a maximal superficial sequence for $\q$ and $\la(\m \q/ \m J) =1. $ Therefore if we consider $\mathbb{M}=\{\q^n\} $ the $\q$-adic filtration on $A$ and $I=\m$ the maximal ideal of $A, $ the filtration $\mathbb{M} $ has Goto almost minimal multiplicity since $\mathbb{M}^I=\{ \m \q^n \} $ has almost minimal multiplicity.  Since $ x^2 \in \q^2 : \q, $ but $x^2 \not \in \q, $ it follows that depth  $gr_{\q}(A) = 0.$ In this case $ \depth  F_{\q}(A) =1 <r-1 $ (cfr. Example 4.5 in \cite{JV1}). }} \end{example}

It is   possible to see the above theorem  as consequence of the following  more general result.

 \begin{theorem}{\label{SF}}  Let  $\mathbb{M}$ be the  $\q$-adic filtration on a Cohen-Macaulay module $M $  of dimension $r$ and let  $I$ be  an ideal containing $\q. $      Assume 
 \vskip 2mm
 1. $\ \depth\ gr_{\mathbb{M}}(M) \ge r-2 $ 
 \vskip 2mm 
 2. $\lambda(I \q^2 M/J I \q M) \le 1$ and $I \q M \cap JM= IJM $ for some ideal $J$ generated by a maximal superficial sequence for   $\mathbb{M}^I. $ 
 \vskip 2mm \noindent Then   $ \depth F_{I}({\mathbb{M}}) \ge r-1.$
 \end{theorem}
 \begin{proof} Since $I  \q M \cap JM= IJM $  and we assume $\lambda (I \q^2M/J I \q M) \le 1, $ by    Theorem \ref{p}  applied to $\mathbb{M}^I, $ we get  $\depth   gr_{\mathbb{M}^{I}}(M) \ge r-1. $  The result follows now  by Proposition \ref{depth}.
   \end{proof}
 
 We remark that    in the classical case of the $\q$-adic filtration of $A $ and $I=\m, $  the assumption    $\m \q \cap J=\m J $ is always satisfied.
 
 \vskip 3mm
 
In Theorem {\ref{G-almost}} we   discussed $\depth  F_{I}({\mathbb{M}}) $ when $\mathbb{M}^I $ has  almost minimal multiplicity. A natural question arises about  the depth of  $F_{I}({\mathbb{M}}) $ when  $\mathbb{M}$  has almost minimal multiplicity,  that is $$e_0(\mathbb{M})=
(1-r)\lambda(M/M_1)+\lambda(M_1/M_2) +1 $$ or equivalently $\lambda(M_2/J M_1)=1 $ for every ideal $J$ generated by a maximal $\mathbb{M}$-superficial sequence.
By   Corollary {\ref{1.6}}, this assumption guarantees  that  $\depth   gr_{\mathbb{M}}(M) \ge r-1.$ 
Examples show that
$F_I(\mathbb{M}) $ is not necessarily Cohen-Macaulay and it is natural to ask whether $\depth  F_{I}({\mathbb{M}})  \ge r-1. $

  We have the analogous result of  Theorem {\ref{G-almost}}.

 \vskip 2mm
 \begin{theorem} {\label{GG-almost}}  Let  $\mathbb{M}$ be the  $\q$-adic filtration on a Cohen-Macaulay module $M $  of dimension $r$ and let  $I$ be  an ideal containing $\q. $    Assume $\mathbb{M}$   has  almost minimal multiplicity  and  $\depth gr_{\mathbb{M}^I}(M) \ge r-1. $ 
 
\noindent  Then $ \depth  F_I(\mathbb{M}) \ge r-1.$
 \end{theorem}

In a recent paper, A.V. Jayanthan, T. Puthenpurakal  and J. Verma (Theorem 3.4. \cite{JPV}) proved a criterion for the Cohen-Macaulayness of  $F_\m(\q) $ when  $\q$ has almost minimal multiplicity giving an answer to a question raised by G.Valla.
We give here a  proof by using our approach.

\begin{theorem} Let $\q$ be an $\m$-primary ideal of a local      Cohen-Macaulay ring $(A, \m)$  of dimension $r. $        Assume $\q$   has  almost minimal multiplicity and let $J$ be an ideal generated by a maximal   superficial sequence for $\q. $
   Then  the following conditions are equivalent
   \vskip 2mm  \noindent
   1.  $\m \q^2 =J \m \q   $ \vskip 2mm  \noindent
   2.  $F_{\m}( \q)  $ is Cohen-Macaulay \vskip 2mm \noindent 
   3.  $P_{F_\m(\q)}(z)=\frac{1+\lambda(\q/J +\q\m)z+z^2+\dots +z^s}{(1-z)^r}  $ for some integer $s\ge2.$
     \end{theorem}
 \begin{proof} As usual, denote by   $ \mathbb{M}^{\m}=\{\m\q^n\} $ the  filtration on $A.$    Since  $\q$   has  almost minimal multiplicity, then $\lambda(\q^2/J\q) =1.$ Hence, by Corollary \ref{1.6}, $ \depth gr_{\q}(A) \ge r-1. $  
Now,  if $\m \q^2 =J \m \q,   $   by using the Valabrega-Valla criterion,  we have that  gr$_{\mathbb{M}^{\m}}(A) $ is Cohen-Macaulay and hence  $F_{\m}(\q)  $ is Cohen-Macaulay by Proposition \ref{depth} proving 1. implies 2. 
  
  \noindent Assume now that  $F_{\m}(\q)  $  is Cohen-Macaulay. If $J=(a_1,\dots,a_r), $  we recall that 
  the corresponding   classes in $\q/\m\q $  form    a system of parameters for $F_\m(\q)  $  and hence a regular sequence on $F_\m(\q).  $  It follows that 
 $$P_{F_\m(\q)}(z)=\frac{ 1}{(1-z)^r}P_{F_\m(\q)/JF_\m(\q)}(z)=\frac{1}{(1-z)^r}\sum_{i\ge 0 } \ \lambda(\q^i/J \q^{i-1} + \m \q^i)z^i.  $$   Since $\lambda(\q^2/J\q) =1, $  then $\lambda(\q^{i+1}/J\q^i) \le 1 $ for every $i\ge 1.$ Let $s\ge 2$   the least integer such that $\q^{s+1}=J \q^s.$ Since $\lambda(\q^i/J\q^{i-1})=1  $ for $i=2, \dots, s $  and  hence $ \m \q^i \subseteq J \q^{i-1} $ for every $i \ge 2, $ we get
  $$P_{F_\m(\q)}(z)=\frac{1+\lambda(\q/J +\q\m)z+z^2+\dots +z^s}{(1-z)^r} .$$  It follows that 2. implies 3.  Actually 2. is equivalent to 3. In fact  if $P_{F_\m(\q)}(z)=\frac{ 1}{(1-z)^r}P_{F_\m(\q)/JF_\m(\q)}(z), $ then  $F_{\m}( \q)  $ is Cohen-Macaulay.
  
  \noindent 
 We prove now  3. implies 1., that is $v_2( \mathbb{M}^{\m} )=\la(\m \q^2/J \m \q)=0.$  Since  $F_{\m}( \q)  $ is Cohen-Macaulay and   depth gr$_{\q}(A) \ge r-1,$   by Proposition \ref{depth}, depth gr$_{\mathbb{M}^{\m}} (A) \ge r-1. $
Hence, by Theorem \ref{e1=CM},  $e_1({\mathbb{M}^{\m}})= \sum_{i \ge 0} v_i({\mathbb{M}^{\m}}).$
  Then we have to  prove $e_1({\mathbb{M}^{\m}})= v_0({\mathbb{M}^{\m}})+v_1({\mathbb{M}^{\m}})=e_0({\mathbb{M}^{\m}})-1 +\lambda(\m\q/J\m).$

   Now, by (\ref{ffi}), we know that $e_1(\mathbb{M}^{\m})=e_0(\mathbb{M}^{\m}) + e_1(\q) -f_0(\q). $ Since $e_1(\q)=\sum_{i\ge s} v_i(\q)= \lambda(\q/J) +s-1 $ and $f_0= 1+\lambda(\q/J +\q\m) +s-1,$ we have  $ e_1(\mathbb{M}^{\m})= e_0(\mathbb{M}^{\m}) + \lambda(\q/J) +s-1-(1+\lambda(\q/J +\q\m) +s-1)=e_0(\mathbb{M}^{\m})-1 +\lambda(\m\q/J\m), $ as required.
  \end{proof}
\vskip 2mm 
 
  \vskip 1cm

\section{The Hilbert coefficients of the Fiber Cone}

  The formulas in (\ref{ffi}) give information on the Hilbert coefficients of the Fiber Cone by means of the theory of   Hilbert functions  of filtered modules. First we get a short proof of a recent result by A. Corso (see  \cite{C}).

\begin{theorem}{\label{f0}} Let $\mathbb{M} $ be  a good $\q$-filtration on a module $M$ and let $I$ be an ideal containing $\q$ such that $M_{n+1} \subseteq I M_n.$ Let $J$ be the ideal generated by a maximal   $  \mathbb{M}^I $-superficial sequence for $\q  $ and denote by $ \mathbb{N} $ the corresponding filtration $\{J^nM\}.$   Then  $$f_0(\mathbb{M }) \le  \text{min} \{ e_1(\mathbb{M }) -  e_0(\mathbb{M }) -e_1(\mathbb{N })+ \lambda (M/I M) +
\lambda(M/I M_1 +JM),$$ $$  e_1(\mathbb{M })   - e_1(\mathbb{N }) + \lambda(M/I M) \}.$$
 
\end{theorem}
\begin{proof} Since $f_0(\mathbb{M })=  e_0(\mathbb{M })+ e_1(\mathbb{M })- e_1(\mathbb{M^{I} }) $ by (\ref{ffi}), it is enough to apply  Theorem \ref{MeN} to  $e_1(  \mathbb{M}^I) $ for $s=1,2.$

\end{proof}

If we apply   the above result  to  the case   $ \mathbb{M}=\{\q^n\} $ with $I=\m, $  we easily obtain the following  bound on the multiplicity of the  fiber cone $F_\m (\q) $ (see  \cite{C}, Theorem 3.4.).

\begin{corollary} {\label{c}}Let $\q$ be  an $\m$-primary ideal of a local ring $(A,\m) $ of dimension $r.$ Let $J$ be the ideal generated by a maximal superficial sequence for $\q,$ then
$$f_0(\q) \le \text{min} \{e_1(\q) - e_0(\q) -e_1(J) + \lambda(A/\q) + \mu(\q) -r +1 ,  e_1(\q)   -e_1(J)   +1\}. $$ 
\end{corollary}

If $A$ is Cohen-Macaulay, then  $e_1(J)=0$ because $J$ is generated by a regular sequence and we are able to characterize  the extremal cases.   The following result   generalizes  Proposition 2.2. and Theorem 2.5.   in  \cite{CPV}.

\begin{corollary} Let $\q$ be an $\m$-primary ideal of a local Cohen-Macaulay ring $(A,\m) $ of dimension $r.$
 Then 
$$f_0(\q) \le  e_1(\q) - e_0(\q)   + \lambda(A/\q) + \mu(\q) -r +1 \le  e_1(\q)     +1 . $$ In particular
\vskip 2mm
1. If $f_0(\q)=e_1(\q)     +1, $ then   $\m \q=
\m J $ for every maximal superficial sequence $J$ for $\q.$ If, in addition, $\lambda (\q^2 \cap J/J\q) \le 1$ for some $J,$ then $ \depth  gr_{\q}(A) \ge r-1 $ and $F_{\m}(\q) $ is Cohen-Macaulay.
\vskip 2mm
2. If $ f_0(\q) =  e_1(\q) - e_0(\q)   + \lambda(A/\q) + \mu(\q) -r +1,  $ then $F_{\m}(\q) $ is unmixed. 
\end{corollary}

\begin{proof} The first inequality follows by Corollary \ref{c}.  We prove now that $ e_1(\q) - e_0(\q)   + \lambda(A/\q) + \nu(\q) -r +1 \le  e_1(\q)     +1. $ If $ J$ is an ideal generated by a maximal superficial sequence for $\q,$ then 
$ e_0(\q)   - \lambda(A/\q) - \nu(\q) +r  = \lambda (\q/J)  - \lambda(\q/\q \m)   +\lambda (J/J\m) =\lambda(\q/J\m)  - \lambda(\q/\q \m)  \ge 0.$ 

\noindent In particular if $f_0(\q)=e_1(\q)     +1, $ then $\q \m = J \m $ and hence the associated graded module to the $\q$-filtration $\mathbb{M}^{\m}= \{\m \q^n\} $ is Cohen-Macaulay by the Valabrega-Valla criterion. Now $\q^2 \subseteq \m \q= \m J  \subseteq J,$ then $\lambda (\q^2 \cap J/J\q)=  \lambda (\q^2/J\q)  \le 1. $  Hence by Theorem \ref{p}, depth $gr_{\q}(A)\ge r-1$ and 1.  follows by Proposition \ref{depth}.

\noindent Now,  from the proof of Theorem \ref{f0},  $ f_0(\q) =  e_1(\q) - e_0(\q)   + \lambda(A/\q) + \nu(\q) -r +1   $  if and only if $e_1(\mathbb {M}^{\m})= 2 e_0(\mathbb {M}^{\m}) - 1 - \lambda (A/\m \q +J)= 2 e_0(\mathbb {M}^{\m}) - 2 h_0(\mathbb {M}^{\m})- h(\mathbb {M}^{\m}) $ and hence, by Theorem \ref{elv},  $gr_{ \mathbb {M}^{\m}}(A) $
is Cohen-Macaulay. Because we have a canonical  injective map from  $F_{\m}(\q) $ to $ gr_{ \mathbb {M}^{\m}}(A) $ the result follows.
\end{proof}

\noindent We remark that, in the above result,  the assumption $\lambda (\q^2 \cap J/J\q) \le 1$ is satisfied for example if 
$A$ is Gorenstein (see Proposition 2.2. in \cite{C}). 

\vskip 3mm

\section{Further numerical invariants: the $g_i  s$}
\vskip 2mm

We give  now  short proof of several recent results proved  in \cite{JV1}, \cite{JV2} and \cite{JPV}. First   we need to relate the numerical  invariants already considered  with  those introduced  by A.V. Jayanthan and J. Verma. They  write the polynomial  $ p^1_{   \mathbb{M}^I}(X) $  of degree $r $ by using  the unusual binomial basis  $\{ \binom{X+r-i-1}{r-i}\ \ : \ \ i=0, \dots,d \}. $  The integers $ g_i ({   \mathbb{M}^I})  $ are uniquely determined 
 $$  p^1_{   \mathbb{M}^I}(X)=\sum_{i=0}^{r }(-1)^i g_i({   \mathbb{M}^I}) \binom{X+r-i-1}{r-i}.  $$ 
     They have the advantage of leading to  more compact formulas than those in  (\ref{ffi}). 
 By (\ref{P}),   it is easy to check that
   \begin{equation*}  e_0(\mathbb{M})=g_0({   \mathbb{M}^I}) \ \ \ \  \text{and} \ \ \ \  
g_i({   \mathbb{M}^I})= e_{i}(\mathbb{M}) - f_{i-1}(\mathbb{M})
 \end{equation*}
for every $ i=1,\dots, r.$
Moreover from the following  equalities
 $$  p^1_{   \mathbb{M}^I}(X)=\sum_{i=0}^{r }(-1)^i g_i ({   \mathbb{M}^I}) \binom{X+r-i-1}{r-i}=$$$$= \sum_{i=0}^{r}(-1)^i e_i(   \mathbb{M}^I)\binom{X+r-i}{r-i  }   $$ 
  we obtain 
\begin{equation*}   e_0(  \mathbb{M}^I )=g_0 ({   \mathbb{M}^I})\ \ \ \  \text{and} \ \ \ \  
e_{i}(  \mathbb{M}^I) = g_{i-1}({   \mathbb{M}^I}) +g_i ({   \mathbb{M}^I})
 \end{equation*}
for every $ i=1,\dots, r.$
Then
\begin{eqnarray}  \label{g}   \fbox{$  g_i({   \mathbb{M}^I}) = \sum_{j=0}^i (-1)^{i-j} e_j( {\mathbb{M}^I})
$ } \end{eqnarray}

\noindent In particular 
 $$ g_1({   \mathbb{M}^I})= e_{1}(  \mathbb{M}^I)\ -\ e_{0}(  \mathbb{M}^I)\  \ \ \ \text{and} \  \ \ \ \ g_2({   \mathbb{M}^I})=e_{2}(  \mathbb{M}^I)- e_{1}(  \mathbb{M}^I)+ e_{0}(  \mathbb{M}^I).$$

\noindent From the  equality (\ref{g}) it is clear   that  the integers  $g_i({   \mathbb{M}^I})$  have  good  behaviour modulo  $  \mathbb{M}^I$-superficial elements for $\q $ (see Lemma 3.5., \cite{JV2}). 
 \vskip 3mm
  
With almost no further effort   we can obtain and generalize several  results in \cite{JV2}. It will be useful to  recall that we have  $ w_0(  \mathbb{M}^I)=v_0(  \mathbb{M}^I)= \la ( M/JM)- \lambda(M/I M) $ and, for $n\ge 1, $  
  $w_n(  \mathbb{M}^I)= \lambda(IM_n+JM /JM ) $ and $ v_n(  \mathbb{M}^I)= \lambda(IM_n /JI M_{n-1}) .$ 
 
  \begin{theorem} \label{g1} {\rm{(\cite{JV2}, Proposition 4.1 and Theorem 4.3)}}  Let $\mathbb{M} $ be  a good $\q$-filtration of a Cohen-Macaulay  module $M$ of dimension $r $ and let $I $ be an ideal containing $\q$ such that $M_{n+1} \subseteq I M_n. $ Let $ J$ be the ideal generated by a maximal   $  \mathbb{M}^I $-superficial sequence for $\q.  $ Then
 \vskip 2mm
$$g_1({   \mathbb{M}^I}) \ge \sum_{n \ge 1} w_n(  \mathbb{M}^I )  - \lambda(M/IM) $$ The equality holds  if and only if  $gr_{\mathbb{M }^I } (M)$ is Cohen-Macaulay.
 \vskip 3mm
 \noindent Moreover if $\depth \ gr_{\mathbb{M}}(M) \ge r-1,$ then
 \vskip 2mm
 \noindent $g_1 ({   \mathbb{M}^I})= \sum_{n \ge 1} w_n(  \mathbb{M}^I)  - \lambda(M/I M) $ if and only if  $   F_I(\mathbb{M })  $ is Cohen-Macaulay.   
  \end{theorem}
 
 \begin{proof} It is enough to recall that $ g_1({   \mathbb{M}^I})= e_{1}(  \mathbb{M}^I)  - e_{0}(  \mathbb{M}^I).$ The result  follows by  Theorem \ref{e1=wj}   applied to the  filtration $   \mathbb{M}^I.   $
The last part  is a consequence of the  previous result  and   Proposition \ref{depth}.    
\end{proof}
 
\vskip 2mm
The next result extends and completes Theorem \ref{e1=CM}  (see also \cite{JV2}).

 \begin{theorem} \label{g1<} Let $\mathbb{M} $ be  a good $\q$-filtration of a Cohen-Macaulay  module $M$ of dimension $r $ and let $I$ be an ideal containing $\q$ such that $M_{n+1} \subseteq IM_n. $ Let $ J$ be the ideal generated by a maximal   $  \mathbb{M}^I $-superficial sequence for $\q.  $ Then we have
 \vskip 2mm  \noindent 
a) $ g_1({   \mathbb{M}^I})\le  \sum_{n \ge 1} v_n(  \mathbb{M}^I)  - \lambda(M/I M) $   
\vskip 2mm \noindent 
b)  $ g_2({   \mathbb{M}^I})\le  \sum_{n \ge 2} (n-1) v_n(  \mathbb{M}^I)  + \lambda(M/I M) $   

\vskip 2mm \noindent 
If $\depth \ gr_{\mathbb{M}}(M) \ge r-1 $ the following conditions are equivalent: 
\vskip 2mm
\noindent
1.  $ \depth   F_I(\mathbb{M }) \ge r-1. $ 
 \vskip 2mm
\noindent  
2. $g_1({   \mathbb{M}^I}) = \sum_{n \ge 1} v_n(  \mathbb{M}^I) - \lambda(M/I M) $ 
\vskip 2mm
\noindent
3.  $g_i({   \mathbb{M}^I}) = \sum_{n \ge i }  {{n-1 }\choose{i-1} }v_n(  \mathbb{M}^I)+ (-1)^i \lambda(M/I M) $ for every $i \ge 1.$

\end{theorem}
 
 \begin{proof} It is enough to recall that $ g_i({   \mathbb{M}^I}) =\sum_{j=0}^i (-1)^{i-j} e_{i}(  \mathbb{M}^I) .$  Now a) follows by  Theorem \ref{e1=CM} a).    Further    b) follows   by Theorem \ref{e1=CM} b) 
 and Northcott's inequality (Theorem \ref{MeN} for $s=1$) always applied  to the  filtration $\mathbb{M}^I.  $ 
The last part  is a consequence of Theorem \ref{e1=CM} c)  and  Proposition \ref{depth}.    
\end{proof}

  \vskip 2mm
  
  \begin{remark} {\rm{ By Theorem \ref{g1},  it follows  that  $$g_1({   \mathbb{M}^I}) \ge -\lambda(M/IM)   $$  We remark that if $g_1({   \mathbb{M}^I})= -\lambda(M/I M)$, then  $   F_I(\mathbb{M }) $ does  not necessarily have maximal depth. In fact, again by Theorem \ref{g1} it follows that   $gr_{\mathbb{M }^{I}} (M)$ is Cohen-Macaulay, but nothing is known about    $gr_{\mathbb{M }}(M).  $  
The following example taken from     \cite{G} has  minimum   $g_1({   \mathbb{M}^I}), $ nevertheless $F_\m(\q) $  is not Cohen-Macaulay. 
\vskip 2mm
Consider  $A=k[[ x^4,x^3y,x^2y^2,xy^3,y^4]] $   a  subring of the formal power series ring $k[[x,y]] $  and let $\mathbb{M}$ be the $\q$-adic filtration with $\q=( x^4,x^3y, xy^3,y^4)  $ and let $I=\m.$ We have $g_1({   \mathbb{M}^I})=-1$ since $\m\q=\m J  $ where   $ J=(x^4,y^4)A .$ In this case   $F_\m(\q) $ has depth $1, $ hence it  is not Cohen-Macaulay.}}
\end{remark}

In the case of the $\q$-adic filtration on $M$ it is possible to characterize the ideals $\q$ for which  $g_1$ is minimal. The following result generalizes Proposition 6.1. \cite{JV2}).

 \vskip 2mm
\begin{proposition}  Let  $\mathbb{M}$ be the  $\q$-adic filtration on a Cohen-Macaulay module $M $  and let  $I $ be  an ideal containing $\q. $ Then $\mathbb{M}^I  $ has minimal multiplicity    
if and only if  $g_1({   \mathbb{M}^I})=- \lambda(M/I M).  $
\end{proposition}
\begin{proof} We recall that $\mathbb{M}^I $ has  minimal multiplicity if and only if $IM_1=I \q M =J I M  $ for every ideal $J$ generated by a maximal superficial sequence for $\mathbb{M}^I.$   Then  by Valabrega-Valla's criterion, gr$_{\mathbb{M }^{I}} (M)$ is Cohen-Macaulay. Now   the result  follows by Theorem \ref{g1}, 1. since   $ w_n(  \mathbb{M}^{I})=0 $   for every $n \ge 1.$  Conversely if $g_1({   \mathbb{M}^I})=- \lambda(M/I M),    $  then  in   Theorem \ref{g1}, 1. we have the equality and    $w_n(  \mathbb{M}^{I})=0 $ for every $n \ge 1. $    Hence gr$_{\mathbb{M}^{I}} (M)$ is Cohen-Macaulay and $I M_n \subseteq J \cap I M_n=J I M_{n-1} $ for every $n \ge 1.$ In particular $I M_1=I\q M =JIM,$ as required.
\end{proof}

\include{Cap66}

\chapter{Applications to the Sally module}

W.V. Vasconcelos enlarged the list of blowup algebras by introducing the {\it{ Sally module}}.     Let $(A, \m) $ be a commutative local ring and let $\q$ be an ideal of $A, $ then 
the   { Sally module} $S_J(\q) $  of $\q$ with respect to a minimal reduction $J $ is a graded  module on the Rees algebra $\R (J)= \oplus_{n\ge 0} J^n$ which is defined in terms of the exact sequence
$$0 \to  \q \R(J) \to \q \R(\q) \to S_J(\q):= \oplus_{n \ge 1} \q^{n+1}/J^n\q \to 0.$$
A motivation for its name is the work of Sally where the underlining  philosophy  is that  it is reasonable to expect to recover some properties of $\R(\q) $ (or $ gr_{\q}(A)$)  starting from the better structure of $\R(J).$

Vasconcelos  proved that if $A$ is Cohen-Macaulay,  then $\dim  S_J(\q) = \dim  A, $ provided $ S_J(\q)$ is not the trivial module.  
\vskip 3mm
We extend the definition  to modules. 
   As usual $\mathbb{M} $ denotes a good $\q$-filtration of a module $M $ of dimension $r $ and let $J$ be the ideal generated by a maximal $\mathbb{M}$-superficial sequence for $\q.$ We define   
   $$ S_J(\mathbb{M}):= \oplus_{n \ge 1} M_{n+1}/J^nM_1 $$
to be    the Sally module of $\mathbb{M}$ with respect to $J.$ 
   
   As we saw for  the fiber cone, this graded $\R (J)$-module  is closely  related to  the associated graded modules  with respect  to  different  filtrations. We consider the $J$-good filtration  induced by the submodule $M_1 $ of $M:$
   
    \begin{equation*} \mathbb{E} : \{E_0=M,  E_{n+1}=J^nM_1  \ \ \ \forall n \ge 0 \}\end{equation*}

\vskip 2mm
 The aim  of this chapter is to relate the numerical invariants  of $S_J(\mathbb{M})$ to those  of the associated graded modules of  the filtrations $\mathbb{M} $ and $\mathbb{E}. $ As in the previous chapter, we will rediscover easily properties of the Sally module by using the general theory on the associated graded modules developed  in the previous chapters. 

\section{Depth of the Sally module}
\vskip 2mm
The Sally module $S_J(\mathbb{M}) $ fits in two exact sequences of graded $\R (J)$-modules with 
$gr_{\mathbb{M}}(M)  $ and $  gr_{\mathbb{E}} (M).$

\begin{proposition} \label{depthS} Let $\mathbb{M} $ be  a good $\q$-filtration of a module $M$ and let $J$ be the ideal generated by a maximal $\mathbb{M}$-superficial sequence for $\q.$   
Then
 $\depth  gr_{\mathbb{M}}(M)  \ge \text{min}\{  \depth  \ S_J(\mathbb{M})  - 1, \depth  \ gr_{\mathbb{E}}(M)\}. $
 \end{proposition} 
  
\begin{proof} Let $N:= \oplus_{n\ge 0} M_n/J^n M_1,$ we have the following homogeneous exact sequences of
$\R(J)$-graded modules:
$$0\to gr_{\mathbb{E}} (M) \to N  \longrightarrow \ S_J(\mathbb{M})(-1)   \to 0$$ 
$$0\to \ S_J(\mathbb{M})  \to  \ N  \longrightarrow \ gr_{\mathbb{M}}(M) \to 0$$

\noindent It is enough to remark that we have the following exact sequences of the homogeneous components  of degree $n $ :
$$0\to J^{n-1}M_1/J^{n}M_1  \to  M_{n }/J^{n-1} M_{1}  \to  M_{n }/J^{n-1}M_1 \to 0$$ 
$$0\to \ M_{n+1 }/J^{n}M_1  \to  M_n/J^nM_1 \longrightarrow   M_n/M_{n+1} \to 0.$$
We remark that $M_{n }/J^{n-1}M_1= (S_J(\mathbb{M})(-1))_n.$

\noindent  The comparison between the depths  follow  from  standard facts (see for example \cite{BH}).
\end{proof}
 
\vskip 4mm

 \section{The Hilbert function of the Sally module}
 \vskip 2mm   
\noindent    From now on we shall assume $\lambda(M/\q M) $ finite.  

\vskip 2mm If $M$ is Cohen-Macaulay, then the filtration 
$\mathbb{E}  $ is well understood. In fact, since $E_2=JE_1, $  by Theorem \ref{nor} and Corollary \ref{ind},  
$gr_{\mathbb{E}}(M) $ is Cohen-Macaulay with minimal multiplicity and  hence 
  \begin{equation} \label{PCM} P_{ \mathbb{E} }(z)=  \frac{   h_0(\M) +  (e_0(\mathbb{M})- h_0(\M))z }{(1-z)^r} .\end{equation}
($ h_0(\M)= h_0(\mathbb{E}), e_0(\mathbb{M})=e_0(\mathbb{E})$). 
In particular  \begin{equation} \label{ECM} e_1(\mathbb{E}) = e_0(\M)-h_0(\M)  \ \ \ {\rm{ and}} \ \ \ 
e_i(\mathbb{E}) =0 \ {\rm{ for \ every}} \  i \ge2.  \end{equation}

\noindent Since  $ \lambda(M_{n+1}/J^nM_1) $ is finite for every $n, $ we may define 
the Hilbert function of the Sally module 
$$H_{S_J(\mathbb{M})}(n)= \lambda(M_{n+1}/J^nM_1) $$
and we denote by $e_i(S_J(\mathbb{M}))$ the corresponding Hilbert coefficients.
  
\noindent  Starting from the exact sequences of  Proposition \ref{depthS}, it is easy to get the following equality on the Hilbert series
 \begin{equation} \label{Pz} (z-1) P_{S_J(\mathbb{M})}(z) = P_{ \mathbb{M}}(z)\ - \  P_{{\mathbb{E}}}(z) 
\end{equation}

 \noindent Several   results easily  follow  from  the above equality.   
 
 \begin{proposition} \label{ME}Let $\mathbb{M} $    be a good $\q$-filtration of a  module $M  $ of dimension $r $ and let  $J$ be the ideal generated by a maximal $\mathbb{M}$-superficial sequence for $\q. $ 
 We have
 \vskip 2mm
 1. $\dim \ S_J(\mathbb{M}) =r  $ if and only if $e_{1}(\mathbb{M})\ >   e_{1}(\mathbb{E}).$ 
 \vskip 2mm
2. If  $\dim \ S_J(\mathbb{M}) =r, $ then for every $i \ge 0$ we have 
 \begin{equation} \label{eS}  e_i(S_J(\mathbb{M})) = e_{i+1}(\mathbb{M})\ -\  e_{i+1}(\mathbb{E})
 \end{equation}
 
 \end{proposition}
 
 From  (\ref{eS}) we deduce  that   the coefficients of the Sally module have a  good behavior with  $\mathbb{M}$-superficial sequence for $\q. $  We  remark that both $ \mathbb{M}$ and  $\mathbb{E} $ are $J$-good filtrations, hence by Remark \ref{conti}  it is possible to find in $J$ a sequence of elements which are superficial for both. 
\vskip 2mm
 The following result was  proved in \cite{C} in the particular case of   the $\q$-adic filtration on $A.$ We present here a direct proof in the general setting.
 
   \begin{corollary} \label{e0S}  Let $\mathbb{M} $    be a good $\q$-filtration of a  module $M  $ of dimension $r $ and let  $J$ be the ideal generated by a maximal $\mathbb{M}$-superficial sequence for $\q. $ Let
   $ \mathbb{N}=\{J^nM\}$ be the $J$-adic filtration on $M $  and assume $ \dim S_J(\mathbb{M}) =r, $ then
   $$e_0(S_J(\mathbb{M})) \le e_1(\mathbb{M})-e_1(\mathbb{N})-   e_0(\mathbb{M})+ h_0(\mathbb{M}). $$
   \end{corollary}
   \begin{proof}  It follows from  (\ref{eS}) and   Proposition \ref{E}.
   \end{proof}

   \vskip 4mm
   
  We note  that next  result rediscovers several results  known in the case of the $\q$-adic filtration on the Cohen-Macaulay ring $A$ (see \cite{VP2}, Corollary 1.2.9., Proposition 1.2.10, Proposition 1.3.3; Corollary 2.7 \cite{RVV} ).

\vskip 2mm

 \begin{proposition} \label{PS} Let $\mathbb{M} $    be a good $\q$-filtration of a Cohen-Macaulay module $M  $ of dimension $r $ and let  $J$ be the ideal generated by a maximal $\mathbb{M}$-superficial sequence for $\q, $ then 
 \vskip 2mm
 1. $\dim \ S_J(\mathbb{M})=r$ if and only if  $e_{1}(\mathbb{M})\  > e_0(\mathbb{M}) - h_0(\mathbb{M}).  $ 
 \vskip 2mm
 2. If  $ \dim S_J(\M)=r,$ then 
    $e_0(S_J(\M))=e_{1}(\M)-e_0(\M)+h_0(\M) $ 
\vskip 2mm

3.  $ \depth \  gr_{\mathbb{M}}(M)  \ge  \depth  \ S_J(\mathbb{M})  - 1$ 
 \vskip 2mm
 4. $   (z-1) P_{S_J(\mathbb{M})}(z) = P_{\mathbb{M}}(z) -   \frac{\lambda(M/M_1) +  (e_0(\mathbb{M})- (\lambda(M/M_1))z }{(1-z)^r}.$   
 \vskip 2mm
 5. $H_{S_J(\mathbb{M})}(n)  $ is not decreasing.
\end{proposition}
  \begin{proof} Assertions 1. and 2. follow  by  Proposition \ref{ME} and  (\ref{ECM}). 
Since $\  gr_{\mathbb{E}}(M) $ is Cohen-Macaulay,   then  3.  follows by Proposition \ref{depthS} because     $ \text{min}\{  \text{depth} \ S_J(\mathbb{M})   -  1, \text{depth} \ gr_{\mathbb{E}}(M)\} =
\text{depth}\ S_J(\mathbb{M})  - 1.$ 

The assertion 4. follows from  (\ref{Pz})  and (\ref{PCM}). Finally 5. follows by 4. and Theorem \ref{nCM}.
 
 \end{proof}
 
  \vskip 3mm

  In our general setting we get  the following result due, in the classical case,  to W. Vasconcelos in  \cite{VW1}, section 5.2.
\begin{corollary} \label{ps} Let $L$ be a submodule of a Cohen-Macaulay module $M  $ of dimension $r $  and  let $ \mathbb{M}= \mathbb{M}_L$    be a good $\q$-filtration induced by $L.$   Let  $J$ be an  ideal generated by a maximal $\mathbb{M}$-superficial sequence for $\q, $ then 
\vskip 2mm
1. $\dim \ S_J(\mathbb{M})=r$ provided it is not the trivial module.
\vskip 2mm
2. $e_0(S_J(\mathbb{M}))= e_{1}(\mathbb{M})- \lambda(L/JM) $ and $e_i(S_J(\mathbb{M}))= e_{i+1}(\mathbb{M}) $ for every $i >0.$
\end{corollary}
\begin{proof} By Proposition \ref{PS}, $\dim \ S_J(\mathbb{M})=r$ provided   $e_{1}(\mathbb{M})\  > e_0(\mathbb{M}) - h_0(\mathbb{M}).  $ By Theorem \ref{nor} and Corollary \ref{ind},  this means $M_2\neq JM_1.$ 
On the other hand,   $M_{n+1}=\q^nL $ for every $n \ge 0, $  hence it is easy to see  that $S_J(\mathbb{M})$ is the trivial module if and only if $M_2=JM_1 $ and 1. follows.  Now 2. is a consequence of Proposition \ref{ME} since $e_{1}(\mathbb{E})= e_0(\mathbb{E}) - h_0(\mathbb{E})=e_0(\mathbb{M}) - h_0(\mathbb{M})=\lambda(L/JM)   $ and $e_{i}(\mathbb{E})=0$ for $i\ge 2.$
\end{proof}
\vskip 2mm 

\begin{remark} {\rm{ If $\q$ is an $\m$-primary ideal of a   local Cohen-Macaulay ring $(A, \m)   $   of dimension $r, $    the value of $e_1(\q) $ has a strong   influence on   the structure of the Sally module. By Corollary \ref{ps} and Theorem \ref{nor}, if
$e_1(\q)= e_0(\q) - \lambda(A/\q),  $ then  $S_J(\q) $ is the trivial module.  The case $e_1(\q)= e_0(\q) - \lambda(A/\q) +1 $ is much more difficult.      Recently S. Goto, K. Nishida, K. Ozeki 
 proved  that, under this  assumption,  there exists a positive integer $c \leq r$ such that 
$$  S_J(\q) \simeq   (x_1, \dots, x_c) \subseteq [\R(J)/\m \R(J)] \simeq  A/\m [x_1, \dots, x_r] $$ as  graded $\R(J)$-modules.  When this is the case $c=v_1(\q)= \lambda(\q^2/J\q) $ (see \cite{GNO2}, Theorem 1.2). 
By using this surprising information and Proposition \ref{PS}, 3., we easily   obtain
$$ P_{\q}(z)= \frac{ \lambda(A/\q) + (e_0(\q) -\lambda(A/\q) -c) z + \sum_{i=2}^{c+1} (-1)^i {{c+1}\choose i} z^i}{(1-z)^r} $$

We remark that,  if $r=2, $ we obtain the Hilbert series described in   Theorem \ref{goto2}.  Very recently S. Goto and K. Ozeki announced an extension of the above result  relaxing the requirement for the Cohen-Macaulyness of $A. $ \cite{GO}.  
}}
\end{remark}
\vskip 5mm
\noindent Under the assumptions of Corollary \ref{ps},  we remark that if $J=(a_1, \dots, a_r) $ and $S_J(\mathbb{M}) $ is not the trivial module,
 then the ideal $JT=(a_1T, \dots, a_rT )$ in  the Rees algebra ${\cal R}(J)=A[JT] $ is generated by a system of parameters for $S_J(\mathbb{M}). $
 In fact $S_J(\mathbb{M})/JT S_J(\mathbb{M})= \oplus_{n\ge1} M_{n+1}/JM_n $  which is an Artinian module. 
 \par \noindent In particular $S_J(\mathbb{M}) $ is Cohen-Macaulay if and only if $a_1T, \dots, a_rT$ is a regular sequence on $S_J(\mathbb{M}). $
 \vskip 2mm
\begin{theorem} \label{vaz} {\rm{(\cite{VP2}, Theorem 2.1.6 and Corollaries 2.1.7, 2.1.8, 2.1.9)}}  Let $L$ be a submodule of a Cohen-Macaulay module $M  $ of dimension $r $  and  let $ \mathbb{M}= \mathbb{M}_L$    be a good $\q$-filtration induced by $L.$   Denote  by $J$ the  ideal generated by a maximal $\mathbb{M}$-superficial sequence for $\q, $ then 
$$ e_0(S_J(\mathbb{M})) \le \sum_{j\ge 1} v_j(\mathbb{M}).$$
The following facts are equivalent:
\vskip 2mm
1. $ e_0(S_J(\mathbb{M})) = \sum_{j\ge 1} v_j(\mathbb{M})  $ 
\vskip 2mm
2. $ e_1(\mathbb{M}) = \sum_{j\ge 0} v_j(\mathbb{M})  $
\vskip 2mm
3. $ \depth gr_{\mathbb{M}}(M) \ge r-1$
\vskip 2mm
4. $ P_{ \mathbb{M}}(z)  =   \frac{\lambda(M/M_1) + \sum_{j\ge 1} (v_{j-1}(\mathbb{M}) - v_{j}(\mathbb{M})) z^j }{(1-z)^r}$
\vskip 2mm
5.  $ P_{S_J(\mathbb{M})} (z)  =   \frac{  \sum_{j\ge 1} v_{j}(\mathbb{M})  z^j }{(1-z)^r}$
\vskip 2mm
6.  $ S_J(\mathbb{M}) $ is Cohen-Macaulay.
 
\end{theorem}
\begin{proof} By Corollary \ref{ps} we have $e_0(S_J(\mathbb{M}))= e_{1}(\mathbb{M})- \lambda(L/JM). $ Hence, by Theorem \ref{e1=CM}, we get  $e_0(S_J(\mathbb{M}))\le \sum_{j\ge 0} v_j(\mathbb{M})-\lambda(L/JM)=\sum_{j\ge 1} v_j(\mathbb{M}).$
The equality holds if and only if $ e_{1}(\mathbb{M})=\sum_{j\ge 0} v_j(\mathbb{M}).$ Hence, by  Theorem \ref{e1=CM},  the assertions 1., 2., 3., 4.   are equivalent.  By Proposition \ref{PS} (4.), 4. and 5. are equivalent.  Since 6. implies 3. by Proposition \ref{PS} (3.),  it is enough  to prove that 5. implies 6.

We may assume $S_J(\mathbb{M})$ has dimension $r.$ We recall that $S_J(\mathbb{M})$ is a $\R (J)$-module and we have $S_J(\mathbb{M})/JT  S_J(\mathbb{M})=\oplus_{n\ge1} M_{n+1}/JM_n.$  From  5. we deduce that   $ P_{S_J(\mathbb{M})} (z)=
\frac{1}{(1-z)^r} P_{S_J(\mathbb{M})/JT  S_J(\mathbb{M})}(z).$ Then $JT$ is generated by a regular sequence of length  $r=$dim$S_J(\mathbb{M})$ and hence $ S_J(\mathbb{M}) $ is Cohen-Macaulay.
\end{proof}  

\noindent In the particular case of the $\m$-adic filtration on $A, $   we  can   easily give   a partial  extension of the above result without assuming the Cohen-Macaulyness of $A $ (see Theorem 3.2. \cite{RV5}).
  
 \begin{theorem} \label{vaznCM} Let $(A, \m) $ be a local ring of dimension $r$ and let $J$ be an ideal generated by  a maximal superficial sequence for $\m. $ If  $\dim \ S_J(\m) =r, $ then
 $$  e_0(S_J(\m)) \le \sum_{j\ge 0} v_j(\m) -e_0(\m) +1 $$ 
The following facts are equivalent:

\vskip 2mm 
1. $ e_0(S_J(\m)) =  \sum_{j\ge 0} v_j(\m) -e_0(\m) +1 $ 
\vskip 2mm
 
2. $ e_1(\m) -e_1(J) = \sum_{j\ge 0} v_j(\m)  $
\vskip 2mm
 
 3. $A$ is Cohen-Macaulay  and $ \depth gr_{\m}(A) \ge r-1$
 \vskip 2mm
 
 \end{theorem}
 
 \begin{proof} It follows by   Corollary \ref{e0S} and  Theorem \ref{th2}.
 \end{proof}

\end{document}